\documentclass[%
  a4paper,
  onecolumn,
]{mypreprint}

\usepackage[english]{babel}

\usepackage{graphicx}
\usepackage{tikz}
\usepackage{pgfplots}
  \pgfplotsset{compat = 1.13}
  \usepgfplotslibrary{external}
  \tikzset{external/system call = {%
      pdflatex \tikzexternalcheckshellescape
      -halt-on-error
      -interaction=batchmode
      -jobname "\image" "\texsource"}}
  \tikzexternaldisable

\usepackage{filemod}

\usetikzlibrary{calc}

\usepackage{caption}
\usepackage{subcaption}

\usepackage{amsmath}
\usepackage{amssymb}
\usepackage{amsthm}

\usepackage{scalerel}

\usepackage{enumitem}
\usepackage{siunitx}
\usepackage{pgfplotstable}
\usepackage{booktabs}
\usepackage[shortcuts]{extdash}
\usepackage{csquotes}

\usepackage{cleveref}
\usepackage{todonotes}

\newcommand{\mv}[1]{\ensuremath{\boldsymbol{\mathbf{#1}}}}
\newcommand{\re}[1]{\ensuremath{\operatorname{Re}\left(#1\right)}}
\newcommand{\im}[1]{\ensuremath{\operatorname{Im}\left(#1\right)}}

\renewcommand{\rm}[1]{\ensuremath{\mathrm{#1}}}

\renewcommand{\i}{\ensuremath{\mathrm{i}}}

\newcommand{\C}{\ensuremath{\mathbb{C}}}
\newcommand{\N}{\ensuremath{\mathbb{N}}}
\newcommand{\R}{\ensuremath{\mathbb{R}}}

\renewcommand\widehat[1]{\hstretch{2}{\hat{\hstretch{.5}{#1}}}}

\DeclareMathOperator{\mspan}{span}
\newcommand{\trans}{\ensuremath{\mkern-1.5mu\mathsf{T}}}
\newcommand{\herm}{\ensuremath{\mathsf{H}}}

\newcommand{\stK}{\ensuremath{\boldsymbol{\mathcal{K}}}}
\newcommand{\stF}{\ensuremath{\boldsymbol{\mathcal{F}}}}
\newcommand{\stG}{\boldsymbol{\ensuremath{\mathcal{G}}}}
\newcommand{\sys}{\ensuremath{\Sigma}}
\newcommand{\tf}{\ensuremath{\mv{H}}}
\newcommand{\tfs}{\ensuremath{H}}

\newcommand{\stKr}{\ensuremath{\skew1\widehat{\stK}}}
\newcommand{\stFr}{\ensuremath{\skew1\widehat{\stF}}}
\newcommand{\stGr}{\ensuremath{\skew1\widehat{\stG}}}
\newcommand{\rom}{\ensuremath{\widehat{\Sigma}}}
\newcommand{\tfr}{\ensuremath{\widehat{\mv{H}}}}
\newcommand{\Mr}{\ensuremath{\widehat{\mv{M}}}}
\newcommand{\Cr}{\ensuremath{\widehat{\mv{C}}}}
\newcommand{\Kr}{\ensuremath{\widehat{\mv{K}}}}
\newcommand{\Fr}{\ensuremath{\widehat{\mv{F}}}}
\newcommand{\Gr}{\ensuremath{\widehat{\mv{G}}}}
\newcommand{\xr}{\ensuremath{\hat{\mv{x}}}}
\newcommand{\yr}{\ensuremath{\hat{\mv{y}}}}

\newcommand{\Hinf}{\ensuremath{\mathcal{H}_{\infty}}}
\newcommand{\Linf}{\ensuremath{\mathcal{L}_{\infty}}}
\newcommand{\Htwo}{\ensuremath{\mathcal{H}_{2}}}

\newcommand{\twonorm}[1]{\ensuremath{\left\lVert#1\right\rVert_2}}

\newcommand{\hinfnorm}[1]{\ensuremath{\left\lVert#1\right\rVert_{\Hinf}}}
\newcommand{\linfnorm}[1]{\ensuremath{\left\lVert#1\right\rVert_{\Linf}}}

\newcommand{\hinfrelerr}{\hinfnorm{\tf-\tfr} / \hinfnorm{\tf}}

\newcommand{\morscore}{\mbox{MORscore}}
\newcommand{\morscores}{\mbox{MORscores}}

\newcommand{\ubraces}[1]{\ensuremath{\left[#1\right]}}
\newcommand{\errlabel}{Relative $\Linf$-error}

\newcommand{\soP}{\textit{p}}
\newcommand{\soPM}{\textit{pm}}
\newcommand{\soPV}{\textit{pv}}
\newcommand{\soVP}{\textit{vp}}
\newcommand{\soVPM}{\textit{vpm}}
\newcommand{\soV}{\textit{v}}
\newcommand{\soFV}{\textit{fv}}
\newcommand{\soSO}{\textit{so}}

\theoremstyle{plain}\newtheorem{proposition}{Proposition}

\newcommand{\moravg}{\textsf{avg}}
\newcommand{\morequi}{\textsf{equi}}
\newcommand{\morlinf}{$\mathsf{L}_{\infty}$}
\newcommand{\morminrel}{\textsf{minrel}}
\newcommand{\morbt}{\textsf{SOBT}}
\newcommand{\morsp}{\textsf{sp}}
\newcommand{\moraaa}{\textsf{soa}}
\newcommand{\tsimag}{\textsf{tsimag}}
\newcommand{\tsreal}{\textsf{tsreal}}
\newcommand{\osimaginput}{\textsf{osimaginput}}
\newcommand{\osimagoutput}{\textsf{osimagoutput}}
\newcommand{\osrealinput}{\textsf{osrealinput}}
\newcommand{\osrealoutput}{\textsf{osrealoutput}}

\newcommand{\lmoravg}{\textsf{avg}}
\newcommand{\lmorequi}{\textsf{equi}}
\newcommand{\lmorlinf}{$\mathsf{L}_{\infty}$}
\newcommand{\lmorminrel}{\textsf{minrel}}
\newcommand{\lmorbt}{\textsf{SOBT}}
\newcommand{\lmorsp}{\textsf{sp}}
\newcommand{\lmoraaa}{\textsf{soa}}
\newcommand{\ltsimag}{\textsf{tsimag}}
\newcommand{\ltsreal}{\textsf{tsreal}}
\newcommand{\losimaginput}{\textsf{osimaginput}}
\newcommand{\losimagoutput}{\textsf{osimagoutput}}
\newcommand{\losrealinput}{\textsf{osrealinput}}
\newcommand{\losrealoutput}{\textsf{osrealoutput}}

\definecolor{TUMBlau}{RGB}{0,101,189}
\definecolor{TUMBlauDunkel}{RGB}{0,82,147}
\definecolor{TUMBlauHell}{RGB}{152,198,234}
\definecolor{TUMBlauMittel}{RGB}{100,160,200}
\definecolor{TUMElfenbein}{RGB}{218,215,203}
\definecolor{TUMGruen}{RGB}{162,173,0}
\definecolor{TUMOrange}{RGB}{227,114,34}
\definecolor{TUMDunkelGrau}{gray}{0.8}
\definecolor{TUMGrau}{gray}{0.5}
\definecolor{TUMHellGrau}{gray}{0.2}
\definecolor{TUMLila}{RGB}{105,008,090}
\definecolor{TUMRed}{RGB}{156,013,022}

\newcommand{\mylinewidth}{1.25pt}
\newcommand{\mymarkscale}{1.5}

\pgfplotscreateplotcyclelist{TUMcolorlist}{%
	TUMBlau,mark=*, line width=\mylinewidth, mark options={scale=\mymarkscale,fill=white,solid}\\%
	TUMOrange,mark=triangle*,line width=\mylinewidth,mark options={scale=\mymarkscale,fill=white,solid}\\%
	TUMGruen,mark=square*,line width=\mylinewidth,mark options={scale=\mymarkscale,fill=white,solid}\\%
	TUMLila,mark=otimes*,line width=\mylinewidth,mark options={scale=\mymarkscale,fill=white,solid}\\%
	TUMRed,mark=diamond*,line width=\mylinewidth,mark options={scale=\mymarkscale,fill=white,solid}\\%
	TUMGrau,mark=pentagon*,line width=\mylinewidth,mark options={scale=\mymarkscale,fill=white,solid}\\%
	TUMBlau,mark=otimes*,line width=\mylinewidth,mark options={scale=\mymarkscale,fill=white,solid}\\%
	TUMOrange,mark=diamond*,line width=\mylinewidth,mark options={scale=\mymarkscale,fill=white,solid}\\%
	TUMGruen,mark=pentagon*,line width=\mylinewidth,mark options={scale=\mymarkscale,fill=white,solid}\\%
	TUMLila,mark=*,line width=\mylinewidth,mark options={scale=\mymarkscale,fill=white,solid}\\%
	TUMRed,mark=triangle*,line width=\mylinewidth,mark options={scale=\mymarkscale,fill=white,solid}\\%
	TUMGrau,mark=square*,line width=\mylinewidth,mark options={scale=\mymarkscale,fill=white,solid}\\%
}

\pgfplotscreateplotcyclelist{TUMcolorlist_simple}{%
	TUMBlau\\%
	TUMOrange\\%
	TUMGruen\\%
	TUMLila\\%
	TUMRed\\%
	TUMGrau\\%
}

\begin{document}


\title{Structured model order reduction for vibro-acoustic problems
  using interpolation and balancing methods}

\author[$\ast$]{Quirin Aumann}
\affil[$\ast$]{Max Planck Institute for Dynamics of Complex Technical Systems,
  Sandtorstra{\ss}e 1, 39106 Magdeburg, Germany.\authorcr
  \email{aumann@mpi-magdeburg.mpg.de }, \orcid{0000-0001-7942-5703}}

\author[$\dagger$]{Steffen W. R. Werner}
\affil[$\dagger$]{Courant Institute of Mathematical Sciences,
  New York University, New York, NY 10012, USA.
  \email{steffen.werner@nyu.edu}, \orcid{0000-0003-1667-4862}}

\shorttitle{Model order reduction for vibro-acoustic problems}
\shortauthor{Q. Aumann, S. W. R. Werner}
\shortdate{2022-01-17}
\shortinstitute{}

\keywords{vibro-acoustic system, second-order system, model order reduction,
  structured interpolation, structure-preserving balanced truncation}

\msc{}

\abstract{%
	Vibration and dissipation in vibro-acoustic systems can be assessed using 
  frequency response analysis.
	Evaluating a frequency sweep on a full-order model can be very 
  costly, so model order reduction methods are employed to compute
  cheap-to-evaluate surrogates.
	This work compares structure-pre\-ser\-ving model reduction methods based on
  rational interpolation and balanced truncation with a specific focus on their 
  applicability to vibro-acoustic systems.
	Such models typically exhibit a second-order structure and their material 
  properties as well as their excitation may be depending on the driving
  frequency.
	We show and compare the effectiveness of all considered methods by applying 
  them to numerical models of vibro-acoustic systems depicting structural 
  vibration, sound transmission, acoustic scattering, and poroelastic problems.
}

\novelty{%
  We give overviews about common system types arising from modeling of
  vibro-acoustic structures as well as state-of-the-art structure-preserving
  model reduction approaches based on interpolation and balanced truncation.
  The presented reduction methods are then employed and compared in
  efficiency on benchmark examples of the different vibro-acoustic
  structures.
}

\maketitle


\section{Introduction}%
\label{sec:intro}

The numerical simulation of structural vibration, acoustic wave propagation, and 
their combination, often termed \emph{vibro-acoustic problems}, is an important 
tool for many engineering applications.
Especially, the prevention of unwanted noise or vibration is important in
practice and many methods to dissipate surplus vibration energy have been 
established~\cite{Mar02, Dud08, ClaVSetal13, XieZJetal19}.
These problems are typically evaluated by frequency sweeps.
In the Laplace (frequency) domain, a vibro-acoustic system is described
by linear systems of equations of the form
\begin{align} \label{eqn:sys}
  \sys\colon \begin{cases}
    \begin{aligned}
      \left( s^{2} \mv{M}(s) + s \mv{C}(s) + \mv{K}(s) \right) \mv{x}(s) & = 
      \mv{F}(s) \mv{u}(s),\\
      \mv{y}(s) & = \mv{G} \mv{x}(s),
    \end{aligned}
    \end{cases}
\end{align}
with frequency-dependent matrix-valued functions $\mv{M}, \mv{C}, \mv{K} \colon 
\C \rightarrow \C^{n \times n}$, describing the internal dynamics and 
representing mass, damping and stiffness, respectively, $\mv{F}\colon \C 
\rightarrow \C^{n \times m}$, for the external forcing via the input 
$\mv{u}(s)$, and the constant matrix $\mv{G} \in \R^{p \times n}$, describing 
the quantities of interest as linear combinations of the system states.
Therein, the controls $\mv{u}\colon \C \rightarrow \C^{m}$ are used to steer the
system state $\mv{x}\colon \C \rightarrow \C^{n}$ to obtain the desired behavior
of the outputs $\mv{y}\colon \C \rightarrow \C^{p}$.
Under the assumption that the system~\cref{eqn:sys} is regular, i.e., there
exists an $s \in \C$ for which the frequency-dependent functions can be 
evaluated and the center term for the (linear) dynamics, $\left( s^{2} 
\mv{M}(s) + s \mv{C}(s) + \mv{K}(s) \right)$, is invertible, the 
input-to-output behavior of~\cref{eqn:sys} is given by its transfer function
\begin{align} \label{eqn:tf}
  \tf(s) & = \mv{G} \left( s^{2} \mv{M}(s) + s \mv{C}(s) + \mv{K}(s) 
  \right)^{-1} \mv{F}(s).
\end{align}

In practical applications, there is a demand for highly accurate models.
Consequently, the number of equations in~\cref{eqn:sys} quickly becomes very
large ($n \gtrapprox 10^{5}$), which makes evaluation of~\cref{eqn:sys} rather
expensive in terms of computational resources such as time and memory.
Model order reduction is a remedy to construct cheap-to-evaluate surrogate 
systems.
In this paper, we consider structure-preserving model order reduction, 
i.e., the computed reduced-order model retains the original 
physically-inspired system structure of the full-order model or at least an 
interpretable equivalent.
This has been shown to yield more accurate approximations if the reduced 
system's order is low or allows even a physical re-interpreation of the 
reduced-order quantities.
In case of vibro-acoustic systems~\cref{eqn:sys}, the reduced-order
model should have the form
\begin{align} \label{eqn:rom}
  \rom\colon \begin{cases}
    \begin{aligned}
      \left( s^{2} \Mr(s) + s \Cr(s) + \Kr(s) \right) \xr(s) & = \Fr(s)
        \mv{u}(s),\\
      \yr(s) & = \Gr \xr(s),
    \end{aligned}
    \end{cases}
\end{align}
with $\Mr, \Cr, \Kr \colon \C \rightarrow \C^{r \times r}$, $\Fr\colon \C
\rightarrow \C^{r \times m}$ and $\Gr \in \C^{p \times r}$, and a much smaller
number of internal states and defining equations $r \ll n$.
Its corresponding (reduced-order) transfer function is given by
\begin{align} \label{eqn:tfr}
  \tfr(s) & = \Gr \left( s^{2} \Mr(s) + s \Cr(s) + \Kr(s) \right)^{-1} \Fr(s).
\end{align}
To serve as surrogate, the reduced-order system~\cref{eqn:rom} needs to
approximate the original system's input-to-output behavior, i.e., for the same
input, the outputs of~\cref{eqn:sys,eqn:rom} should be close to each other:
\begin{align*}
  \left\lVert \mv{y} - \yr \right\rVert \leq \tau \cdot 
    \left\lVert \mv{u} \right\rVert,
\end{align*}
with a tolerance $\tau$, in some appropriate norm and for all admissible
inputs $\mv{u}$.
For frequency domain models, as they occur in vibro-acoustics, the approximation
of the input-to-output behavior amounts to the approximation of the system's
transfer function in frequency regions of interest such that
\begin{align*}
  \left\lVert \tf - \tfr \right\rVert \leq \tau,
\end{align*}
with the transfer functions from~\cref{eqn:tf,eqn:tfr}.

The majority of system-theoretic model reduction techniques for dynamical
systems can be classified as either based on system modes, interpolation or
balancing of energy.
Classically, modal approaches are used to solve structural and acoustic problems 
efficiently.
These methods consider the solution of the system's underlying eigenproblem 
to find poles of interest, which are retained in the reduced-order system.
Modal methods are the classical, well-established approaches for the efficient
reduction of vibro-acoustic systems and, therefor, we will not consider them
further in this work.
See~\cite{BesTLetal13} for a review on modal techniques for
vibro-acoustic systems, \cite{RouDL17}~for the special case of
systems with frequency-dependent material properties, and~\cite{Wer21} for a new
structured dominant pole-based approach for systems with modal damping.

Interpolation or moment-matching methods aim for low-order approximations that
match the original transfer function behavior at certain expansion points.
A general framework for structure-preserving interpolation of linear systems
has been proposed in~\cite{BeaG09}.
This approach can be immediately used for structured systems of the
form~\cref{eqn:sys}, and has been applied to vibro-acoustic systems
in~\cite{HetTF12}.
For efficient numerical computations, ideas for Arnoldi-like algorithms
have been extended to the second-order system
case~\cite{BaiS05, LuSB16} for application to structural 
and vibro-acoustic systems~\cite{BaiS05a, VanADetal17, XieZJetal20}.
Further related methods for vibro-acoustic systems with poroelastic materials
can be found in, e.g.,~\cite{DecDMetal20,RumGD14, Rum18}.
An important aspect in interpolatory model order reduction is the choice of 
appropriate interpolation points.
The \emph{iterative rational Krylov algorithm (IRKA)}~\cite{GugAB08} 
computes interpolation point locations automatically.
The resulting first-order systems are locally optimal in the sense of the
approximation error's $\Htwo$-norm.
Heuristic extensions for second-order systems with constant coefficient matrices 
have been proposed in, e.g.,~\cite{Wya12}, and proven to provide good 
results in practice.

One of the most successful model reduction approaches for unstructured
first-order systems is the balanced truncation method~\cite{Moo81}.
This approach considers the energy behavior of the system to identify parts of
the state which contribute only marginally to the input-to-output behavior.
Different extensions of heuristic nature have been proposed for second-order
systems~\cite{MeyS96, ReiS08, ChaLVetal06}, which yield good results
in practical applications~\cite{SaaSW19}.
More recently, new formulas have been proposed to limit the approximation
behavior of the second-order balanced truncation methods in frequency or time
domain to ranges of interest~\cite{BenW21b}.
Extensive comparisons can be found in~\cite{Wer21}.

In this work, we present and compare different methods for reducing the
numerical complexity of vibro-acoustic systems.
Three main types of systems, which are encountered in the vibro-acoustic
setting, are grouped regarding their damping and coupling behavior, and
numerical benchmarks are given for each type.
This allows a structured assessment of the applicability and the approximation
quality of the employed methods.
Here, we only consider interpolatory and balancing-related methods for
second-order systems and do not include modal methods into our comparisons.
Reviews and comparisons for such methods are given
in~\cite{BesTLetal13, SaaSW19, RouDL17, Wer21}.
We identify, which reduction algorithms are applicable to a wide range of
different systems without the need of problem specific modifications.
The approximation quality is assessed by computing the MORscore~\cite{Him21},
which yields a single number for each method on each benchmark example.

The rest of the article is organized as follows:
In \Cref{sec:vibro}, we give an overview about numerical modeling of
vibro-acoustic systems and group different application cases regarding the form
of their transfer functions.
The following \Cref{sec:mor} reviews the concepts of model order reduction
methods we will use for the comparison.
In \Cref{sec:numerics}, numerical models for all proposed types of
vibro-acoustic systems are introduced and the performance of applicable
reduction methods is compared.
\Cref{sec:conclusions} concludes the article.

\section{Types of vibro-acoustic systems}%
\label{sec:vibro}

In this section, we describe the different system types occurring in the
modeling of vibro-acoustic problems that we consider in this paper.


\subsection{Structural dynamics}

Structural vibration in spatial domains is often modeled using the
\emph{finite element method (FEM)}.
A discretization with finite elements leads to a system with the following
transfer function
\begin{align} \label{eqn:sotf}
	\tf(s) & = \mv{G} \left( s^{2} \mv{M} + s \mv{C} + \mv{K} \right)^{-1} 
	\mv{F},
\end{align}
where the typically symmetric matrices $\mv{M}$ and $\mv{K}$ resemble the mass 
and stiffness of the model, respectively.
The viscous damping, introduced by the damping matrix $\mv{C}$, is often
proportionally described, e.g., by Rayleigh damping.
In this case, the corresponding damping matrix is a linear combination of the 
system's mass and stiffness matrices and is given by
\begin{align*}
	\mv{C}_{\rm{R}} = \alpha \mv{M} + \beta \mv{K},
\end{align*}
where the coefficients $\alpha$ and $\beta$ control the frequency region in
which the structure is damped.
Therein, $\alpha$ models the damping effect that the surrounding medium has on
the structure, while $\beta$ accounts for the material damping.

The effect of proportional damping is not constant over the complete frequency
range and the coefficients have to be tuned individually for each problem.
A damping effect being constant for all frequencies can be described by
structural damping, often introduced by a complex stiffness matrix.
The corresponding frequency-dependent damping matrix is given by
\begin{align*}
	\mv{C}_{\rm{S}}(s) = \frac{\eta}{s} \i \mv{K},
\end{align*}
where $\eta$ is the structural loss factor, for which standard values are
available for various materials.
Discrete damping elements, for example dashpot dampers, can also be introduced
in $\mv{C}$.
This results in a matrix structure, which is not proportional to mass or
stiffness.


\subsection{Acoustical modeling}

Acoustic wave propagation is often modeled via the Helmholtz equation, which can
also be discretized by finite elements.
For bounded problems, for example inside a cavity, this leads to a transfer
function of the form
\begin{align}\label{eqn:vap_d}
	\tf(s) & = \mv{G} \left( \frac{s^{2}}{c^2} \mv{M} + s \mv{C} + \mv{K} 
	\right)^{-1} \mv{F}(s),
\end{align}
with $c$, the wave speed in the acoustic medium.
The acoustic mass matrix $\mv{M}$ represents the compressibility of the medium 
and the acoustic stiffness matrix $\mv{K}$ its mobility.
Damping is introduced, for example, by admittance boundary conditions such that
the acoustic damping matrix $\mv{C}$ is not proportional to $\mv{M}$ or 
$\mv{K}$.
The three material matrices are typically symmetric.
The frequency-dependent acoustic source term $\mv{F}(s)$ introduces
velocity or pressure sources into the system.
Its frequency dependency is either linear or quadratic.

Unbounded problems depicting acoustic wave propagation in the open space can be
modeled using different methods, such as absorbing boundary conditions (ABC),
perfectly matched layers (PML), or infinite elements
(IE)~\cite{Ihl98}.
In these methods, the model is truncated at an arbitrarily chosen boundary,
which is then treated in a way that Sommerfeld's radiation condition is
fulfilled.
An example for an ABC is the Dirichlet-to-Neumann (DtN) condition, which imposes
the analytical solution of the exterior domain on the arbitrary model boundary,
leading to densely populated matrices~\cite{Giv99}.
IEs map the solution of decaying waves traveling outwards of the modeled domain
to a finite set of nodes on which the numerical integration can be
performed~\cite{Ast08}.
PMLs are absorbing layers added to the boundaries of the modeled domain. 
These layers ensure that waves can enter this region without
being reflected at the boundary and are decaying inside the layer so they do not
travel back into the domain of interest~\cite{BerM20}.
Introducing these conditions leads to a system with a damping matrix depending 
on the driving frequency of the system and transfer function
\begin{align}
	\tf(s) & = \mv{G} \left( \frac{s^{2}}{c^2} \mv{M} + \mv{C}\left(s\right) + 
	\mv{K} \right)^{-1} \mv{F}(s).
\end{align}
The matrices resulting from a discretization with finite elements are complex
symmetric and may for special cases be rewritten such that the frequency 
dependency can be separated.
In case of PML, for example, it can be tuned for a specific frequency, such 
that the frequency dependency vanishes~\cite{VerMR19}.


\subsection{Vibro-acoustic systems}

Coupling a vibrating structure with the surrounding acoustic fluid results in a
vibro-acoustic system.
The coupling is active in both ways: The vibrating structure radiates energy
into the adjacent medium and is, in turn, excited by waves traveling through the
acoustic fluid.
Such systems are often used for modeling the sources of noise in machines or
vehicles and to find ways to dissipate unwanted vibration or acoustic energy 
by introducing suitable damping.
Examples for these mechanisms include poroelastic materials~\cite{Bio56},
constrained-layer damping~\cite{AmiA09}, or mechanical
joints~\cite{MatBKetal20}.
The dissipation mechanisms are often governed by frequency-dependent
complex-valued functions $\tilde{\phi}_{i}(s)$.
Their influence on the model is given by corresponding constant matrices
$\mv{C}_{i}$ resulting in transfer functions of the form
\begin{align*}
	\tf(s) & = \mv{G} \left( s^{2} \mv{M} + s \mv{C} + \mv{K} + 
	  \sum\limits_{i = 1}^{k} \tilde{\phi}_{i}(s) \mv{C}_{i}  \right)^{-1}
	  \mv{F}(s),
\end{align*}
for some $k \in \N$.
The coupling between the solid and fluid phases introduces off-diagonal terms 
in $\mv{M}$ and $\mv{K}$, making them non-symmetric.
The frequency dependency of the input is only given, if the system is excited by
an acoustic source.
Structural excitation can be modeled by a constant input $\mv{F}$.


\subsection{Model problems}%
\label{subsec:vibro_categories}

Although most vibro-acoustic systems are depicted as second-order systems, the
different damping and coupling methods described above influence the structure
of the transfer function, which may hinder the application of a specific
reduction method.
Therefor, we identify three types of model problems with varying properties,
which will be reduced using applicable model order reduction techniques in the
following:
\begin{description}
  \item[Case A]
    A structural or interior vibro-acoustic system with proportional damping and
    no acoustic source following
    \begin{align*}
  	  \tf(s) & = \mv{G} \left( s^{2} \mv{M} + s \left(d_{1} \mv{M} + d_{2} 
  	    \mv{K} \right) + \mv{K} \right)^{-1} \mv{F},
    \end{align*}
  	where the damping factors $d_{1}$ and $d_{2}$ represent either Rayleigh
    or hysteretic damping.
    The resulting system matrices may be complex-valued if hysteretic damping is
    applied and non-symmetric for an interior vibro-acoustic problem.
  \item[Case B]
    An interior acoustic or vibro-acoustic system with acoustic source following
    \begin{align*}
      \tf(s) & = \mv{G} \left( \frac{s^{2}}{c^2} \mv{M} + s \mv{C} + \mv{K} 
        \right)^{-1} \mv{F}(s).
    \end{align*}
    An exterior radiation problem can also be modeled with this system type, if,
    for example, a PML is tuned to a single frequency.
  \item[Case C]
    An interior vibro-acoustic system with acoustic source and
    frequency-dependent material properties following
    \begin{align} \label{eqn:tfcasec}
      \tf(s) & = \mv{G} \left( s^{2} \mv{M} + s \mv{C} + \mv{K} + 
        \sum\limits_{i = 1}^{k} \phi_{i}(s) \mv{C}_{i} \right)^{-1} \mv{F}(s),
    \end{align}
    where the frequency dependency can be described in an affine representation.
    The parameter dependency on the wave velocity $c$ is included in the entries
    of $\mv{M}$ corresponding to the acoustic fluid.
    Again, exterior problems can be modeled, if the method ensuring Sommerfeld's
    radiation condition can be represented by either a constant matrix or linear
    combination of matrices and corresponding frequency-dependent functions.
\end{description}

\section{Model order reduction for vibro-acoustic systems}%
\label{sec:mor}

In this section, model order reduction methods for vibro-acoustic problems are
considered.
We begin with the general concept of projection-based model reduction and outline
afterwards methods based on interpolation/moment matching and balancing.


\subsection{Model reduction by projection}%
\label{subsec:proj}

Consider linear systems of the form~\cref{eqn:sys}, which are described by
their transfer functions~\cref{eqn:tf}.
The goal of model reduction methods is now the construction of~\cref{eqn:rom}
such that the input-to-output behavior of~\cref{eqn:sys} is approximated.
The resulting main question is how to construct~\cref{eqn:rom}.
In projection-based model reduction, two (constant) truncation matrices 
$\mv{V}, \mv{W} \in \C^{n \times r}$ are constructed as bases of underlying 
projection spaces $\mspan(\mv{V})$ and $\mspan(\mv{W})$ such that the 
reduced-order quantities can be computed by
\begin{align} \label{eqn:proj}
  \begin{aligned}
      \Mr(s) & = \mv{W}^{\herm} \mv{M}(s) \mv{V}, &
      \Cr(s) & = \mv{W}^{\herm} \mv{C}(s) \mv{V}, &
      \Kr(s) & = \mv{W}^{\herm} \mv{K}(s) \mv{V},\\
      \Fr(s) & = \mv{W}^{\herm} \mv{F}(s), &
      \Gr & = \mv{G} \mv{V};
  \end{aligned}
\end{align}
see, e.g.,~\cite{AntBG20, BeaG09}.
In practice, the matrix functions in~\cref{eqn:proj} have particular
realizations from which the reduced-order models are built, e.g., the matrix
function describing the mass of the system can always be realized in a
frequency-affine decomposition
\begin{align} \label{eqn:freqaff}
  \mv{M}(s) & = \sum\limits_{k = 1}^{n_{\rm{M}}} g_{k, \rm{M}}(s) \mv{M}_{k},
\end{align}
with constant matrices $\mv{M}_{k} \in \C^{n \times n}$, scalar 
frequency-dependent functions $g_{k, \rm{M}}\colon \C \rightarrow \C$ and 
usually $n_{\rm{M}} \ll n$.
The reduced-order matrix function is then given by
\begin{align} \label{eqn:romfreqaff}
  \Mr(s) & = \mv{W}^{\herm} \mv{M}(s) \mv{V} = \sum\limits_{k = 1}^{n_{\rm{M}}}
    g_{k,\rm{M}}(s) \mv{W}^{\herm} \mv{M}_{k} \mv{V} = 
    \sum\limits_{k = 1}^{n_{\rm{M}}} g_{k,\rm{M}}(s) \Mr_{k},
\end{align}
with reduced-order matrices $\Mr \in \C^{r \times r}$.
In a similar way, the reduced matrix functions for the other terms
in~\cref{eqn:proj} are given.
Since the scalar frequency-dependent functions
in~\cref{eqn:freqaff,eqn:romfreqaff} are the same, model order reduction by
projection preserves the internal system structure in reduced-order models.
Therefor, the reduced-order matrices in~\cref{eqn:romfreqaff} can be used to
replace their full-order counterparts to give a realization of the reduced-order
model.

As a particular example, consider the classical second-order
system~\cref{eqn:sotf}.
The reduced-order system is then given by the matrices
\begin{align*}
  \begin{aligned}
    \Mr & = \mv{W}^{\herm} \mv{M} \mv{V}, &
      \Cr & = \mv{W}^{\herm} \mv{C} \mv{V}, &
      \Kr & = \mv{W}^{\herm} \mv{K} \mv{V}, &
      \Fr & = \mv{W}^{\herm} \mv{F}, &
      \Gr & = \mv{G} \mv{V}.
  \end{aligned}
\end{align*}
However, with these specific construction rules for reduced-order models, the
main question in projection-based model reduction is relocated to the actual
construction of $\mv{V}$ and $\mv{W}$.
The following sections give a short overview about the model reduction methods
that will be used in the numerical experiments of this paper.


\subsection{Moment matching and interpolation}%
\label{subsec:interp}

A classical approach for choosing $\mv{V}$ and $\mv{W}$ as of \Cref{subsec:proj}
is by moment matching (interpolation) of the system's transfer function.
With the observation that classical linear systems, e.g.,
second-order systems of the form~\cref{eqn:sotf}, have rational transfer
functions, the idea of moment matching roots in the theory of Pad{\'e}
approximation~\cite{Bak75, BulV86}.
Thereby, rational Hermite interpolants of minimal degree in numerator and
denominator are constructed.
On the other hand, the moment matching method considers the representation of
the transfer function in terms of a power series, for which the first
coefficients are matched; see, e.g.,~\cite{GraL83, DeVS87}.
In a more general setting, all these ideas amount to the construction
of the reduced-order model such that its corresponding transfer function solves
a Hermite interpolation problem of the form
\begin{align} \label{eqn:interp}
  \frac{\rm{d}^{i_{j}}}{\rm{d}s^{i_{j}}}\tf(\sigma_{j}) & =
    \frac{\rm{d}^{i_{j}}}{\rm{d}s^{i_{j}}}\tfr(\sigma_{j}),
\end{align}
for $j = 0, \ldots, k$ and $0 \leq i_{j} \leq \ell_{j}$, in the interpolation
points $\sigma_{1}, \ldots, \sigma_{k} \in \C$.
The most important result, in the computational sense, is that this
interpolation~\cref{eqn:interp} can be performed using the projection framework 
from \Cref{subsec:proj}; see, e.g.,~\cite{Gri97, BeaG09}; which makes the
moment matching approach well suited for use in many different applications.
In the following, a quick overview about particular moment matching methods
related to the system structure occurring in acoustic problems
(\Cref{sec:vibro}) is given.
For a more general introduction to interpolation and moment matching see
also~\cite{AntBG20}.


\subsubsection{Second-order Arnoldi method}%
\label{subsec:arnoldi}

Krylov subspaces are often employed as valid choices for the projection
subspaces $\mspan(\mv{V})$ and $\mspan(\mv{W})$.
Their orthogonal bases $\mv{V}$ and $\mv{W}$ can, for example, be computed via 
an Arnoldi method~\cite{SuC91, BaiS05a, LuSB16}.
For simplicity, we consider in this theoretic overview only systems with a
single input and a single output (SISO), i.e., $\mv{F} = \mv{f} \in \C^{n}$
and $\mv{G} = \mv{g} \in \C^{n}$.
For the application of the method to multi-input/multi-output (MIMO) systems
see, e.g.,~\cite{ChuTL10}.
The generalized $r$-th \emph{second-order Krylov subspace} is defined by two
matrices $\mv{A}, \mv{B} \in \C^{n \times n}$ and a vector $\mv{v}_{0} \in 
\C^{n}$ such that
\begin{align*}
	\mathcal{S}_{r} \left(\mv{A},\mv{B}; \mv{v}_{0}\right) & = 
	  \mspan\left( \mv{v}_{0}, \mv{v}_{1}, \dots, \mv{v}_{r-1} \right),
\end{align*}
where the recursively related vectors $\mv{v}_{k}$ are given by
\begin{align} \label{eqn:ks}
	\begin{aligned}
		\mv{v}_{1} &= \mv{A} \mv{v}_{0},\\
		  \mv{v}_{i} &= \mv{A} \mv{v}_{i-1} + \mv{B} \mv{v}_{i-2}.
	\end{aligned}
\end{align}
The vectors in~\cref{eqn:ks} are also known as the \emph{Krylov sequence} based
on $\mv{A}, \mv{B}, \mv{v}_{0}$ and form the sought truncation matrix
$\mv{V} = \left[\mv{v}_0, \mv{v}_1, \dots, \mv{v}_{r-1}\right]$.
In order to find a basis, which projects the original system~\cref{eqn:sotf}
onto a low-dimensional subspace such that the reduced system matches the
first $r$ moments of the original system, the matrices and starting vector
defining the Krylov sequence are set according to $\mv{A} = -\mv{K}^{-1}\mv{C}$,
$\mv{B} = -\mv{K}^{-1} \mv{M}$ and $\mv{v}_{0} = \mv{K}^{-1} \mv{f}$.
The truncation matrix $\mv{W}$ can be computed in the same way considering the
adjoint problem.

If approximation around a specified frequency $s_{0}$ other than zero is
desired, the transfer function can be shifted about this shift $s_{0}$ as in
\begin{align*}
	\tfs(s) & = \mv{g} \left( (s - s_{0})^{2} \mv{M} + (s - s_{0}) 
	  \widetilde{\mv{C}} + \widetilde{\mv{K}} \right)^{-1} \mv{f},
\end{align*}
with $\widetilde{\mv{C}} = 2 s_{0} \mv{M} + \mv{C}$ and $\widetilde{\mv{K}} = 
s_{0}^{2} \mv{M} + s_{0} \mv{C} + \mv{K}$.
The corresponding $r$-th second-order Krylov subspace is then given by
\begin{align*}
	\mathcal{S}_{r} \left( -\widetilde{\mv{K}}^{-1} \widetilde{\mv{C}},
    -\widetilde{\mv{K}}^{-1} \mv{M}; \widetilde{\mv{K}}^{-1} \mv{f} \right).
\end{align*}

The choice of the subspace dimension $r$ and location of the shift has a large
influence on the approximation quality of the reduced model.
To obtain a reduced-order model with a high accuracy for a wide range of
frequencies, it usually is beneficial to not only increase the size of the
Krylov subspace to match higher-order moments, but to combine multiple subspaces
with different shifts into a global basis.


\subsubsection{Structure-preserving interpolation}
\label{subsubsec:structure_preserving_interpolation}

A more general approach for moment matching of structured linear systems is
described in~\cite{BeaG09}.
Any matrix-valued function of the form
\begin{align} \label{eqn:sttf}
  \tf(s) & = \stG(s) \stK(s)^{-1} \stF(s),
\end{align}
with $\stG\colon \C \rightarrow \C^{p \times n}$, representing the outputs,
$\stK\colon \C \rightarrow \C^{n \times n}$, for the internal dynamics, and
$\stF\colon \C \rightarrow \C^{n \times m}$, describing the system inputs,
can be structure-preserving interpolated by projection.
Given the two basis matrices $\mv{V}, \mv{W} \in \C^{n \times r}$ of underlying 
right and left projection spaces, the reduced-order model is computed similarly
to~\cref{eqn:proj} by
\begin{align} \label{eqn:stproj}
  \begin{aligned}
    \stGr(s) & = \stG(s) \mv{V}, &
      \stKr(s) & = \mv{W}^{\herm} \stK(s) \mv{V}, &
      \stFr(s) & = \mv{W}^{\herm} \stF(s).
  \end{aligned}
\end{align}
The case of acoustic systems~\cref{eqn:sys} is a special instance
of~\cref{eqn:sttf}, and particular realizations of~\cref{eqn:stproj} are
computed using the same idea of frequency-affine decompositions as
in~\cref{eqn:freqaff,eqn:romfreqaff}.
With the same arguments, projection methods based on~\cref{eqn:stproj} are
guaranteed to preserve the internal system structure.
To give a concise overview, the following proposition states the most important
result from~\cite{BeaG09} to solve~\cref{eqn:interp} for systems of the
form~\cref{eqn:sttf}.

\begin{proposition}[Structured {interpolation~\cite[Theorem~1]{BeaG09}}]%
  \label{prp:strint}
  Let $\tf$ be the transfer function of a linear system, described
  by~\cref{eqn:sttf}, and $\tfr$ the reduced-order transfer function constructed
  by projection~\cref{eqn:stproj}.
  Let the matrix functions $\stG$, $\stK^{-1}$, $\stF$ and $\stKr^{-1}$
  be analytic in the interpolation point $\sigma \in \C$, and let
  $k, \theta \in \N_{0}$ be derivative orders.
  \begin{enumerate}[label = (\alph*)]
    \item If $\mspan \big( \frac{\rm{d}^{j}}{\rm{d}s^{j}}
      (\stK^{-1} \stF) (\sigma) \big) \subseteq \mspan(\mv{V})$, for
      $j = 0, \ldots, k$, then it holds
      \begin{align*}
        \frac{\rm{d}^{j}}{\rm{d}s^{j}} \tf(\sigma)
          & = \frac{\rm{d}^{j}}{\rm{d}s^{j}} \tfr(\sigma),
          & \text{for}~j = 0, \ldots, k.
      \end{align*}
    \item If $\mspan \big( \frac{\rm{d}^{i}}{\rm{d}s^{i}}
      (\stK^{-\herm} \stG^{\herm}) (\sigma) \big) \subseteq \mspan(\mv{W})$, for
      $i = 0, \ldots, \theta$, then it holds
      \begin{align*}
        \frac{\rm{d}^{i}}{\rm{d}s^{i}} \tf(\sigma)
          & = \frac{\rm{d}^{i}}{\rm{d}s^{i}} \tfr(\sigma),
          & \text{for}~i = 0, \ldots, \theta.
      \end{align*}
    \item If $\mv{V}$ and $\mv{W}$ are constructed as in Parts~(a) and~(b),
      then, additionally, it holds
      \begin{align*}
        \frac{\rm{d}^{j}}{\rm{d}s^{j}} \tf(\sigma)
          & = \frac{\rm{d}^{j}}{\rm{d}s^{j}} \tfr(\sigma),
          & \text{for}~j = 0, \ldots, k + \theta + 1.
      \end{align*}
  \end{enumerate}
\end{proposition}

\Cref{prp:strint} shows that only linear systems of equations need to be solved
to construct an interpolating structure-preserving reduced-order model.
However, the usual question that remains for interpolation methods is the choice
of interpolation points.
Over time, there have been many different attempts for heuristics strategies
for how to choose interpolation points.
The following ones will be used in this paper.

A very classical choice are points on the frequency axis $\i \R$.
Depending on the frequency range of interest that is considered, the points are
often chosen either linearly or logarithmically equidistant.
This idea serves usually an overall reasonable approximation behavior but easily
misses features of the system, which are not close enough to the interpolation
points.
In engineering sciences, this is then supplemented by educated guesses of points
in frequency regions, which may be of certain importance due to additional
knowledge about the modeling of the system.

As a more sophisticated choice of interpolation points, methods have been
developed to minimize the approximation error in different systems norms.
The worst case approximation error is described by the $\Hinf$-norm.
Therefor, large-scale computation methods~\cite{SchV18, AliBMetal20} and error
estimators~\cite{FenB19} can be used to determine successively new
interpolation points minimizing the $\Hinf$-error in a greedy
algorithm~\cite{BedBDetal20}.
On the other hand, there is the \emph{iterative rational Krylov algorithm
(IRKA)} for unstructured first-order systems, which solves the best
approximation problem in the $\Htwo$-norm~\cite{GugAB08}.
There is no extension for the general setting~\cref{eqn:sttf} in a
structure-preserving fashion.
However, the \emph{transfer function IRKA (TF-IRKA)} can be used to construct an
optimal unstructured (first-order) $\Htwo$-approximation for arbitrary transfer
functions~\cite{BeaG12}.
It has been shown to be very effective to use the final optimal interpolation
points from TF-IRKA in the structure-preserving interpolation setting
(\Cref{prp:strint}).
This basically resembles some of the ideas from~\cite{Wya12} for an extension
of the IRKA method to second-order systems.

A common drawback of interpolation methods is their global error behavior.
While these methods are exact in the interpolation points, the surrounding error
behavior may strongly differ depending on the actual transfer function.
While it might help to also match several derivatives in the interpolation
points using \Cref{prp:strint}, a quite often used approach is the averaging or
approximation of the constructed subspaces.
The general idea is to compute the solution to the linear systems in
\Cref{prp:strint} for a large amount of interpolation points, and then to
approximate the resulting large projection spaces by lower-order ones.
This can be done using, e.g., the pivoted QR decomposition or the singular value
decomposition (SVD); see, e.g.,~\cite{Wer21} for more details.
It has been shown in~\cite{BenGP19} that this oversampling procedure can be
used to recover the controllability and observability subspaces of general
systems like~\cref{eqn:sttf}.
Using some additional matrix products of the resulting basis matrices with the
system matrices, the authors of~\cite{BenGP19} can compute minimal
realizations of (linear) dynamical systems.
Additional truncation of the resulting reduced-order matrices has been shown to
still yield good results and is closely related to the idea of approximate
subspaces.


\subsection{AAA approximation of frequency-dependent functions}%
\label{sec:aaa}

The contributions of general nonlinear frequency-dependent functions acting on
parts of the system denoted by a set of constant matrices $\mv{C}_{i}$ have to 
be considered for problems of type~\cref{eqn:tfcasec} from Case~C in
\Cref{subsec:vibro_categories}.
To efficiently apply the structure-preserving interpolation framework from
\Cref{subsubsec:structure_preserving_interpolation} to such models,
the matrix functions in~\cref{eqn:sttf} can be represented by a series
expansion about an interpolation point $s_{0}$:
\begin{align} \label{eqn:seriesexp}
  \begin{aligned}
	  \stK\left(s + s_{0}\right) & = \sum\limits_{\ell = 0}^{\infty}
      s_{0}^{\ell} \stK_{\ell}, &
      \stF \left(s + s_{0}\right) & = \sum\limits_{\ell = 0}^{\infty}
      s_{0}^{\ell} \stF_{\ell}, &
      \stG\left(s + s_{0}\right) = \sum\limits_{\ell = 0}^{\infty}
      s_{0}^{\ell} \stG_{\ell}.
  \end{aligned}
\end{align}
While the series expansion factors can be computed straightforwardly for
polynomials, general nonlinear functions require special treatment.

In the following approach, the required series expansion factors are computed
from AAA approximations of the frequency-dependent functions.
AAA computes a rational interpolant $r(s)$ of a complex-valued function
$\phi(s)$ given function evaluations only~\cite{NakST18}.
The approximant $r(s)$ is originally given in barycentric form
\begin{align*}
	  r\left(s\right) & = \frac{\sum\limits_{j=1}^{q} \frac{w_{j} f_{j}}
      {s - s_{j}}}{\sum\limits_{j=1}^{q}\frac{w_{j}}{s-s_{j}}},
\end{align*}
where $q \geq 1$ is the order of approximation, $w_{j}$ are weights,
$f_{j}$ are data points, and $s_{j}$ are support points.
The barycentric form can be written in matrix notation as
\begin{align*}
	r\left(s\right) = \mv{a} \left( \mv{D} + s\mv{E} \right)^{-1} \mv{b}
\end{align*}
with
\begin{align*}
\begin{aligned}
	\mv{a} & =
		\begin{bmatrix}
			w_{1}f_{1} & \cdots & w_{q}f_{q}
		\end{bmatrix}, &
	\mv{b} & =
		\begin{bmatrix}
			1& 0& \cdots & 0
		\end{bmatrix}^{\trans},\\
	\mv{D} & =
		\begin{bmatrix}
			w_{1} & w_{2} & \cdots & w_{q-1} & w_{q} \\
			-s_{1} & s_{2} & & &\\
			& -s_{2} & \ddots & & \\
			& & \ddots & s_{q-1} & \\
			& & & -s_{q-1} & s_{q}
		\end{bmatrix}, &
	\mv{E} & =
		\begin{bmatrix}
			0 & 0 & \cdots & 0 & 0 \\
			1 & -1 &&& \\
			& 1 & \ddots && \\
			&& \ddots & -1 & \\
			&&& 1 & -1
		\end{bmatrix}.
\end{aligned}
\end{align*}
Shifting about $s_{0}$ yields
\begin{align} \label{eqn:aaa_matrix}
	r(s) & = \mv{a} \left( \mv{I}_{q} - \left(s - s_{0}\right) 
	  \widetilde{\mv{D}} \right)^{-1} \tilde{\mv{b}},
\end{align}
with $\widetilde{\mv{D}} = -\left( \mv{D} + s_0 \mv{E} \right)^{-1} \mv{E}$, 
$\tilde{\mv{b}} = \left( \mv{D} + s_0 \mv{E} \right)^{-1} \mv{b}$, and the 
$q$-dimensional identity matrix $\mv{I}_{q}$; see~\cite{LieMPetal21}. 
The involved matrix inverse is computationally inexpensive, as $q$ is typically
small.
This procedure has to be performed for all $k$ functions $\phi_{i}(s)$, yielding
a set of two vectors and a matrix $\{\mv{a}_{i}, \tilde{\mv{b}}_{i},
\widetilde{\mv{D}}_{i} \}_{i = 1}^{k}$ for each $\phi_{i}(s)$.
Expanding~\cref{eqn:aaa_matrix} into a Neumann series for all $k$ functions
$\phi_{i}(s)$ from~\cref{eqn:tfcasec} yields the sought after expansion factors
for~\cref{eqn:seriesexp}. 
For a Case C transfer function~\cref{eqn:tfcasec} and $\ell=2$, the series 
expansion of $\mathcal{K}_{\ell}$ in terms of~\cref{eqn:seriesexp} is 
given by
\begin{align*}
  \begin{aligned}
	  \stK_{0} &= s_0^2 \mv{M} + s_0 \mv{C} + \mv{K} + 
	  \sum\limits_{i = 0}^{k} \mv{a}_{i} \tilde{\mv{b}}_{i} \mv{C}_{i}, 
	  & \stK_{1} & = 2s_0 \mv{M} + \mv{C} + 
	  \sum\limits_{i = 0}^{k} \mv{a}_{i} \widetilde{\mv{D}}_{i} 
	  \tilde{\mv{b}}_{i} \mv{C}_{i},\\
	  \stK_{2} & = \mv{M} +
	  \sum\limits_{i = 0}^{k} \mv{a}_{i} \widetilde{\mv{D}}_{i}^2 
	  \tilde{\mv{b}}_{i} \mv{C}_{i}.
  \end{aligned}
\end{align*}
The matrix functions $\stF$ and $\stG$ in~\cref{eqn:seriesexp} are handled
analogously.

This approach can be used to fit transfer functions with non-polynomial
frequency-dependent terms into the standard second-order
Arnoldi scheme from \Cref{subsec:arnoldi}.
Here, only the expansion terms up to second order are considered and all
higher-order terms are truncated.
Enough expansion points have to be considered in the reduced model in order to
depict the frequency dependency as quadratic functions with reasonable accuracy
between the shifts.
A similar approach based on a Taylor series was successfully applied to
poroacoustic problems in~\cite{XieZJetal20}.


\begin{table}[t]
  \centering
  \caption{Second-order balanced truncation formulas. The $\ast$ denotes
    factors of the SVD not needed, and thus not accumulated in practical
    computations~\cite{BenW21b, Wer21}.}
  \label{tab:sobt}
  \begin{tabular}{cllc}
    \hline\noalign{\medskip}
      \multicolumn{1}{c}{\textbf{Type}}
      & \multicolumn{1}{c}{\textbf{SVD(s)}}
      & \multicolumn{1}{c}{\textbf{Truncation}}
      & \multicolumn{1}{c}{\textbf{Reference}}\\
    \noalign{\smallskip}\hline\noalign{\smallskip}
      \soV
        & $\mv{U} \mv{\Sigma} \mv{T}^{\trans} = \mv{L}_{\rm{v}}^{\trans} \mv{M} 
          \mv{R}_{\rm{v}}$
        & $\mv{W} = \mv{L}_{\rm{v}} \mv{U}_{1} \mv{\Sigma}_{1}^{-\frac{1}{2}},~
          \mv{V} = \mv{R}_{\rm{v}} \mv{T}_{1} \mv{\Sigma}_{1}^{-\frac{1}{2}}$
        & \cite{ReiS08}
        \vphantom{\Large $P^{\frac{1}{2}}_{\frac{1}{2}}$}\\
    \noalign{\smallskip}\hline\noalign{\smallskip}
      \soFV
      & $\ast\, \mv{\Sigma} \mv{T}^{\trans} = \mv{L}_{\rm{p}}^{\trans} 
        \mv{R}_{\rm{p}}$
      & $\mv{W} = \mv{V},~
        \mv{V} = \mv{R}_{\rm{p}} \mv{T}_{1} \mv{\Sigma}_{1}^{-\frac{1}{2}}$
      & \cite{MeyS96}%
      \vphantom{\Large $P^{\frac{1}{2}}_{\frac{1}{2}}$}\\
    \noalign{\smallskip}\hline\noalign{\smallskip}
      \soVPM
      & $\mv{U} \mv{\Sigma} \mv{T}^{\trans} = \mv{L}_{\rm{p}}^{\trans} 
        \mv{R}_{\rm{v}}$
      & $\mv{W} = \mv{M}^{-\trans} \mv{L}_{\rm{p}} \mv{U}_{1}
        \mv{\Sigma}_{1}^{-\frac{1}{2}},~
        \mv{V} = \mv{R}_{\rm{v}} \mv{T}_{1} \mv{\Sigma}_{1}^{-\frac{1}{2}}$
      & \cite{ReiS08}%
      \vphantom{\Large $P^{\frac{1}{2}}_{\frac{1}{2}}$}\\
    \noalign{\smallskip}\hline\noalign{\smallskip}
    \soPM
      & $\mv{U} \mv{\Sigma} \mv{T}^{\trans} = \mv{L}_{\rm{p}}^{\trans} 
        \mv{R}_{\rm{p}}$
      & $\mv{W} = \mv{M}^{-\trans} \mv{L}_{\rm{p}} \mv{U}_{1}
        \mv{\Sigma}_{1}^{-\frac{1}{2}},~
        \mv{V} = \mv{R}_{\rm{p}} \mv{T}_{1} \mv{\Sigma}_{1}^{-\frac{1}{2}}$
      & \cite{ReiS08}%
      \vphantom{\Large $P^{\frac{1}{2}}_{\frac{1}{2}}$}\\
    \noalign{\smallskip}\hline\noalign{\smallskip}
      \soPV
      & $\mv{U} \mv{\Sigma} \mv{T}^{\trans} = \mv{L}_{\rm{v}}^{\trans} \mv{M} 
        \mv{R}_{\rm{p}}$
      & $\mv{W} = \mv{L}_{\rm{v}} \mv{U}_{1} \mv{\Sigma}_{1}^{-\frac{1}{2}},~
         \mv{V} = \mv{R}_{\rm{p}} \mv{T}_{1} \mv{\Sigma}_{1}^{-\frac{1}{2}}$
      & \cite{ReiS08}%
      \vphantom{\Large $P^{\frac{1}{2}}_{\frac{1}{2}}$}\\
    \noalign{\smallskip}\hline\noalign{\smallskip}
      \soVP
      & {\arraycolsep = 1.4pt
        $\begin{array}{rcl}
           \ast\, \mv{\Sigma} \mv{T}^{\trans} & = & \mv{L}_{\rm{p}}^{\trans} 
           \mv{R}_{\rm{v}}, \vphantom{\text{\Large} P^{\frac{1}{2}}_{1}}\\
           \mv{U} \ast \ast & = & \mv{L}_{\rm{v}}^{\trans} \mv{M} 
           \mv{R}_{\rm{p}}
           \vphantom{\text{\Large} P^{\frac{1}{2}}_{\frac{1}{2}}}
         \end{array}$}
      & $\mv{W} = \mv{L}_{\rm{v}} \mv{U}_{1} \mv{\Sigma}_{1}^{-\frac{1}{2}},~
         \mv{V} = \mv{R}_{\rm{v}} \mv{T}_{1} \mv{\Sigma}_{1}^{-\frac{1}{2}}$
      & \cite{ReiS08}\\
    \noalign{\smallskip}\hline\noalign{\smallskip}
    \soP
      & {\arraycolsep = 1.4pt
        $\begin{array}{rcl}
           \ast\, \mv{\Sigma} \mv{T}^{\trans} & = & \mv{L}_{\rm{p}}^{\trans} 
             \mv{R}_{\rm{p}},
             \vphantom{\text{\Large} P^{\frac{1}{2}}_{1}}\\
           \mv{U} \ast \ast & = & \mv{L}_{\rm{v}}^{\trans} \mv{M} 
             \mv{R}_{\rm{v}}
             \vphantom{\text{\Large} P^{\frac{1}{2}}_{\frac{1}{2}}}
         \end{array}$}
      & $\mv{W} = \mv{L}_{\rm{v}} \mv{U}_{1} \mv{\Sigma}_{1}^{-\frac{1}{2}},~
         \mv{V} = \mv{R}_{\rm{p}} \mv{T}_{1} \mv{\Sigma}_{1}^{-\frac{1}{2}}$
      & \cite{ReiS08}\\
    \noalign{\smallskip}\hline\noalign{\smallskip}
      \soSO
      & {\arraycolsep = 1.4pt
        $\begin{array}{rcl}
          \mv{U}_{\rm{p}} \mv{\Sigma}_{\rm{p}} \mv{T}_{\rm{p}}^{\trans} & = &
            \mv{L}_{\rm{p}}^{\trans} \mv{R}_{\rm{p}},
            \vphantom{\text{\Large} P^{\frac{1}{2}}_{1}}\\
          \mv{U}_{\rm{v}} \mv{\Sigma}_{\rm{v}} \mv{T}_{\rm{v}} & = &
            \mv{L}_{\rm{v}}^{\trans} \mv{M} \mv{R}_{\rm{v}}
            \vphantom{\text{\Large} P^{\frac{1}{2}}_{\frac{1}{2}}}
         \end{array}$}
      & {\arraycolsep = 1.4pt
        $\begin{array}{rclrcl}
           \mv{W}_{\rm{p}} & = & \mv{L}_{\rm{p}} \mv{U}_{\rm{p},1}
             \mv{\Sigma}_{\rm{p},1}^{-\frac{1}{2}},~&
             \mv{V}_{\rm{p}} & = & \mv{R}_{\rm{p}} \mv{T}_{\rm{p},1}
             \mv{\Sigma}_{\rm{p},1}^{-\frac{1}{2}},\\
           \mv{W}_{\rm{v}} & = & \mv{L}_{\rm{v}} \mv{U}_{\rm{v},1}
             \mv{\Sigma}_{\rm{v},1}^{-\frac{1}{2}},~&
             \mv{V}_{\rm{v}} & = & \mv{R}_{\rm{v}} \mv{T}_{\rm{v},1}
             \mv{\Sigma}_{v,1}^{-\frac{1}{2}}
         \end{array}$}
      & \cite{ChaLVetal06}\\
      \noalign{\smallskip}\hline
  \end{tabular}
\end{table}

\subsection{Second-order balanced truncation methods}
\label{sec:sobt}

A different type of model reduction method is given with balanced truncation.
It was developed for first-order systems and utilizes the concepts of
controllability and observability to remove states, which have no big influence
on the input-to-output system behavior~\cite{Moo81}.
There have been several attempts on extending the balanced truncation approach
to second-order systems like~\cref{eqn:sotf};
see~\cite{ReiS08, MeyS96, ChaLVetal06}.
All are based on considering the linear first-order system
\begin{align*}
  \begin{aligned}
    s\mv{E} \mv{x}(s) & = \mv{A} \mv{x}(s) + \mv{B} \mv{u}(s), \\
    \mv{y}(s) & = \mv{D} \mv{x}(s),
  \end{aligned}
\end{align*}
with the system matrices concatenated from the original second-order
terms~\cref{eqn:sotf} in the following way:
\begin{align*}
  \begin{aligned}
    \mv{D} & = \begin{bmatrix} \mv{G} & 0 \end{bmatrix}, &
      \mv{E} & = \begin{bmatrix} \mv{I}_{n} & 0 \\ 
        0 & \mv{M} \end{bmatrix}, &
      \mv{A} & = \begin{bmatrix} 0 & \mv{I}_{n} \\ 
        -\mv{K} & -\mv{C} \end{bmatrix}, &
      \mv{B} & = \begin{bmatrix} 0 \\ \mv{F} \end{bmatrix}.
  \end{aligned}
\end{align*}
Using the solutions of the dual Lyapunov equations
\begin{align*}
  \mv{A} \mv{P} \mv{E}^{\trans} + \mv{E} \mv{P} \mv{A}^{\trans} + \mv{B} 
    \mv{B}^{\trans} & = 0,\\
  \mv{A}^{\trans} \mv{Q} \mv{E} + \mv{E}^{\trans} \mv{Q} \mv{A} + 
    \mv{D}^{\trans} \mv{D} & = 0,
\end{align*}
the truncation matrices for projection-based model reduction in
\Cref{subsec:proj}, can be build using partitioned (low-rank) Cholesky
factorizations
\begin{align*}
  \begin{aligned}
    \mv{P} & = \begin{bmatrix} \mv{R}_{\rm{p}} \\ \mv{R}_{\rm{v}} \end{bmatrix}
      \begin{bmatrix} \mv{R}_{\rm{p}} \\ \mv{R}_{\rm{v}} \end{bmatrix}^{\trans} 
      & \text{and} &&
      \mv{Q} & = \begin{bmatrix} \mv{L}_{\rm{p}} \\
        \mv{L}_{\rm{v}} \end{bmatrix}
      \begin{bmatrix} \mv{L}_{\rm{p}} \\ \mv{L}_{\rm{v}} \end{bmatrix}^{\trans}
  \end{aligned}
\end{align*}
and the formulas from \Cref{tab:sobt}.
The last line of \Cref{tab:sobt} describes a different approach, which
is a projection method on the corresponding first-order realization followed by
structure recovery of the second-order system.
Therein, the reduced-order model is computed by
\begin{align}\label{eqn:sobtproj}
  \begin{aligned}
    \Mr  & = \mv{S} \left( \mv{W}_{\rm{v}}^{\trans} \mv{M} \mv{V}_{\rm{v}} 
      \right) \mv{S}^{-1}, &
    \Cr & = \mv{S} \left( \mv{W}_{\rm{v}}^{\trans} \mv{C} \mv{V}_{\rm{v}} 
      \right) \mv{S}^{-1}, &
    \Kr & = \mv{S} \left( \mv{W}_{\rm{v}}^{\trans} \mv{K} \mv{V}_{\rm{p}} 
      \right),\\
    \Fr & = \mv{S} \left( \mv{W}_{\rm{v}}^{\trans} \mv{F} \right), &
    \Gr & = \mv{G} \mv{V}_{\rm{p}}, &
  \end{aligned}
\end{align}
with $\mv{S} = \mv{W}_{\rm{p}} \mv{V}_{\rm{v}}$ and the truncation matrices
$\mv{W}_{\rm{p}}$, $\mv{W}_{\rm{v}}$, $\mv{V}_{\rm{p}}$, $\mv{V}_{\rm{v}}$
from the last row of \Cref{tab:sobt}.

In cases where the computation of both Gramian factors is not feasible, e.g., if
the system output is the complete state, or when only a one-sided projection is
desired, dominant subspaces can be a useful tool.
Similar to the approximate subspaces for interpolation methods
(\Cref{subsec:interp}), dominant subspaces have been shown to be effective to
approximate the controllability or observability behavior of the system.
Therefor, pivoted QR decompositions or SVDs of one or both of the Gramian
factors are computed and truncated to the desired reduced order.
Then, the first $r$ rows of the resulting truncated orthogonal bases form the
truncation matrices for model reduction by projection.

\section{Numerical experiments}%
\label{sec:numerics}

In the following, we apply the model order reduction methods described in
\Cref{sec:mor} to different vibro-acoustic systems categorized as
Cases~A,~B and~C in \Cref{subsec:vibro_categories}.
The setup of the numerical comparison is described in detail in the upcoming
subsection.

The experiments reported here have been executed on single nodes of the Leibniz
supercomputing centre's Linux cluster \mbox{CoolMUC-2} running on 
SUSE\textsuperscript{\textregistered} Linux Enterprise Server~15~SP1.
Each node is equipped with a $28$ core Intel\textsuperscript{\textregistered}
Haswell based CPU and for each experiment we had access to $56$\,GB main memory.
All algorithms and experiments have been implemented in MATLAB
9.8.0.1451342 (R2020a Update~5).
The MATLAB toolboxes M-M.E.S.S. version~2.0.1~\cite{BenKS21, SaaKB20} and
SOLBT version~3.0~\cite{BenW21a} have been used in some of the experiments.
The models and data have been created with Kratos
Mutliphysics~8.1~\cite{DadRO10, FerBRetal20}.

\begin{center}%
  \setlength{\fboxsep}{5pt}%
  \fbox{%
  \begin{minipage}{.92\linewidth} \small
    \textbf{Code and data availability}\newline
    The used data, the source code of the implementations used for the numerical
    experiments and the computed results are available at
    \vspace{.5\baselineskip}
    \begin{center}
        \href{https://doi.org/10.5281/zenodo.5836047}%
          {\texttt{doi:10.5281/zenodo.5836047}}
    \end{center}
    \vspace{.5\baselineskip}
    under the BSD-2-Clause license, authored by Quirin Aumann and
    Steffen W. R. Werner.
  \end{minipage}}
\end{center}


\subsection{Experimental setup}

For a clear comparison of the various different model reduction methods, we are
using the \morscore~\cite{Him21, Wer21}.
In principle, it compresses the behavior of error-per-order graphs into single
scalar values, which can then be easily compared for different methods and error
measures.
Given a relative error graph $(r, \varepsilon(r)) \in \N_{0} \times (0, 1]$,
which relates the reduced order $r$ to the relative approximation error
$\varepsilon(r)$, the \morscore{} is the area below the normalized error graph
$(\varphi_{r}, \varphi_{\varepsilon(r)})$ where
\begin{align*}
  \begin{aligned}
    \varphi_{r}\colon r \mapsto \frac{r}{r_{\max}} &&
      \text{and} &&
      \varphi_{\varepsilon(r)}\colon \varepsilon(r) \mapsto \frac{\log_{10}\big(
      \varepsilon(r) \big)}{\lfloor \log_{10}\big( \epsilon \big) \rfloor}.
  \end{aligned}
\end{align*}
Thereby, $r_{\max}$ denotes the maximum reduced order for the comparison and
$\epsilon$ is a reasonably small tolerance describing the approximation accuracy
that shall be reached.

For the following computations of the \morscore, we approximate the relative
error under the $\Linf$-norm via
\begin{equation}
	\varepsilon(r) = \frac{\max\limits_{\omega \in [\omega_{\min},\omega_{\max}]}
    \left\lVert \tf(\omega \i) - \tfr_{r}(\omega \i) \right\lVert_{2}}
    {\max\limits_{\omega \in [\omega_{\min},\omega_{\max}]}
    \left\lVert \tf(\omega \i) \right\lVert_{2}}
    \approx \frac{\linfnorm{\tf - \tfr_{r}}}{\linfnorm{\tf}},
\end{equation}
where $\tfr_{r}$ is the transfer function of a reduced-order model of size
$r$.
For simplicity, we denote the errors in plots with the $\Linf$-norm.

We use the following structure-preserving methods in our numerical
comparison:
\begin{description}
  \item[\morequi] is the structure-preserving interpolation framework
    from~\cite{BeaG09} with linearly equidistant interpolation points on
    the imaginary axis in a frequency range of interest
    $[\omega_{\min}, \omega_{\max}] \i$.
    If suitable, the interpolation points are supplemented by educated guesses
    with high impact transfer function behavior.
  \item[\moravg] is the approximate/averaged subspace approach, described in
    \Cref{subsubsec:structure_preserving_interpolation}.
    The linear systems necessary for interpolation as in \Cref{prp:strint} are
    solved and collected into a matrix with $q \geq r$ columns, which is then
    approximated by a pivoted QR decomposition of order $r$, from which we use
    the resulting orthogonal basis to obtain the final reduced-order model.
  \item[$\boldsymbol{\mathsf{L}_{\infty}}$] denotes greedy interpolation
    algorithms based on minimizing the $\mathcal{L}_{\infty}$-norm. Here, new 
    expansion points for the reduced-order model are iteratively chosen at
    locations where the difference between the transfer functions of the original
    model and the reduced-order model are maximal.
    Due to the unstable numerical behavior of currently available
    $\mathcal{L}_{\infty}$ computation routines, we save in a presampling step
    the evaluations of the transfer function on the imaginary axis and compute
    the $\mathcal{L}_{\infty}$-error of the reduced-order models with respect to
    this discrete set.
  \item[\morminrel] follows the minimal realization algorithm
    from~\cite{BenGP19}.
    Dominant subspaces of a potentially minimal realization are computed by
    choosing a reasonable large number of expansion points $q > r$ on the
    imaginary axis, from which reduced-order models are computed using the SVD.
  \item[\morbt] are the second-order balanced truncation formulas according
    to \Cref{tab:sobt}.
\end{description}
Further potential model reduction methods for this comparison could be based
on (TF\=/)\allowbreak IRKA.
However, in our numerical experiments, the used implementations of
(TF\=/)\allowbreak IRKA did not converge and took an exceptional amount of 
computation time for the different requested reduced orders.
Therefor, results of these model reduction methods are not presented here.

The interpolation-related methods \moravg{}, \morlinf{} and \morminrel{} require
a presampling step in which a reduction basis of order $q \geq r$ or the
evaluation of the transfer function at chosen expansion points are computed
before the reduction methods are applied.
We employ three different approaches to compute this database for the following
experiments:
\begin{itemize}
  \item For the standard method, classical transfer function interpolation in
    $q$ different expansion points is employed, which are considered in a
    frequency range of interest $[\omega_{\min}, \omega_{\max}] \i$.
    Thereby, $q$ solutions of linear systems of order $n$ are required for each
    presampling basis $\mv{V}_{\rm{pre}}$ and $\mv{W}_{\rm{pre}}$.
    This method is not specially denoted in the final naming of the model
    reduction methods.
  \item The second method, denoted by \morsp{}, is based on the
    higher-order structure-preserving Hermite interpolation scheme
    summarized in \Cref{subsubsec:structure_preserving_interpolation}.
    The interpolation order $\ell$ at each expansion point is set, such that all 
    derivatives of the transfer function factors would vanish for $\ell+1$.  
    This means, that standard second-order systems described as Case~A and 
    Case~B have an interpolation order of $\ell=2$, Case~C systems may have 
    a higher $\ell$.
    Following the approach from
    \Cref{subsubsec:structure_preserving_interpolation}, this results in an 
    $\left(\ell+1\right)$-term recurrence, and $\left(\ell+1\right)m$ or
    $\left(\ell+1\right)p$ columns for $\mv{V}_{\rm{pre}}$ or 
    $\mv{W}_{\rm{pre}}$, respectively, are computed per expansion point 
    requiring only one decomposition of an $n\times n$ matrix each.
    Consequently, the order of models reduced with \morlinf{} is 
    always a factor of $\max \left\{\left(\ell+1\right)m, \left(\ell+1\right)p 
    \right\}$, as all columns associated with the chosen shift are selected for 
    the reduction basis.
    For \morsp{}, the involved derivatives of the matrix-valued functions are
    computed analytically.
  \item Last, we employ an arbitrary Hermite interpolation order at each
    expansion point using the second-order Arnoldi method presented in
    \Cref{subsec:arnoldi}, further denoted by \moraaa{}.
    Similarly to \morsp{}, this implies that multiple columns of the reduction
    basis are computed in the same step and that the reduced models for
    \morlinf{} grow by a factor of the interpolation order.
    The AAA-based method from \Cref{sec:aaa} is used to obtain a
    second-order representation of systems with a Case C transfer function,
    such that the second-order Arnoldi method is applicable.
\end{itemize}

We also compare different projections to assess of controllability and 
observability of the various systems.
Apart from the two-sided projection, we also consider one-sided projections
regarding the system input and output, respectively.
The input projection is obtained by setting $\mv{W}=\mv{V}$, the 
output projection by setting $\mv{V} = \mv{W}$.
Where applicable, we also compare complex- and real-valued projection bases.
A real-valued basis is obtained from the initial complex-valued basis by 
$\mv{V} = \begin{bmatrix} \re{\mv{V}} & \im{\mv{V}} \end{bmatrix}$ and 
$\mv{W} = \begin{bmatrix} \re{\mv{W}} & \im{\mv{W}} \end{bmatrix}$.
The considered projections are:
\begin{description}
	\item[\tsimag] Two-sided projection $\mv{W} \neq \mv{V}$ with $\mv{V}, 
	  \mv{W} \in \C^{n \times r}$,
	\item[\tsreal] Two-sided projection $\mv{W} \neq \mv{V}$ with $\mv{V}, 
	  \mv{W} \in \R^{n \times r}$,
	\item[\osimaginput] Single-sided projection $\mv{W} = \mv{V}$ with $\mv{V}, 
	  \mv{W} \in \C^{n \times r}$,
	\item[\osrealinput] Single-sided projection $\mv{W} = \mv{V}$ with $\mv{V}, 
	  \mv{W} \in \R^{n \times r}$,
	\item[\osimagoutput] Single-sided projection $\mv{V} = \mv{W}$ with 
	  $\mv{V}, \mv{W} \in \C^{n \times r}$,
	\item[\osrealoutput] Single-sided projection $\mv{V} = \mv{W}$ with 
	  $\mv{V}, \mv{W} \in \R^{n \times r}$.
\end{description}
Consequently, a real-valued projection incorporating \morequi{} yields reduced 
models with an even dimension $r$ only.
In order to obtain a real-valued projection using \moravg{}, \morlinf{} and 
\morminrel{}, their presampled bases $\mv{V}_{\rm{pre}}$ and $\mv{W}_{\rm{pre}}$ 
are modified such that they are of dimension $n \times 2q$ before the methods 
compute the truncated projection matrices $\mv{V}, \mv{W} \in \R^{n \times r} $.
To allow a one-sided projection regarding the system input for \morbt{}, only 
the controllability Gramian is computed and used as left and right 
projection matrices in orthogonalized and truncated form.
A one-sided projection regarding the system output is analogously possible by
computing the observability Graminan.
No complex-valued bases are considered for \morbt{} since the method is only
applied to real-valued systems, for which it preserves the realness.
However, the \morbt{} methods are in the following assigned to the imaginary
bases methods since no additional realification process is employed.


\subsection{Vibration of a plate with distributed mass-spring systems}

\begin{figure}[t]
	\centering
	\begin{subfigure}[b]{.4\textwidth}
		\centering
		\def\svgwidth{\textwidth}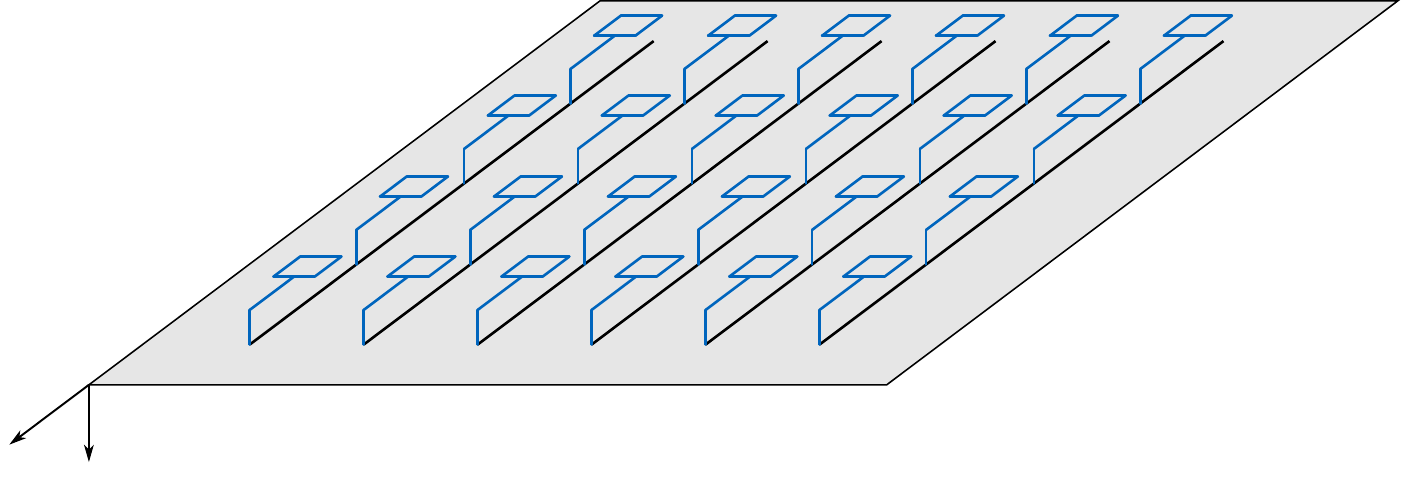
		\vspace{1.5cm}
		\subcaption{Sketch of the plate with TVAs.}
		\label{fig:plate_sketch}
	\end{subfigure}%
	\hfill%
	\begin{subfigure}[b]{.55\textwidth}
		\centering
  \tikzexternalenable%
  \tikzsetnextfilename{plate48}%
  \filemodCmp{graphics/plate48.tex}{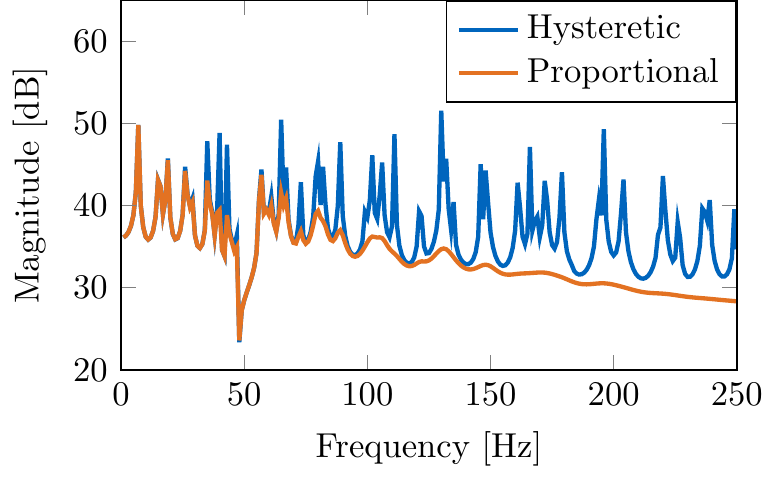}%
  {\tikzset{external/remake next}}{}%
  \begin{tikzpicture}[font = \small]
  \begin{axis}[
    height = .225\textheight,
    width  = .95\textwidth,
    xlabel = {Frequency \ubraces{\si{\hertz}}},
    ylabel = {Magnitude \ubraces{\si{\decibel}}},
    xmin = 0,
    xmax = 250,
    ymin = 20,
    ymax = 65,
    xtick distance = 50,
    ytick distance = 10,
	axis on top,
    legend style = {
      at = {(1, 1)},
      anchor = north east},
    legend cell align = left,
    every axis plot/.append style={line width=\mylinewidth}
  ]

	\addplot[mark=none, color=TUMBlau] table[x=freq, y expr={10*log10(\thisrow{hysteretic}/1e-9)}, col sep=comma] {graphics/data/plate_tf.csv};
	\addlegendentry{Hysteretic}

    \addplot[mark=none, color=TUMOrange] table[x=freq, y expr={10*log10(\thisrow{rayleigh}/1e-9)}, col sep=comma] {graphics/data/plate_tf.csv};
  	\addlegendentry{Proportional}

  \end{axis}
\end{tikzpicture}%
  \tikzexternaldisable%

		\subcaption{Root mean square of displacement for both plates with different
			damping mechanisms.}
		\label{fig:plate_tf}
	\end{subfigure}
	
	\caption{Sketch and transfer function of the vibrating plates. }
\end{figure}

The vibration response of simply supported strutted plates excited by a point 
load is modeled in this example.
The plates have dimensions of $0.8 \times \SI{0.8}{\meter} $, a thickness of
$t=\SI{1}{\milli\meter}$, and are made out of aluminum with the material
parameters $E=\SI{69}{\giga\pascal}$,
$\rho=\SI{2650}{\kilogram\per\cubic\meter}$ and $\nu=\num{0.22}$.
Two damping models are considered: proportional damping with
$\alpha = \num{0.01}$ and $\beta = \num{1e-4}$, and hysteretic damping with
$\eta = \num{0.001}$.
The plates are equipped with arrays of tuned vibration absorbers (TVA)
reducing their vibration response in the frequency region of the TVAs' tuning
frequency $f = \SI{48}{\hertz}$.
The TVAs are placed on the struts of the plates and are modeled as discrete
spring-damper elements with attached point masses.
In total, an extra mass of \SI{10}{\percent} of the plate structure's mass is 
added by the TVAs.
A point load near a corner of the plate with amplitude \SI{0.1}{\newton} excites
the system.
The model is sketched in \Cref{fig:plate_sketch}.
A similar system has been experimentally examined in~\cite{ClaDP16}.
The effect of the TVAs is limited to the frequency region directly adjacent to
the tuning frequency and is clearly visible in the frequency response plot of
the root mean square of the displacement on the plate surface in
\Cref{fig:plate_tf}.
The two damping models have a large influence on the respective transfer
functions.
The discretized system has an order of $n=\num{201900}$ and is evaluated in a 
frequency range of \SIrange{1}{250}{\hertz}.
While the poles of the proportionally damped system are only visible in the
lower frequency region, the hysteretically damped system's transfer function
shows many peaks over the complete frequency range of interest.
As only structural loads excite the system, Case~A transfer functions are used
to describe the output of both systems.
All system matrices are symmetric, respectively complex symmetric for 
the case of hysteretic damping, as no interaction effects between structure and 
fluid are modeled.
In order to evaluate the root mean square of the displacement at all points on
the plate surface, the displacement at these locations needs to be recovered
from the reduced space.
This is done using an output matrix $\mv{G}$ with dimensions $p \times n$, where
$p$ is the number of nodes on the plate surface mapping the result of each node
to an individual output.
Due to its size we only consider projections regarding the system 
input, i.e., \osimaginput{} and \osrealinput.

The models are reduced using all methods described in the beginning of this 
section, except \morbt{} not being applicable to the hysteretically damped 
model, because $\mv{C} = 0$ in this system.
The presampling basis for \morminrel{}, \moravg{} and \morlinf{}
considers $n_{\rm{s}} = 250$ frequency shifts distributed linearly in
$2\pi\i [1, 250]$.
Since the models are described by standard second-order transfer functions,
\morsp{} yields 3 columns for each interpolation point.
Using $n_{\rm{s}} = 80$ shifts linearly distributed in the same range and augmented
by shifts at $2\pi\i [46, 47, 48, 50]$ yields the intermediate reduction basis 
with $q = 250$. 
The additional shifts are introduced to capture the local behavior near the 
tuning frequency of the TVAs. 
A local order of $10$ is chosen for \moraaa{} presampling.
The intermediate basis of order $q = 250$ is computed considering $n_{\rm{s}} = 21$
shifts linearly distributed in $2\pi\i [1, 250]$ and four additional shifts at 
$2\pi\i [46, 47, 48, 50]$.
The expansion point sampling for \morequi{} is modified similarly to the 
presampling methods to account for the high impact of the TVAs on the transfer 
function near their tuning frequency.
The shifts at $2\pi\i [46, 47, 48, 49, 50]$ are always considered, the location
of the remaining shifts are linearly distributed in the frequency range of
interest.
For orders $r < 5$ only the first $r$ extra shifts were considered.

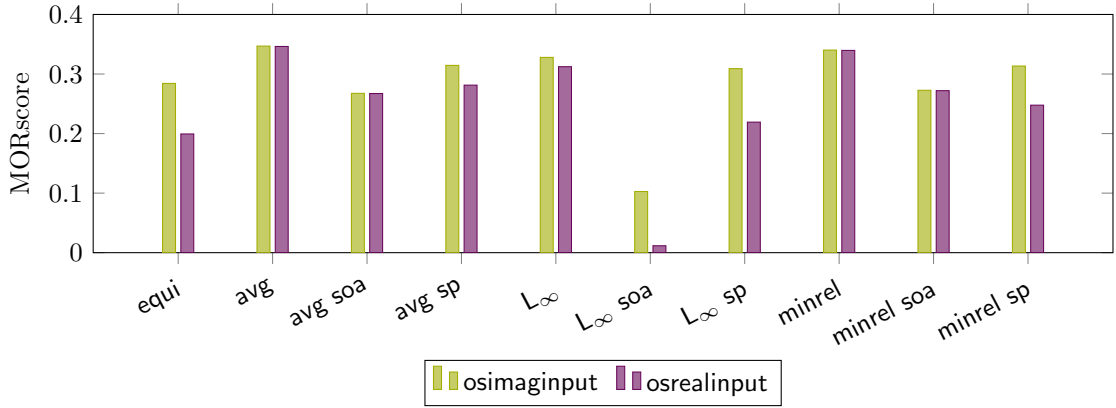
\begin{figure}[t]
	\centering
  \tikzexternalenable%
  \tikzsetnextfilename{hysteretic_morscore}%
  \filemodCmp{graphics/hysteretic_morscore.tex}{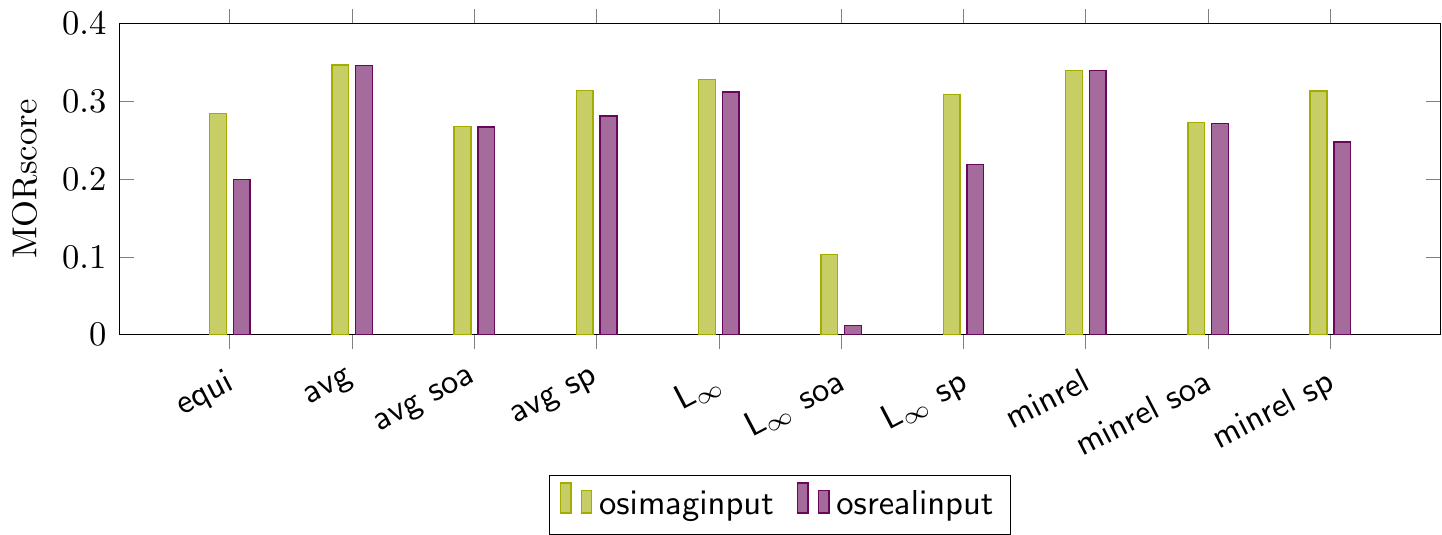}%
  {\tikzset{external/remake next}}{}%
  \begin{tikzpicture}[font = \small]
	\begin{axis}[height=0.2\textheight, width=\textwidth,
		ybar,
		ylabel = {\morscore{}},
		xtick distance = 1,
		ytick distance = .1,
		ymin=0, ymax=.4,
		legend cell align = left,
		legend style={at={(0.5,-0.45)},anchor=north, /tikz/every even column/.append style={column sep=1ex}},
		legend columns=-1,
		no markers,
		every axis plot/.append style={fill},
		cycle list name=TUMcolorlist_simple,
		bar width=0.17cm,
		symbolic x coords={\lmorequi{}, \lmoravg{}, \lmoravg{} \lmoraaa{}, \lmoravg{} \lmorsp{}, \lmorlinf{}, \lmorlinf{} \lmoraaa{}, \lmorlinf{} \lmorsp{}, \lmorminrel{}, \lmorminrel{} \lmoraaa{}, \lmorminrel{} \lmorsp{}},
		xticklabel style={rotate=27.5,anchor=north east},
	]
		\pgfplotsset{cycle list shift=2}

		\addplot+[fill opacity=0.6] table [x=method, y index=1, col sep=comma] {graphics/data/plate_48_hysteretic_res_mor_score_osimaginput.csv};
		\addlegendentry{\losimaginput{}}

		\addplot+[fill opacity=0.6] table [x=method, y index=1, col sep=comma] {graphics/data/plate_48_hysteretic_res_mor_score_osrealinput.csv};
		\addlegendentry{\losrealinput{}}

	\end{axis}
\end{tikzpicture}%
  \tikzexternaldisable%

	\caption{\morscores{} of all employed reduction methods for 
		the hysteretically damped plate, with maximum accuracy
    $\epsilon = \num{1e-6}$ and maximum order $r_{\rm{max}} = 250$.}
	\label{fig:hysteretic_morscore}
\end{figure}

\begin{figure}[t]
	\centering
  \tikzexternalenable%
  \tikzsetnextfilename{hysteretic_err_hinf}%
  \filemodCmp{graphics/hysteretic_err_hinf.tex}{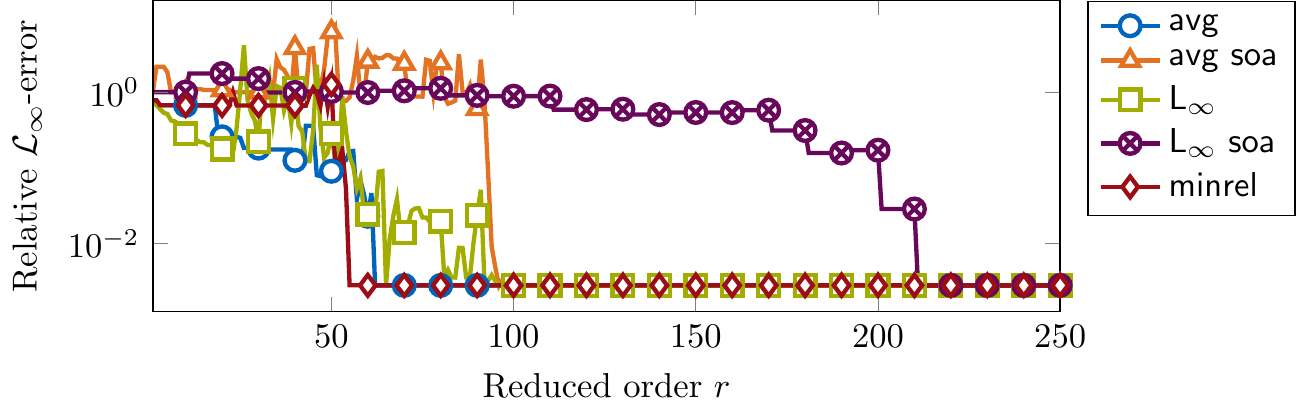}%
  {\tikzset{external/remake next}}{}%
  \begin{tikzpicture}[font = \small]
\begin{axis}[height=0.2\textheight, width=.72\textwidth, ymode=log,
		xlabel = {Reduced order $r$},
		ylabel = {\errlabel},
		xmin=1, xmax=250,
		axis on top,
		legend pos = outer north east, legend cell align = left,
		cycle list name=TUMcolorlist,
		mark repeat={10},
		mark phase=10,
	]

	\addplot+ table[x=r, y=strint_avg, col sep=comma] {graphics/data/plate_48_hysteretic_res_hinfrelerr_osimaginput.csv};
	\addlegendentry{\lmoravg{}}

	\addplot+ table[x=r, y=strint_avg_aaaa, col sep=comma] {graphics/data/plate_48_hysteretic_res_hinfrelerr_osimaginput.csv};
	\addlegendentry{\lmoravg{} \lmoraaa{}}


	\addplot+ table[x=r, y=strint_linf, col sep=comma] {graphics/data/plate_48_hysteretic_res_hinfrelerr_osimaginput.csv};
	\addlegendentry{\lmorlinf{}}

	\addplot+ table[x=r, y=strint_linf_aaaa, col sep=comma] {graphics/data/plate_48_hysteretic_res_hinfrelerr_osimaginput.csv};
	\addlegendentry{\lmorlinf{} \lmoraaa{}}


	\addplot+ table[x=r, y=minrel, col sep=comma] {graphics/data/plate_48_hysteretic_res_hinfrelerr_osimaginput.csv};
	\addlegendentry{\lmorminrel{}}

%
\end{axis}
\end{tikzpicture}%
  \tikzexternaldisable%

	\caption{Relative \Linf-error for reduced models of the plate model
		with hysteretic damping computed by several reduction methods and 
		\osimaginput{} projection. }
	\label{fig:hysteretic_hinf}
\end{figure}
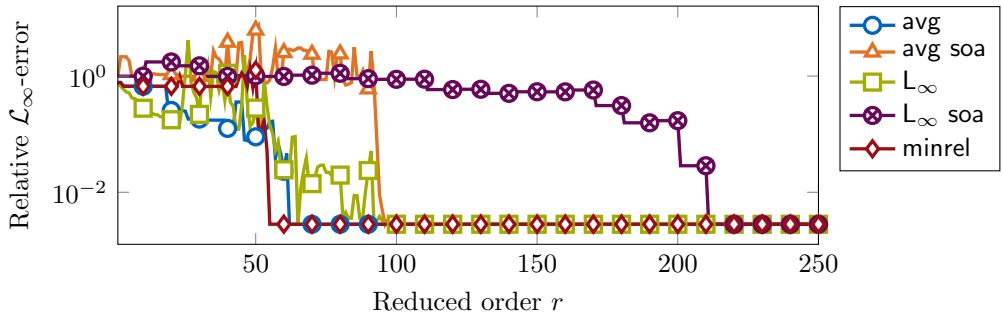

We consider the results for the plate with hysteretic damping.
Despite the high number of weakly damped poles in the transfer function, all
applicable methods are able to compute reasonably accurate reduced-order models.
The \morscores{} referenced to $\epsilon = \num{1e-6}$ and $r_{\rm{max}} = 250$
are given in \Cref{fig:hysteretic_morscore}.
The choice of the tolerance value is motivated by the fact that due to the
conditioning of the example the relative approximation errors do not drop below
$\num{1e-3}$ for any employed method.
Choosing a smaller $\epsilon$ would hinder the proper comparison of the model
order reduction methods.
The projections with complex-valued basis matrices yield good results for all
reduction methods, only \morlinf{}~\moraaa{} falls short.
It has to be noted, that all reduced models computed from a \moraaa{}
presampling need higher reduced orders $r$ to be as accurate as the other
methods.
This comes from the focus on derivatives of the transfer function in a
smaller number of expansion points compared to the other methods.
Employing real-valued projections yields comparable \morscores{}.
When using \morlinf{} in combination with \moraaa{}, the order of the reduced
models is increased in steps of $10$, which is the chosen size of the
second-order Krylov subspaces employed.
Therefor, larger reduced models are constructed in comparison to the other
methods, but less computational effort is required for the presampling process.
\Cref{fig:hysteretic_hinf} shows the (approximate) relative $\Linf$-errors
plotted over the reduced order.
It can be seen that all methods including \morlinf{}~\moraaa{} are able to
compute reduced models of the same accuracy given a large enough
reduced order.

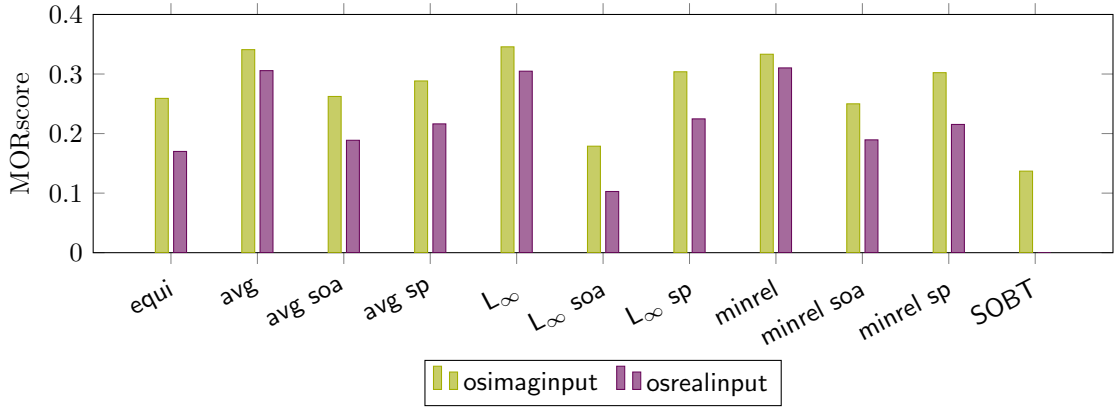
\begin{figure}[t]
	\centering
  \tikzexternalenable%
  \tikzsetnextfilename{rayleigh_morscore}%
  \filemodCmp{graphics/rayleigh_morscore.tex}{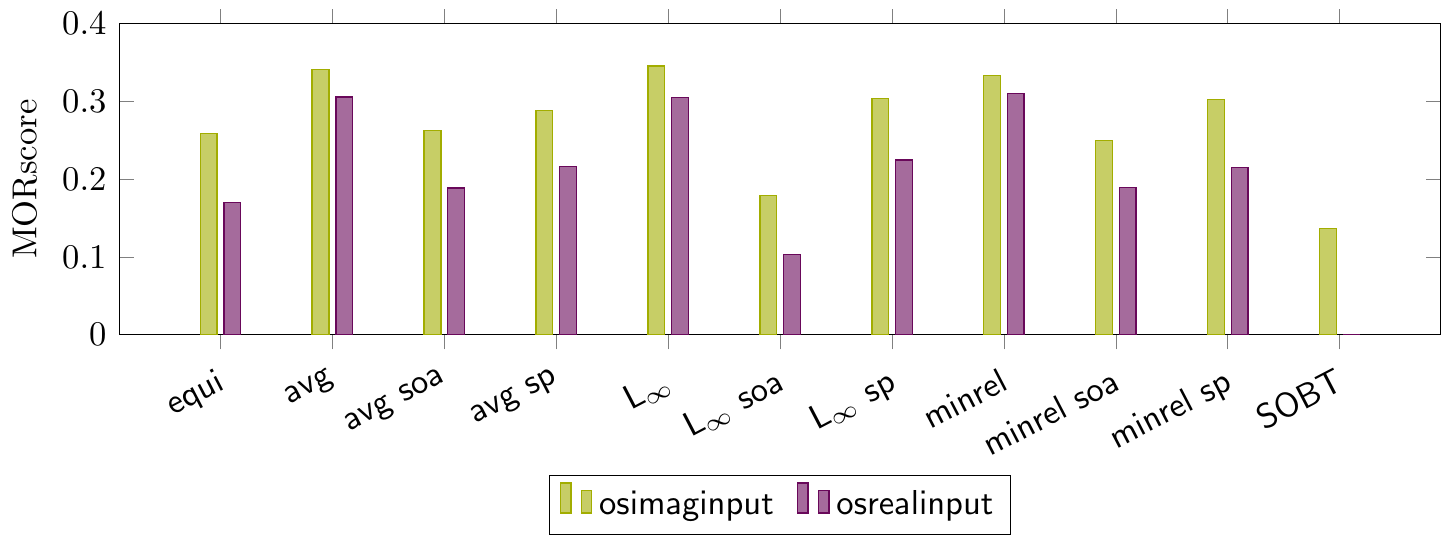}%
  {\tikzset{external/remake next}}{}%
  \begin{tikzpicture}[font = \small]
	\begin{axis}[height=0.2\textheight, width=\textwidth,
		ybar,
		ylabel = {\morscore{}},
		xtick distance = 1,
		ytick distance = .1,
		ymin=0, ymax=.4,
		legend cell align = left,
		legend style={at={(0.5,-0.45)},anchor=north, /tikz/every even column/.append style={column sep=1ex}},
		legend columns=-1,
		no markers,
		every axis plot/.append style={fill},
		cycle list name=TUMcolorlist_simple,
		bar width=0.17cm,
		symbolic x coords={\lmorequi{}, \lmoravg{}, \lmoravg{} \lmoraaa{}, \lmoravg{} \lmorsp{}, \lmorlinf{}, \lmorlinf{} \lmoraaa{}, \lmorlinf{} \lmorsp{}, \lmorminrel{}, \lmorminrel{} \lmoraaa{}, \lmorminrel{} \lmorsp{}, \lmorbt{}},
		xticklabel style={rotate=27.5,anchor=north east},
		xmin = {\lmorequi{}},
		xmax = { \lmorbt{}},
		enlarge x limits=0.09,
		]
		\pgfplotsset{cycle list shift=2}

		\addplot+[fill opacity=0.6] table [x=method, y index=1, col sep=comma] {graphics/data/plate_48_rayleigh_res_mor_score_osimaginput.csv};
		\addlegendentry{\losimaginput{}}

		\addplot+[fill opacity=0.6] table [x=method, y index=1, col sep=comma] {graphics/data/plate_48_rayleigh_res_mor_score_osrealinput.csv};
		\addlegendentry{\losrealinput{}}

	\end{axis}
\end{tikzpicture}%
  \tikzexternaldisable%

	\caption{\morscores{} of all employed reduction for the proportionally damped
    plate, with maximum accuracy $\epsilon = \num{1e-16}$ and maximum order
    $r_{\rm{max}} = 250$.}
	\label{fig:rayleigh_morscore}
\end{figure}

\begin{figure}[t]
	\centering
  \tikzexternalenable%
  \tikzsetnextfilename{rayleigh_err_hinf}%
  \filemodCmp{graphics/rayleigh_err_hinf.tex}{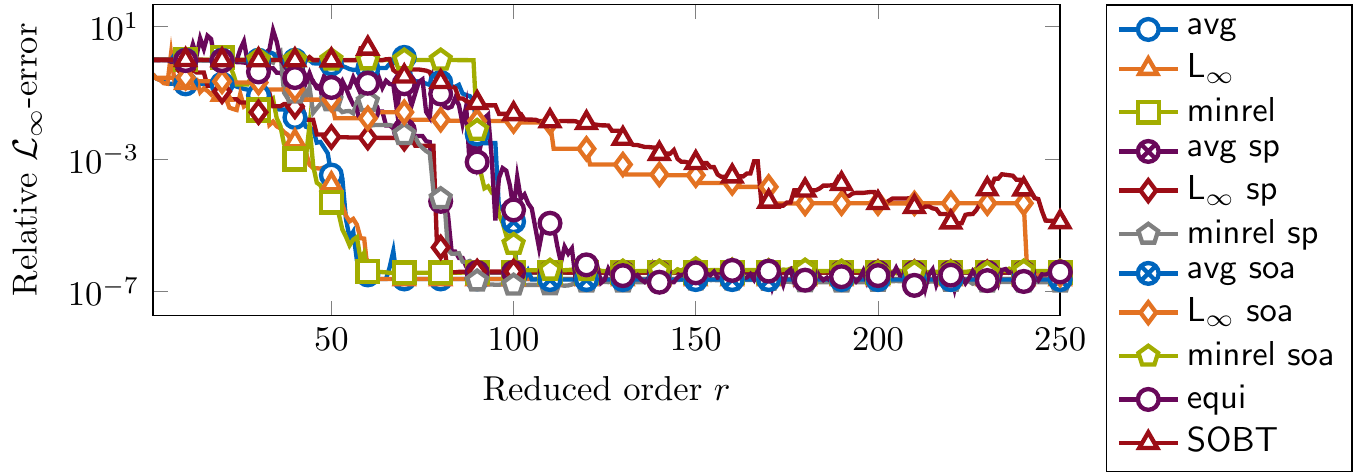}%
  {\tikzset{external/remake next}}{}%
  \begin{tikzpicture}[font = \small]
\begin{axis}[height=0.2\textheight, width=.72\textwidth, ymode=log,
	xlabel = {Reduced order $r$},
	ylabel = {\errlabel},
	xmin=1, xmax=250,
	axis on top,
	legend cell align = left,
	legend style={at={(1.05,1)},anchor=north west},
	cycle list name=TUMcolorlist,
	mark repeat={10},
	mark phase=10,
	]
	\addplot+ table[x=r, y=strint_avg, col sep=comma] {graphics/data/plate_48_rayleigh_res_hinfrelerr_osimaginput.csv};
	\addlegendentry{\lmoravg{}}

	\addplot+ table[x=r, y=strint_linf, col sep=comma] {graphics/data/plate_48_rayleigh_res_hinfrelerr_osimaginput.csv};
	\addlegendentry{\lmorlinf{}}

	\addplot+ table[x=r, y=minrel, col sep=comma] {graphics/data/plate_48_rayleigh_res_hinfrelerr_osimaginput.csv};
	\addlegendentry{\lmorminrel{}}

	\addplot+ table[x=r, y=strint_avg_strprs, col sep=comma] {graphics/data/plate_48_rayleigh_res_hinfrelerr_osimaginput.csv};
	\addlegendentry{\lmoravg{} \lmorsp{}}

	\addplot+ table[x=r, y=strint_linf_strprs, col sep=comma] {graphics/data/plate_48_rayleigh_res_hinfrelerr_osimaginput.csv};
	\addlegendentry{\lmorlinf{} \lmorsp{}}

	\addplot+ table[x=r, y=minrel_strprs, col sep=comma] {graphics/data/plate_48_rayleigh_res_hinfrelerr_osimaginput.csv};
	\addlegendentry{\lmorminrel{} \lmorsp{}}

	\addplot+ table[x=r, y=strint_avg_aaaa, col sep=comma] {graphics/data/plate_48_rayleigh_res_hinfrelerr_osimaginput.csv};
	\addlegendentry{\lmoravg{} \lmoraaa{}}

	\addplot+ table[x=r, y=strint_linf_aaaa, col sep=comma] {graphics/data/plate_48_rayleigh_res_hinfrelerr_osimaginput.csv};
	\addlegendentry{\lmorlinf{} \lmoraaa{}}

	\addplot+ table[x=r, y=minrel_aaaa, col sep=comma] {graphics/data/plate_48_rayleigh_res_hinfrelerr_osimaginput.csv};
	\addlegendentry{\lmorminrel{} \lmoraaa{}}

	\addplot+ table[x=r, y=strint_equi, col sep=comma] {graphics/data/plate_48_rayleigh_res_hinfrelerr_osimaginput.csv};
	\addlegendentry{\lmorequi{}}

	\addplot+ table[x=r, y=osinput, col sep=comma] {graphics/data/plate_48_rayleigh_res_hinfrelerr_sobt.csv};
	\addlegendentry{\lmorbt{}}

\end{axis}
\end{tikzpicture}%
  \tikzexternaldisable%

	\caption{Relative \Linf-errors of reduced models of the plate model 
		with proportional damping computed by several reduction methods with 
		\osimaginput{} projection.}
	\label{fig:rayleigh_hinf}
\end{figure}
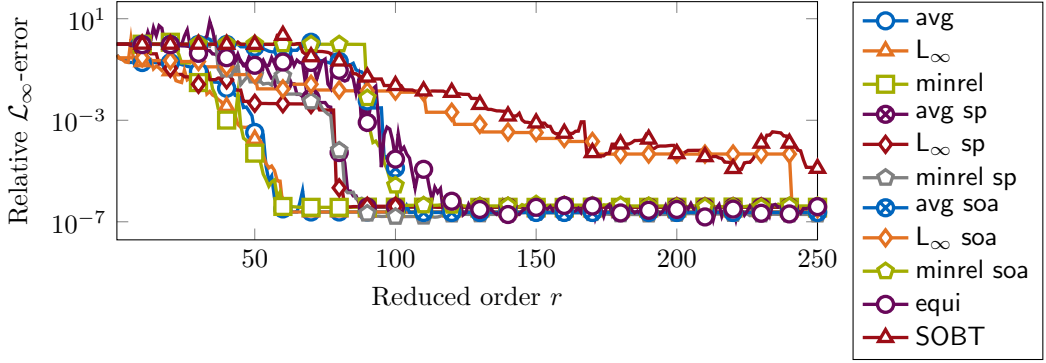

Now, we consider the results for the proportionally damped plate model.
The reduction methods yield models with on overall better accuracy compared to
the model with hysteretic damping, as there are less weakly damped poles in the 
transfer function.
Because of this, the accuracy for computing the \morscore{} is set to
$\epsilon = \num{1e-16}$; again, we consider $r_{\rm{max}} = 250$.
The \morscores{} for all employed methods are given in 
\Cref{fig:rayleigh_morscore}.
Especially \moravg{}, \morlinf{} and \morminrel{} using the standard
presampling method have high \morscores, while the models obtained from methods
considering \moraaa{} or \morsp{} presampling have slightly lower \morscores.
This is acceptable considering the lower computational cost for computing these
presampling bases, especially for \moraaa{}.
Only \morlinf{}~\moraaa{} has a considerably lower \morscore.
It can be seen in the error-per-order plot \Cref{fig:rayleigh_hinf} that the 
reduced model computed by \morlinf{}~\moraaa{} reaches the error level of the 
other methods for a reduced order $r=250$.
Its lower \morscore{} is mainly influenced by the fact that the reduced order
$r$ is again incremented in steps of $10$, i.e., the size of the employed Krylov
subspace.
We can also observe that the presampling methods have a comparable influence on
\moravg{}, \morlinf{} and \morminrel{}.
Using the standard presampling, the best achievable accuracy can be reached with
reduced models of order around $r = 60$, models computed from \morsp{}
presampling require $r = 90$, and using \moraaa{} yields comparable accuracy for
reduced models larger than $r = 100$.
\morlinf{}~\moraaa{} requires a larger reduced order $r$, because $r$ is a
multiple of 10 here.
However, \moravg{} and \morminrel{} in combination with \moraaa{} are comparable
to \morequi{}.
Thus, the presampling subspace computed by \moraaa{} is able to capture the most
important features of the original system's transfer function. 
The reduced models computed with the one-sided \morbt{} do not reach the 
accuracy of the other methods and attain their best approximation error
around $r \geq 160$.

\begin{figure}[t]
	\centering
  \tikzexternalenable%
  \tikzsetnextfilename{rayleigh_tf_red}%
  \filemodCmp{graphics/rayleigh_tf_red.tex}{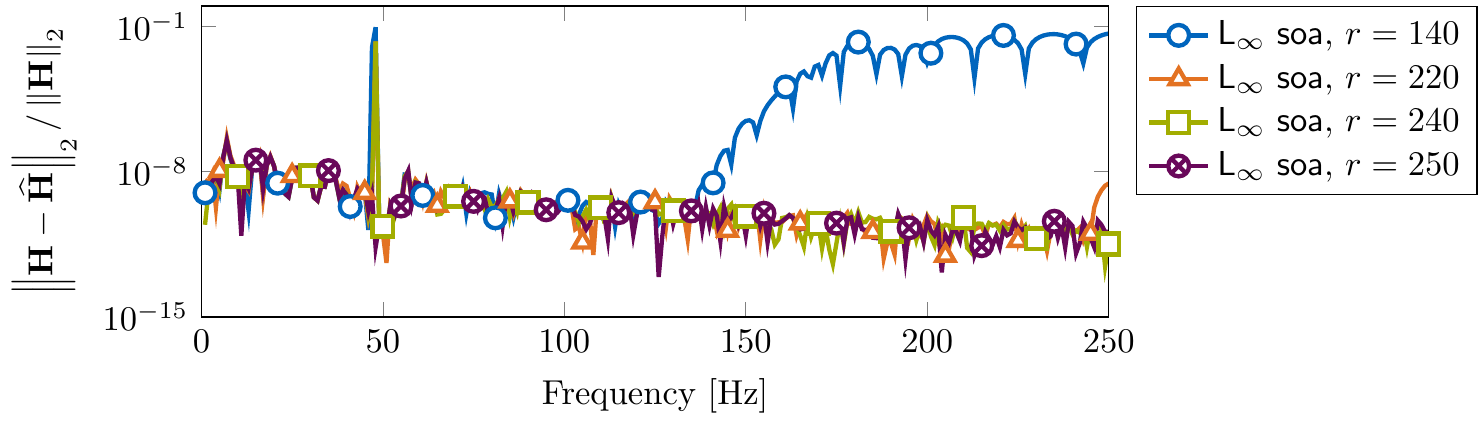}%
  {\tikzset{external/remake next}}{}%
  \begin{tikzpicture}[baseline=(current bounding box.north), font = \small]
	\begin{axis}[height=0.2\textheight, width=.72\textwidth, ymode=log,
			xlabel = {Frequency \ubraces{\si{\hertz}}},
			ylabel = {$\twonorm{\tf - \tfr} / \twonorm {\tf}$}, 
			xmin=0, xmax=250,
			xtick distance = 50,
			ymin=1e-15, ymax=1,
			axis on top,
			legend pos = outer north east,
			legend cell align = left,
			cycle list name=TUMcolorlist,
			mark repeat={20},
		]
		\addplot+ table[x=s, y=hinferr, col sep=comma] {graphics/data/plate_48_rayleigh_res_tf_strint_linf_aaaa_osimaginput_r_140.csv};
		\addlegendentry{\lmorlinf{} \lmoraaa{}, $r=140$}

		\addplot+[mark phase={5}] table[x=s, y=hinferr, col sep=comma] {graphics/data/plate_48_rayleigh_res_tf_strint_linf_aaaa_osimaginput_r_220.csv};
		\addlegendentry{\lmorlinf{} \lmoraaa{}, $r=220$}

		\addplot+[mark phase={10}] table[x=s, y=hinferr, col sep=comma] {graphics/data/plate_48_rayleigh_res_tf_strint_linf_aaaa_osimaginput_r_240.csv};
		\addlegendentry{\lmorlinf{} \lmoraaa{}, $r=240$}

		\addplot+[mark phase={15}] table[x=s, y=hinferr, col sep=comma] {graphics/data/plate_48_rayleigh_res_tf_strint_linf_aaaa_osimaginput.csv};
		\addlegendentry{\lmorlinf{} \lmoraaa{}, $r=250$}

	\end{axis}
\end{tikzpicture}%
  \tikzexternaldisable%

	\caption{Pointwise relative transfer function errors for reduced-order models
    of the proportionally damped plate computed by \morlinf{}~\moraaa{}.
    The error peak near the tuning frequency of the TVAs at $f = \SI{48}{\hertz}$
    is clearly visible for reduced-order models of size $r < 250$.}
	\label{fig:rayleigh_tf_red}
\end{figure}
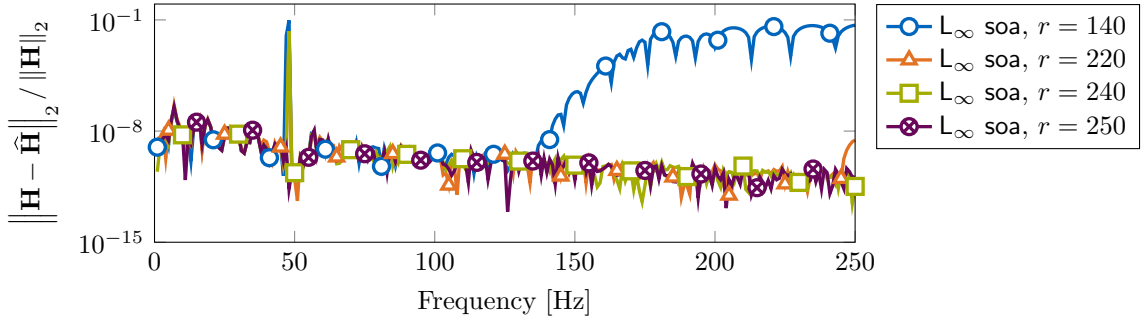

The reason for the stagnation of the approximation error of reduced models 
computed by \morlinf{}~\moraaa{} can be observed in the transfer function error
plot \Cref{fig:rayleigh_tf_red}.
The relatively high error in the frequency region near the tuning frequency of
the TVA at $f=\SI{48}{\hertz}$ is present up to models with $r=240$.
Only at $r=250$, \morlinf{} selects the shift and corresponding subspace
providing enough information to approximate the original transfer function also
in this frequency region.
Thus, the error drops to the level of the models computed using the other
reduction methods.

\begin{figure}[htb]
	\centering
  \tikzexternalenable%
  \tikzsetnextfilename{rayleigh_single_tf_red}%
  \filemodCmp{graphics/rayleigh_single_tf_red.tex}{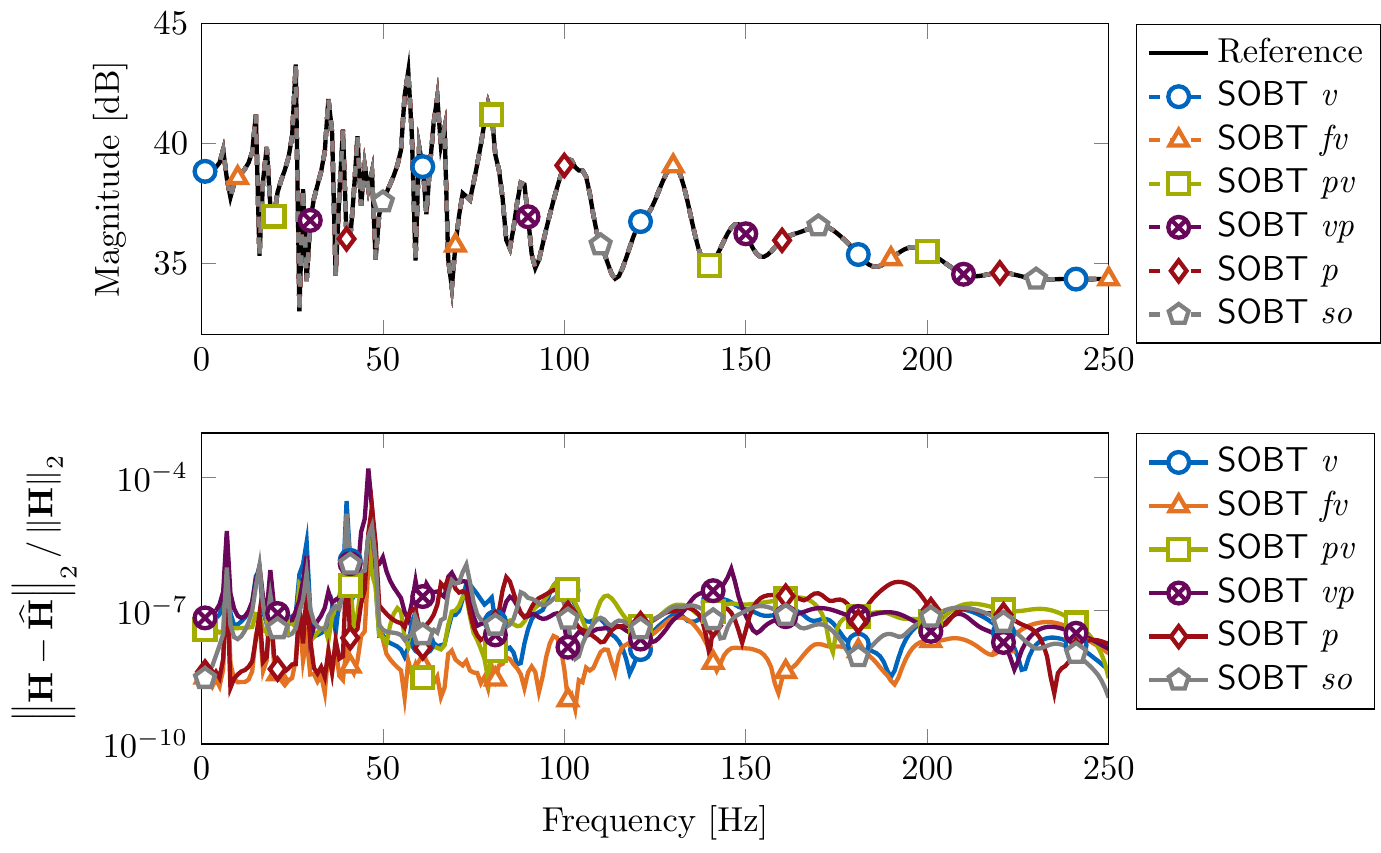}%
  {\tikzset{external/remake next}}{}%
  \begin{tikzpicture}[baseline=(current bounding box.north), font = \small]
\begin{axis}[name=rayligh_single_tf,height=0.2\textheight, width=.72\textwidth, 
	ylabel = {Magnitude \ubraces{\si{\decibel}}},
	xmin=0, xmax=250,
	xtick distance = 50,
	ymin=32, ymax=45,
	axis on top,
	legend cell align = left,
	legend pos = outer north east,
	cycle list name=TUMcolorlist,
	mark repeat={60}
	]
	\addplot[color=black, line width=\mylinewidth] table[x=s, y expr={10*log10(\thisrow{res}/1e-9)}, col sep=comma] {graphics/data/plate_48_rayleigh_single_res_tf_sobt_v.csv};
	\addlegendentry{Reference}
	\pgfplotsset{cycle list shift=-1}

	\addplot+[dashed] table[x=s, y expr={10*log10(\thisrow{res_r}/1e-9)}, col sep=comma] {graphics/data/plate_48_rayleigh_single_res_tf_sobt_v.csv};
	\addlegendentry{\lmorbt{} \soV{}}

	\addplot+[dashed, mark phase={10}] table[x=s, y expr={10*log10(\thisrow{res_r}/1e-9)}, col sep=comma] {graphics/data/plate_48_rayleigh_single_res_tf_sobt_fv.csv};
	\addlegendentry{\lmorbt{} \soFV{}}

	\addplot+[dashed, mark phase={20}] table[x=s, y expr={10*log10(\thisrow{res_r}/1e-9)}, col sep=comma] {graphics/data/plate_48_rayleigh_single_res_tf_sobt_pv.csv};
	\addlegendentry{\lmorbt{} \soPV{}}

	\addplot+[dashed, mark phase={30}] table[x=s, y expr={10*log10(\thisrow{res_r}/1e-9)}, col sep=comma] {graphics/data/plate_48_rayleigh_single_res_tf_sobt_vp.csv};
	\addlegendentry{\lmorbt{} \soVP{}}

	\addplot+[dashed, mark phase={40}] table[x=s, y expr={10*log10(\thisrow{res_r}/1e-9)}, col sep=comma] {graphics/data/plate_48_rayleigh_single_res_tf_sobt_p.csv};
	\addlegendentry{\lmorbt{} \soP{}}

	\addplot+[dashed, mark phase={50}] table[x=s, y expr={10*log10(\thisrow{res_r}/1e-9)}, col sep=comma] {graphics/data/plate_48_rayleigh_single_res_tf_sobt_so.csv};
	\addlegendentry{\lmorbt{} \soSO{}}

\end{axis}
\begin{axis}[at={($(rayligh_single_tf.south)-(0,1cm,)$)},anchor=north,height=0.2\textheight, width=.72\textwidth, ymode=log,
	xlabel = {Frequency \ubraces{\si{\hertz}}},
	ylabel = {$\twonorm{\tf - \tfr} / \twonorm {\tf}$},
	xmin=0, xmax=250,
	xtick distance = 50,
	ymin=1e-10, ymax=1e-3,
	axis on top,
	legend cell align = left,
	legend pos = outer north east,
	cycle list name=TUMcolorlist,
	mark repeat={20},
	]

	\addplot+ table[x=s, y=hinferr, col sep=comma] {graphics/data/plate_48_rayleigh_single_res_tf_sobt_v.csv};
	\addlegendentry{\lmorbt{} \soV{}}

	\addplot+ table[x=s, y=hinferr, col sep=comma] {graphics/data/plate_48_rayleigh_single_res_tf_sobt_fv.csv};
	\addlegendentry{\lmorbt{} \soFV{}}
	
%

	\addplot+ table[x=s, y=hinferr, col sep=comma] {graphics/data/plate_48_rayleigh_single_res_tf_sobt_pv.csv};
	\addlegendentry{\lmorbt{} \soPV{}}

	\addplot+ table[x=s, y=hinferr, col sep=comma] {graphics/data/plate_48_rayleigh_single_res_tf_sobt_vp.csv};
	\addlegendentry{\lmorbt{} \soVP{}}

	\addplot+ table[x=s, y=hinferr, col sep=comma] {graphics/data/plate_48_rayleigh_single_res_tf_sobt_p.csv};
	\addlegendentry{\lmorbt{} \soP{}}

	\addplot+ table[x=s, y=hinferr, col sep=comma] {graphics/data/plate_48_rayleigh_single_res_tf_sobt_so.csv};
	\addlegendentry{\lmorbt{} \soSO{}}

\end{axis}
\end{tikzpicture}%
  \tikzexternaldisable%

	\caption{Comparison of the \morbt{} formulas. Original and reduced transfer 
		functions as well as relative errors for the proportionally damped 
		plate with a single output and reduced order $r = 250$. Note the 
		approximation error maximum around the tuning frequency of the TVA 
		($f = \SI{48}{\hertz}$).}
	\label{fig:rayleigh_single_tf_red}
\end{figure}
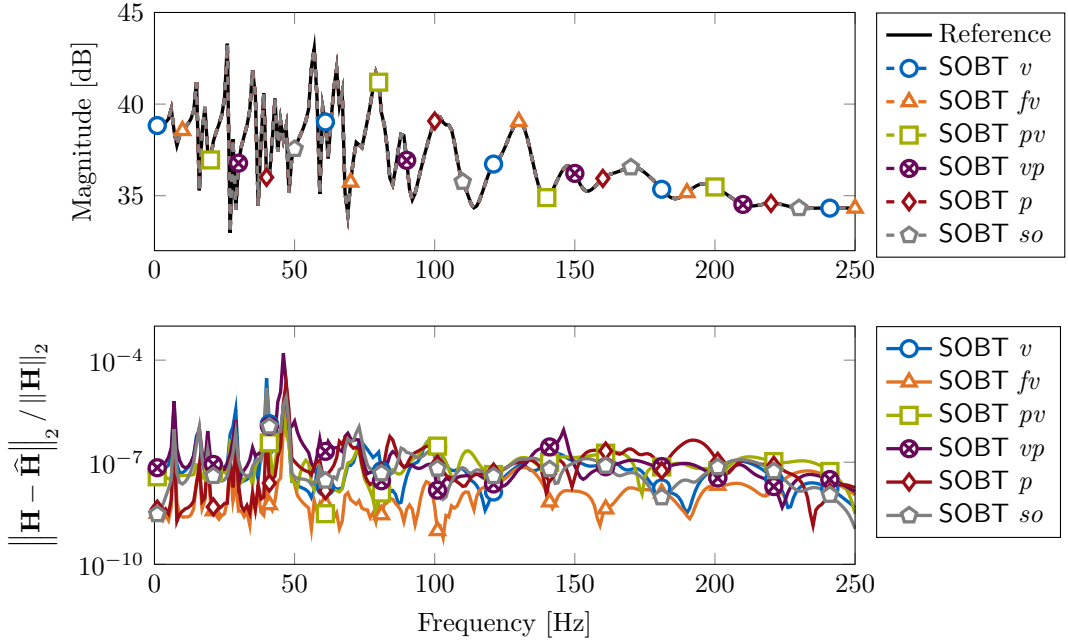

In order to compare the different formulas for \morbt{} given in
\Cref{tab:sobt}, a slightly modified model of the proportionally damped plate is
considered in the following.
Here, the displacement is evaluated at the load location rather than averaged
over the plate's surface.
The resulting SISO system allows the computation of a left projection basis
$\mv{W}$ in reasonable time.
The transfer function of the resulting systems, of the computed reduced-order
models using \morbt{} as well as the relative approximation errors are shown
in \Cref{fig:rayleigh_single_tf_red}.
All formulas, except \soVPM{} \soPM{} which fail at computing reduced-order
models approximating the original transfer function, yield reasonably accurate
reduced models.
Similar to the results above, the very local effect of the TVAs around
$f=\SI{48}{\hertz}$ cannot easily be captured by the reduction method,
resulting in a visible peak in the relative error in this frequency region even
for a reduced order of $r = 250$.

\begin{figure}[t]
	\centering
  \tikzexternalenable%
  \tikzsetnextfilename{rayleigh_single_morscore}%
  \filemodCmp{graphics/rayleigh_single_morscore.tex}{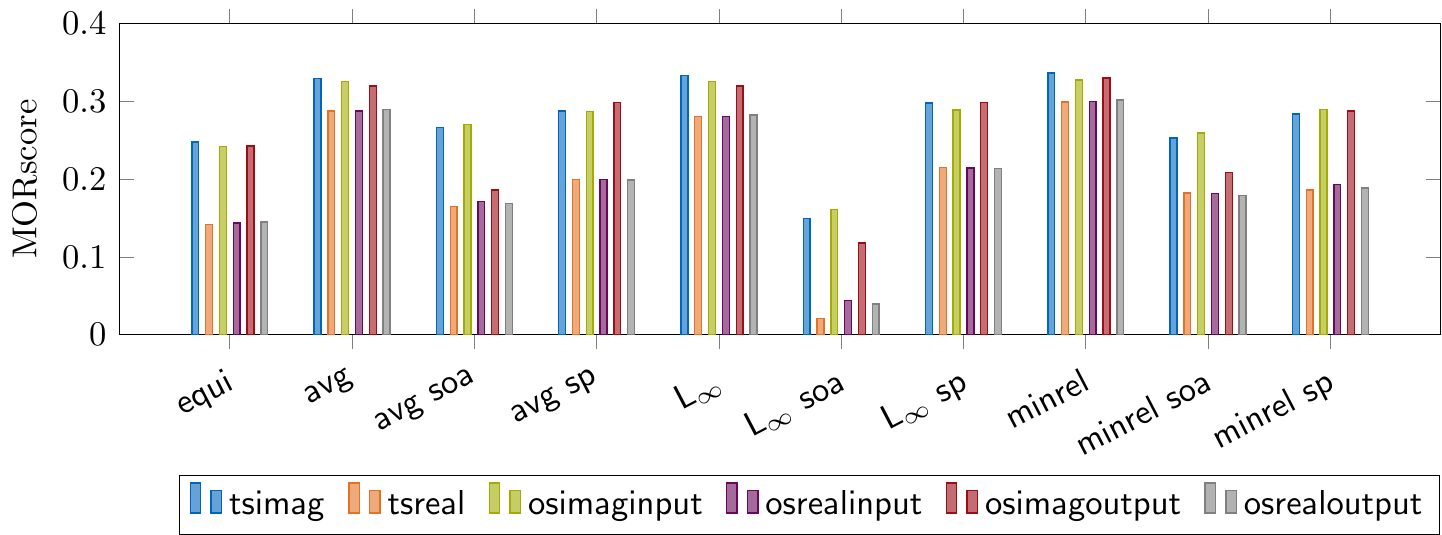}%
  {\tikzset{external/remake next}}{}%
  \begin{tikzpicture}[font = \small]
	\begin{axis}[height=0.2\textheight, width=\textwidth,
		ybar,
		ylabel = {\morscore{}},
		xtick distance = 1,
		ytick distance = .1,
		ymin=0, ymax=.4,
		legend cell align = left,
		legend style={at={(1,-0.45)},anchor=north east, /tikz/every even column/.append style={column sep=1ex}},
		legend columns=-1,
		no markers,
		every axis plot/.append style={fill},
		cycle list name=TUMcolorlist_simple,
		bar width=0.07cm,
		symbolic x coords={\lmorequi{}, \lmoravg{}, \lmoravg{} \lmoraaa{}, \lmoravg{} \lmorsp{}, \lmorlinf{}, \lmorlinf{} \lmoraaa{}, \lmorlinf{} \lmorsp{}, \lmorminrel{}, \lmorminrel{} \lmoraaa{}, \lmorminrel{} \lmorsp{}},
		xticklabel style={rotate=27.5,anchor=north east},
		]
		\addplot+[fill opacity=0.6] table [x=method, y index=1, col sep=comma] {graphics/data/plate_48_rayleigh_single_res_mor_score_tsimag.csv};
		\addlegendentry{\ltsimag{}}

		\addplot+[fill opacity=0.6] table [x=method, y index=1, col sep=comma] {graphics/data/plate_48_rayleigh_single_res_mor_score_tsreal.csv};
		\addlegendentry{\ltsreal{}}

		\addplot+[fill opacity=0.6] table [x=method, y index=1, col sep=comma] {graphics/data/plate_48_rayleigh_single_res_mor_score_osimaginput.csv};
		\addlegendentry{\losimaginput{}}

		\addplot+[fill opacity=0.6] table [x=method, y index=1, col sep=comma] {graphics/data/plate_48_rayleigh_single_res_mor_score_osrealinput.csv};
		\addlegendentry{\losrealinput{}}

		\addplot+[fill opacity=0.6] table [x=method, y index=1, col sep=comma] {graphics/data/plate_48_rayleigh_single_res_mor_score_osimagoutput.csv};
		\addlegendentry{\losimagoutput{}}

		\addplot+[fill opacity=0.6] table [x=method, y index=1, col sep=comma] {graphics/data/plate_48_rayleigh_single_res_mor_score_osrealoutput.csv};
		\addlegendentry{\losrealoutput{}}

	\end{axis}
\end{tikzpicture}%
  \tikzexternaldisable%

	\caption{\morscores{} of all employed reduction for the proportionally damped
    plate with a single output, using the maximum accuracy
    $\epsilon = \num{1e-16}$ and maximum order $r_{\rm{max}} = 250$.}
	\label{fig:rayleigh_single_morscore}
\end{figure}

The \morscores{} of all other methods employed to compute reduced models for the
single output version of the example are shown in
\Cref{fig:rayleigh_single_morscore}.
The results are similar to the ones reported above and all methods produce
accurate reduced models.
Again, \morlinf{}~\moraaa{} has a lower \morscore{}, as the $r$ 
is incremented in steps of $k = 10$.
All reduced-order models capture the transfer function in the critical region
around $f = \SI{48}{\hertz}$ for large enough reduced orders.
As the displacement at the loading point is evaluated in the transfer function,
i.e., input and output vectors are identical such that a two-sided projection is
not beneficial and \tsimag{}, \osimaginput{} and \osimagoutput{} show nearly the
same \morscores{}.
The same holds for the real-valued projections.
We note that \moravg{} and \morminrel{} with classical presampling yield nearly
identical results for complex- and real-valued truncation matrices.

\subsection{Sound transmission through a plate}

\begin{figure}[t]
	\centering
  \begin{subfigure}[b]{.4\textwidth}
    \centering
	  \def\svgwidth{\textwidth}
\begingroup%
  \makeatletter%
  \providecommand\color[2][]{%
    \errmessage{(Inkscape) Color is used for the text in Inkscape, but the package 'color.sty' is not loaded}%
    \renewcommand\color[2][]{}%
  }%
  \providecommand\transparent[1]{%
    \errmessage{(Inkscape) Transparency is used (non-zero) for the text in Inkscape, but the package 'transparent.sty' is not loaded}%
    \renewcommand\transparent[1]{}%
  }%
  \providecommand\rotatebox[2]{#2}%
  \newcommand*\fsize{\dimexpr\f@size pt\relax}%
  \newcommand*\lineheight[1]{\fontsize{\fsize}{#1\fsize}\selectfont}%
  \ifx\svgwidth\undefined%
    \setlength{\unitlength}{467.58313006bp}%
    \ifx\svgscale\undefined%
      \relax%
    \else%
      \setlength{\unitlength}{\unitlength * \real{\svgscale}}%
    \fi%
  \else%
    \setlength{\unitlength}{\svgwidth}%
  \fi%
  \global\let\svgwidth\undefined%
  \global\let\svgscale\undefined%
  \makeatother%
  \begin{picture}(1,0.68920016)%
    \lineheight{1}%
    \setlength\tabcolsep{0pt}%
    \put(0,0){\includegraphics[width=\unitlength,page=1]{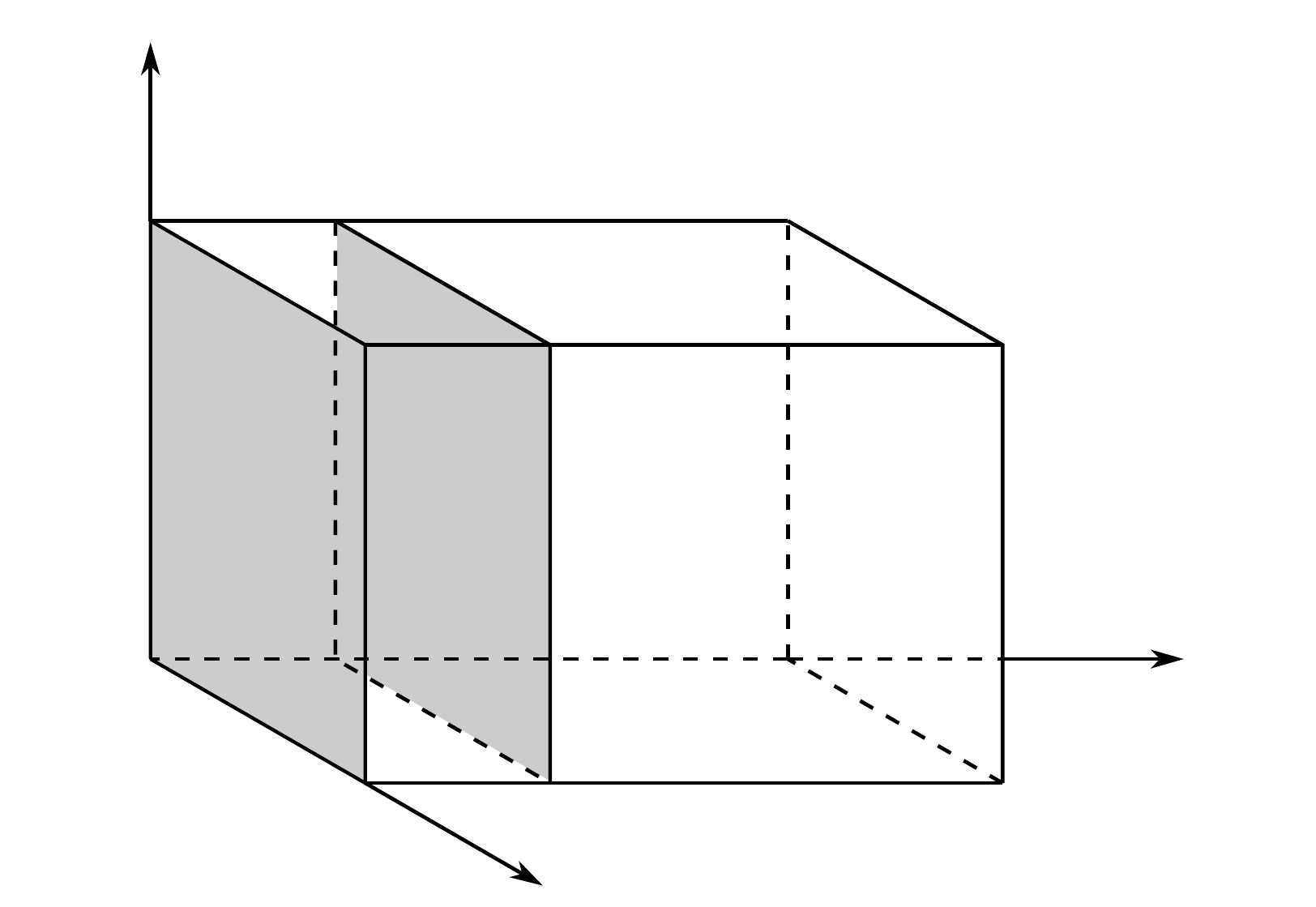}}%
    \put(0.13455111,0.69949799){\color[rgb]{0,0,0}\makebox(0,0)[lt]{\begin{minipage}{0.11915375\unitlength}\raggedright $y$\end{minipage}}}%
    \put(0.89758747,0.24876107){\color[rgb]{0,0,0}\makebox(0,0)[lt]{\begin{minipage}{0.11915375\unitlength}\raggedright $x$\end{minipage}}}%
    \put(0.42326391,0.06223964){\color[rgb]{0,0,0}\makebox(0,0)[lt]{\begin{minipage}{0.11915375\unitlength}\raggedright $z$\end{minipage}}}%
    \put(0,0){\includegraphics[width=\unitlength,page=2]{graphics/guy_model_svg-tex.pdf}}%
    \put(0.0417423,0.33958993){\color[rgb]{0.89019608,0.44705882,0.1254902}\makebox(0,0)[lt]{\begin{minipage}{0.21799067\unitlength}\raggedright $p$\end{minipage}}}%
    \put(0,0){\includegraphics[width=\unitlength,page=3]{graphics/guy_model_svg-tex.pdf}}%
    \put(0.69454874,0.32695388){\color[rgb]{0,0.39607843,0.74117647}\makebox(0,0)[lt]{\lineheight{0}\smash{\begin{tabular}[t]{l}$P_1$\end{tabular}}}}%
  \end{picture}%
\endgroup%

	  \vspace{.05cm}
  	\subcaption{Sketch of the problem with pressure load $p$ and evaluation
      point $P_{1}$.}
	  \label{fig:guy_sketch}
  \end{subfigure}%
  \hfill%
  \begin{subfigure}[b]{.55\textwidth}
    \centering
  \tikzexternalenable%
  \tikzsetnextfilename{guy_damped}%
  \filemodCmp{graphics/guy_damped.tex}{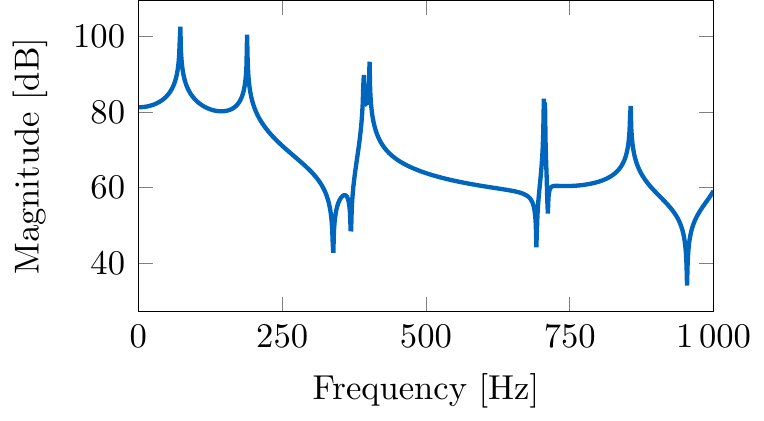}%
  {\tikzset{external/remake next}}{}%
  \begin{tikzpicture}[font = \small]
\pgfkeys{/pgf/number format/.cd,1000 sep={\,}}
\begin{axis}[
    height = .2\textheight,
    width = .9\textwidth,
    xlabel = {Frequency \ubraces{\si{\hertz}}},
    ylabel = {Magnitude \ubraces{\si{\decibel}}},
	xmin = 0,
    xmax = 1000,
	xtick distance = 250,
	axis on top,
	legend pos = outer north east, legend cell align = left,
	no markers
  ]
	
    \addplot[mark=none, color=TUMBlau, line width=\mylinewidth] table[x=freq, y=res, col sep=comma] {graphics/data/guy_damped.csv};
  \end{axis}
\end{tikzpicture}%
  \tikzexternaldisable%

    \subcaption{Transfer function.}
    \label{fig:guy_tf}
  \end{subfigure}

	\caption{Sketch and transfer function of the sound transmission problem. }
\end{figure}

Radiation of vibrating plates and excitation of a structure by an oscillating 
acoustic fluid are modeled in this example.
The system consists of a cuboid acoustic cavity, where one wall is considered a
system of two parallel elastic brass plates with a $\SI{2}{\centi\meter}$ air
gap between them; all other walls are considered rigid.
The plates measure $0.2 \times \SI{0.2}{\meter} $ and have a thickness of 
$t = \SI{0.9144}{\milli\meter}$.
The material parameters $E = \SI{104}{\giga\pascal}$,
$\rho = \SI{8500}{\kilogram\per\cubic\meter}$
$\nu = \num{0.37} $ are considered for brass.
The receiving cavity is $\SI{0.2}{\meter}$ wide; wave speed
$c = \SI{343}{\meter\per\second}$ and density
$\rho = \SI{1.21}{\kilogram\per\cubic\meter}$ are considered for the 
acoustic fluid.
The configuration is based on an experiment conducted in~\cite{Guy81}.
It is sketched in \Cref{fig:guy_sketch} along with the acoustic 
pressure in the cavity resulting from a uniform pressure load $p$ applied to 
the outer plate.
The pressure is measured at the middle point of the wall 
opposite to the elastic plate $P_{1}$.
Energy dissipation inside the structural part of the system is modeled using
proportional damping with $\beta = \num{1e-7}$.
The system is discretized using the finite element method and $n = \num{95480}$
degrees of freedom are required to obtain an accurate result in a frequency
range up to \SI{1000}{\hertz}.
No acoustic sources are present, so the excitation vector is frequency
independent.
Considering the two way coupling between structure and fluid leads to
non-symmetric system matrices.
Thus, a transfer function of Case~A with real-valued matrices is used 
to describe the system.

The standard presampling for \morminrel, \moravg{} and \morlinf{} considers 
$n_{\rm{s}} = 200$ frequency shifts distributed linearly in $2\pi\i [1, 1000]$.
As the quadratic frequency associated with the mass matrix is 
the highest order of $s$ in the transfer function, each shift computed by a 
\morsp{} presampling contributes three columns to the intermediate basis.
Therefor, $n_{\rm{s}} = 67$ shifts, linearly distributed in the same range, are
chosen such that the intermediate reduction basis is of size $q = 201$.
For \moraaa{}, a local order $k = 10$ along with $n_{\rm{s}} = 20$ is chosen,
yielding an intermediate reduction basis of order $q = 200$.
Because the numerical model contains unstable eigenvalues, the required Gramians
for \morbt{} cannot be computed and, thus, the method is not applied.

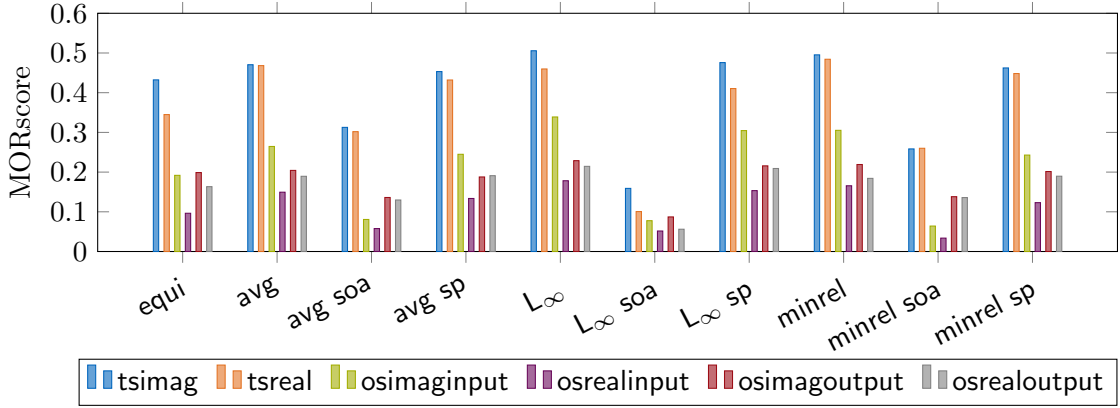
\begin{figure}[t]
	\centering
  \tikzexternalenable%
  \tikzsetnextfilename{guy_morscore}%
  \filemodCmp{graphics/guy_morscore.tex}{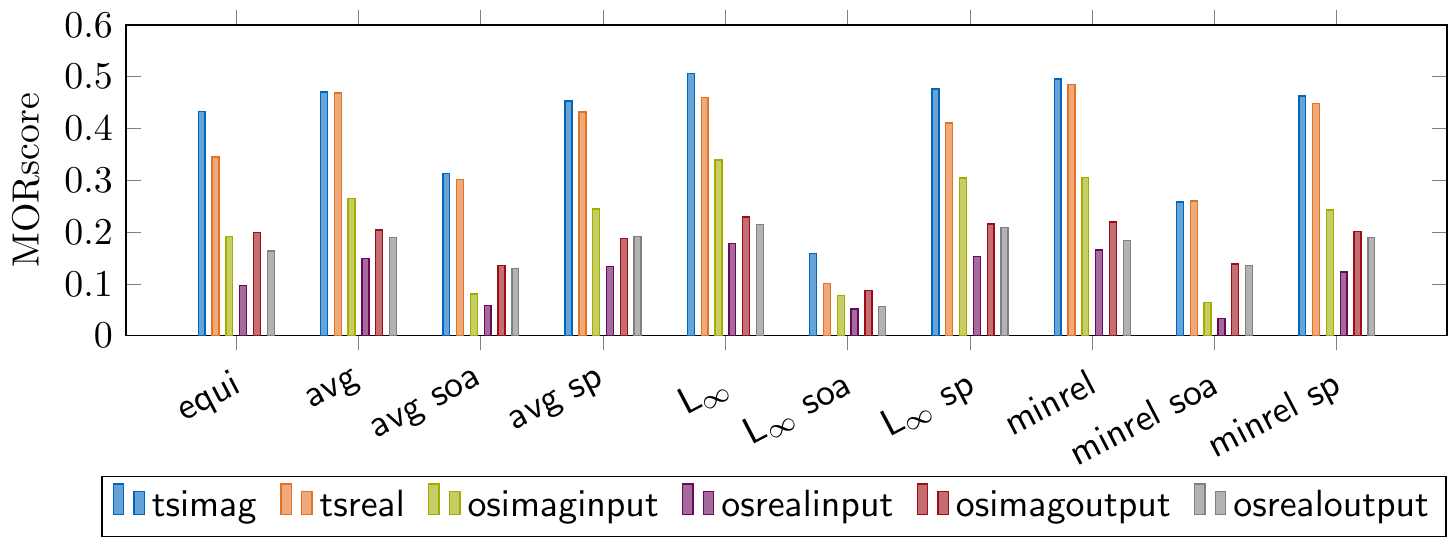}%
  {\tikzset{external/remake next}}{}%
  \begin{tikzpicture}
	\begin{axis}[height=0.2\textheight, width=\textwidth,
		ybar,
		ylabel = {\morscore{}},
		xtick distance = 1,
		ytick distance = .1,
		ymin=0, ymax=.6,
		legend cell align = left,
		legend style={at={(1,-0.45)},anchor=north east, /tikz/every even column/.append style={column sep=1ex}},
		legend columns=-1,
		no markers,
		every axis plot/.append style={fill},
		cycle list name=TUMcolorlist_simple,
		bar width=0.07cm,
		symbolic x coords={\lmorequi{}, \lmoravg{}, \lmoravg{} \lmoraaa{}, \lmoravg{} \lmorsp{}, \lmorlinf{}, \lmorlinf{} \lmoraaa{}, \lmorlinf{} \lmorsp{}, \lmorminrel{}, \lmorminrel{} \lmoraaa{}, \lmorminrel{} \lmorsp{}},
		xticklabel style={rotate=27.5,anchor=north east},
		]
		\addplot+[fill opacity=0.6] table [x=method, y index=1, col sep=comma] {graphics/data/guy_damped_mid_res_mor_score_tsimag.csv};
		\addlegendentry{\ltsimag{}}

		\addplot+[fill opacity=0.6] table [x=method, y index=1, col sep=comma] {graphics/data/guy_damped_mid_res_mor_score_tsreal.csv};
		\addlegendentry{\ltsreal{}}

		\addplot+[fill opacity=0.6] table [x=method, y index=1, col sep=comma] {graphics/data/guy_damped_mid_res_mor_score_osimaginput.csv};
		\addlegendentry{\losimaginput{}}

		\addplot+[fill opacity=0.6] table [x=method, y index=1, col sep=comma] {graphics/data/guy_damped_mid_res_mor_score_osrealinput.csv};
		\addlegendentry{\losrealinput{}}

		\addplot+[fill opacity=0.6] table [x=method, y index=1, col sep=comma] {graphics/data/guy_damped_mid_res_mor_score_osimagoutput.csv};
		\addlegendentry{\losimagoutput{}}

		\addplot+[fill opacity=0.6] table [x=method, y index=1, col sep=comma] {graphics/data/guy_damped_mid_res_mor_score_osrealoutput.csv};
		\addlegendentry{\losrealoutput{}}

	\end{axis}
\end{tikzpicture}%
  \tikzexternaldisable%

	\caption{\morscores{} of all employed reduction and projection methods for 
		the sound transmission problem, with maximum accuracy 
		$\epsilon = \num{1e-16}$ and maximum order $r_{\rm{max}} = 100$.}
	\label{fig:guy_morscore}
\end{figure}

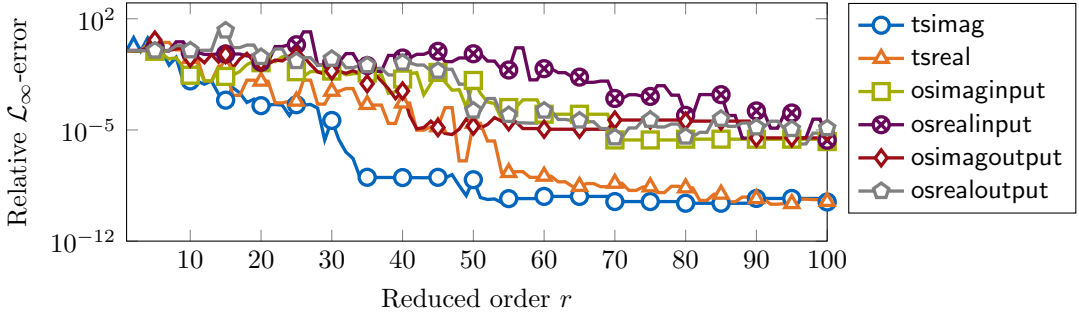
\begin{figure}[t]
	\centering
  \tikzexternalenable%
  \tikzsetnextfilename{guy_err_hinf}%
  \filemodCmp{graphics/guy_err_hinf.tex}{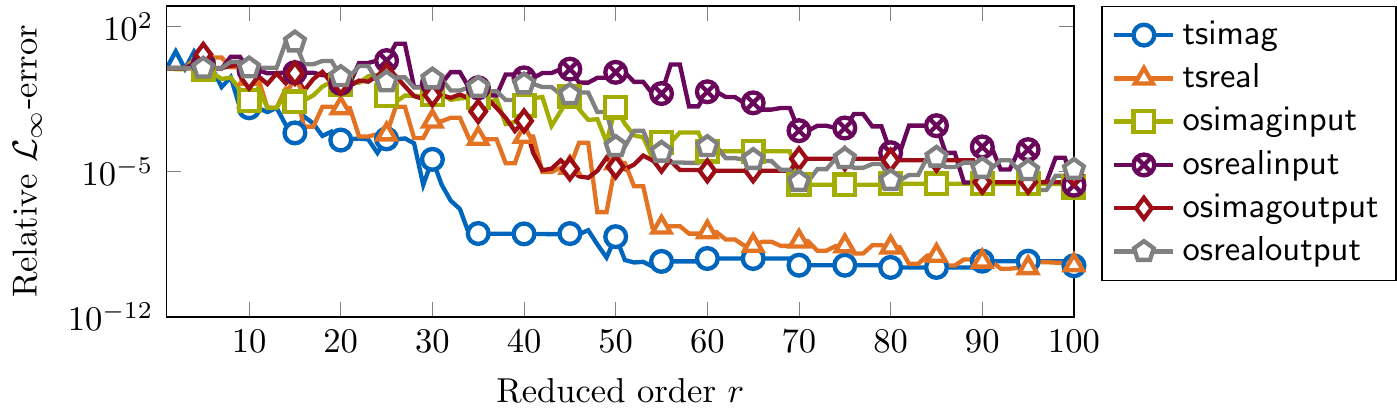}%
  {\tikzset{external/remake next}}{}%
  \begin{tikzpicture}[font = \small]
\begin{axis}[height=0.2\textheight, width=.72\textwidth, ymode=log,
	xlabel = {Reduced order $r$},
	ylabel = {\errlabel},
	xmin=1, xmax=100,
	ymin=1e-12, ymax=1e3,
	axis on top,
	legend pos = outer north east, legend cell align = left,
	cycle list name=TUMcolorlist,
	mark repeat={5},
	mark phase=5,
	]
	\addplot+ table[x=r, y=tsimag, col sep=comma] {graphics/data/guy_damped_mid_res_hinfrelerr_strint_equi.csv};
	\addlegendentry{\ltsimag{}}

	\addplot+ table[x=r, y=tsreal, col sep=comma] {graphics/data/guy_damped_mid_res_hinfrelerr_strint_equi.csv};
	\addlegendentry{\ltsreal{}}

	\addplot+ table[x=r, y=osimaginput, col sep=comma] {graphics/data/guy_damped_mid_res_hinfrelerr_strint_equi.csv};
	\addlegendentry{\losimaginput{}}

	\addplot+ table[x=r, y=osrealinput, col sep=comma] {graphics/data/guy_damped_mid_res_hinfrelerr_strint_equi.csv};
	\addlegendentry{\losrealinput{}}

	\addplot+ table[x=r, y=osimagoutput, col sep=comma] {graphics/data/guy_damped_mid_res_hinfrelerr_strint_equi.csv};
	\addlegendentry{\losimagoutput{}}

	\addplot+ table[x=r, y=osrealoutput, col sep=comma] {graphics/data/guy_damped_mid_res_hinfrelerr_strint_equi.csv};
	\addlegendentry{\losrealoutput{}}

\end{axis}
\end{tikzpicture}%
  \tikzexternaldisable%

	\caption{Comparison of the relative \Linf-errors of reduced models of 
		the sound transmission problem computed with \morequi{} using different 
		projections.}
	\label{fig:guy_hinf}
\end{figure}

The \morscores{} given in \Cref{fig:guy_morscore} show that especially the 
two-sided projections yield very good results with the highest \morscores{} 
observed.
As expected, \morlinf{}~\moraaa{} falls short due to the reduced order being
again incremented in steps of $10$.
But also \moravg{}~\moraaa{} and \morminrel{}~\moraaa{} perform not as good as
the other methods, while still showing a \morscore{} larger than $0.3$, which is
comparable to the other numerical examples.
It can be seen that using one-sided projections has a significant impact on the
approximation quality.
The error comparisons in \Cref{fig:guy_hinf} show that the approximation error
of the one-sided projections stagnates at around \num{1e-5}, while the two-sided
projections yield models with higher accuracy. 
As expected, the real-valued projection yield reduced-order models of comparable
accuracy for higher $r$.

\subsection{Radiation and scattering of a complex geometry}

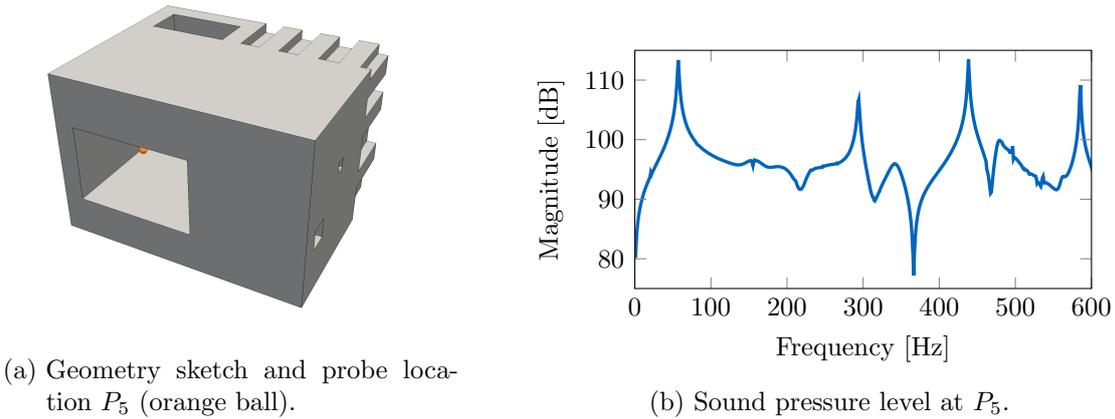
\begin{figure}[t]
  \centering
  \begin{subfigure}[b]{.4\textwidth}
    \centering
	  \includegraphics[width=.95\textwidth]{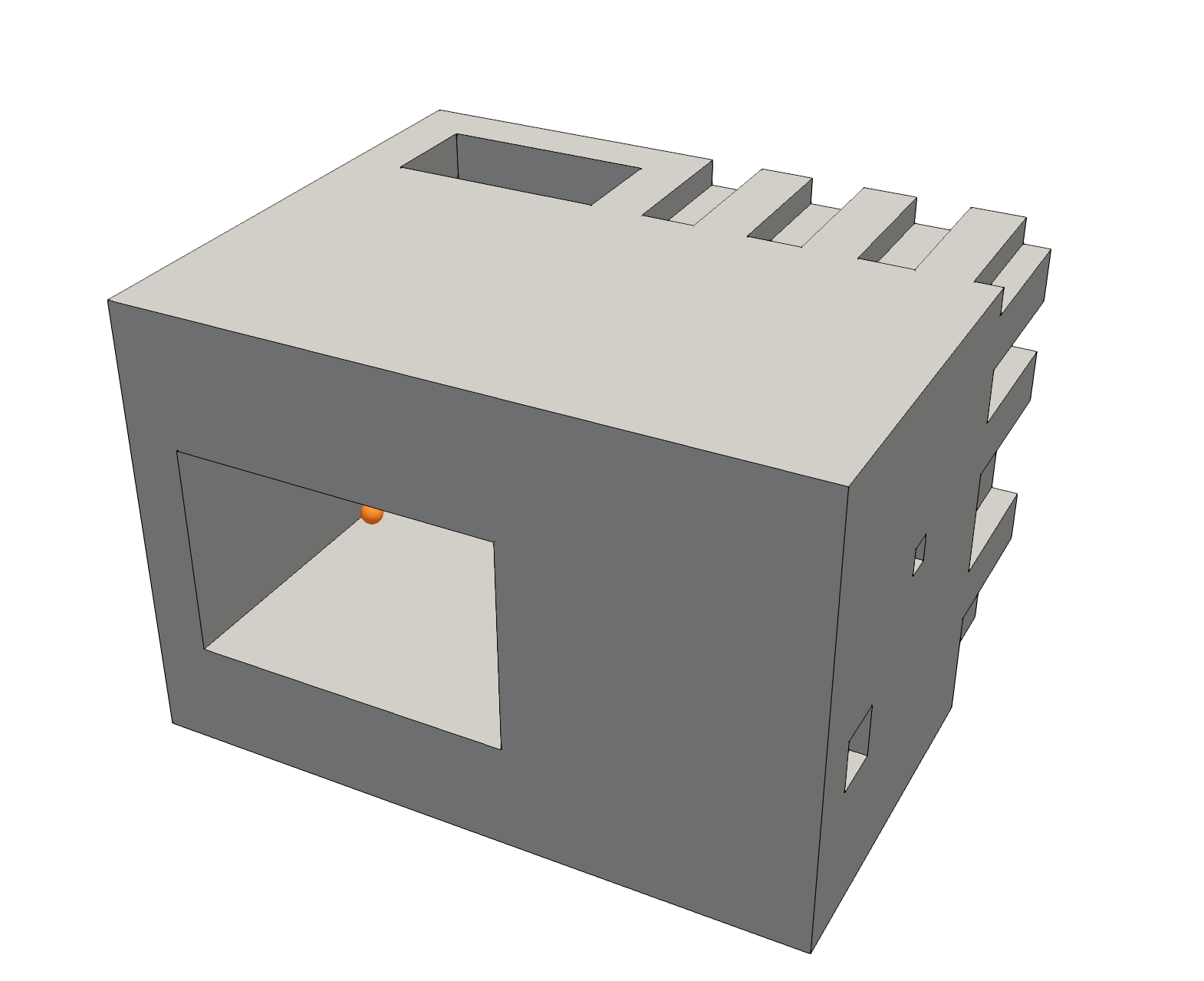}
	  \vspace{.3cm}
    \subcaption{Geometry sketch and probe location $P_5$ (orange ball).}
    \label{fig:radiatterer_sketch}
  \end{subfigure}%
  \hfill%
  \begin{subfigure}[b]{.55\textwidth}
    \centering
  \tikzexternalenable%
  \tikzsetnextfilename{radiatterer_p5}%
  \filemodCmp{graphics/radiatterer_p5.tex}{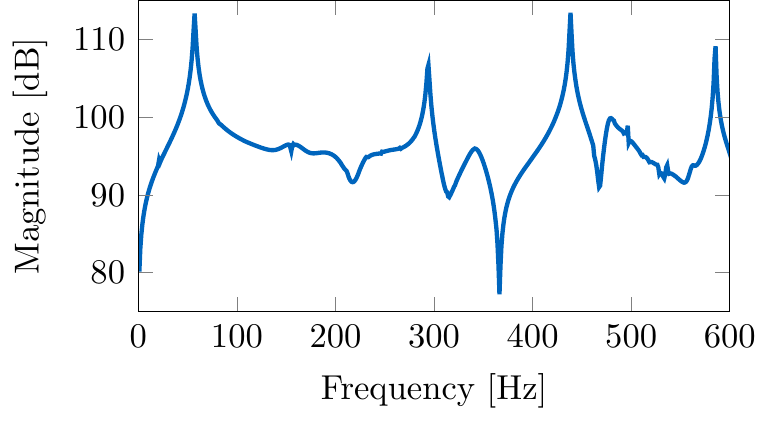}%
  {\tikzset{external/remake next}}{}%
  \begin{tikzpicture}[font = \small]
\begin{axis}[height=0.2\textheight, width=.92\textwidth, 
	xlabel = {Frequency \ubraces{\si{\hertz}}},
	ylabel = {Magnitude \ubraces{\si{\decibel}}},
	xmin=0, xmax=600,
	xtick distance = 100,
	ytick distance = 10,
	ymin=75, ymax=115,
	axis on top,
	legend pos = outer north east, legend cell align = left,
	no markers
	]
	\addplot[mark=none, color=TUMBlau, line width=\mylinewidth] table[x index={0}, y index={1}, col sep=comma] {graphics/data/radiatterer_p5.csv};
\end{axis}
\end{tikzpicture}%
  \tikzexternaldisable%

    \subcaption{Sound pressure level at $P_5$.}
    \label{fig:radiatterer_tf}
  \end{subfigure}

  \caption{Sketch and transfer function of the radiation problem. }
\end{figure}

Now, we consider a complex geometry based on a rigid block with various 
openings, cavities and sharp corners.
The experiment introduced in~\cite{HorKM15} is called \enquote{radiatterer} as
both radiation and scattering effects are taken into account.
The basic shape is a box with dimensions $2.5\times 2.0\times \SI{1.7}{\meter}$,
which is enclosed by an acoustic fluid of size $3.5\times 3.0\times
\SI{2.7}{\meter}$.
The geometry is sketched in \Cref{fig:radiatterer_sketch}; for the exact shape,
see~\cite{HorKM15}.
A normal velocity $v_{\rm{n}} = \SI{0.001}{\meter\per\second}$ 
acts on the complete surface of the geometry and excites the surrounding
acoustic fluid.
The free radiation from the geometry is realized with a PML of thickness 
$d = \SI{0.3}{\meter}$.
It is tuned to $f = \SI{500}{\hertz}$ to eliminate the frequency dependency of 
the PML matrices~\cite{VerMR19}.
Thus, this system is described by a transfer function of Case~B with complex
symmetric matrices and a frequency-dependent excitation vector.
The numerical model has an order of $n = \num{250000}$ and is evaluated in a 
frequency range from 1 to \SI{600}{\hertz}.
The transfer function plotted in \Cref{fig:radiatterer_tf} measures the sound
pressure level at a point inside the large cutout at
$(x, y, z) = (0.6,0.5,0.8)\si{\meter}$.
A reference solution is available in~\cite{Mar16}, where the same problem has
been analyzed using a boundary element method.

The standard presampling considers $n_{\rm{s}} = 200$ frequency shifts distributed 
linearly in $2\pi\i \left[1, 600 \right]$.
Again, \morsp{} yields three columns for the intermediate reduction basis, i.e.,
$n_{\rm{s}} = 67$ linearly distributed shifts are chosen to obtain a basis of
size $q = 201$.
For \moraaa{}, a local order $k = 5$ is used such that $n_{\rm{s}} = 40$ expansion
points yield a presampling basis with order $q = 200$.
A lower local order is chosen for this model as many weakly damped modes are
present in the transfer function and otherwise not enough information about the
full-order model would be available in the intermediate reduction basis.
\morbt{} is not applicable to this problem because of its frequency-dependent
input vector.

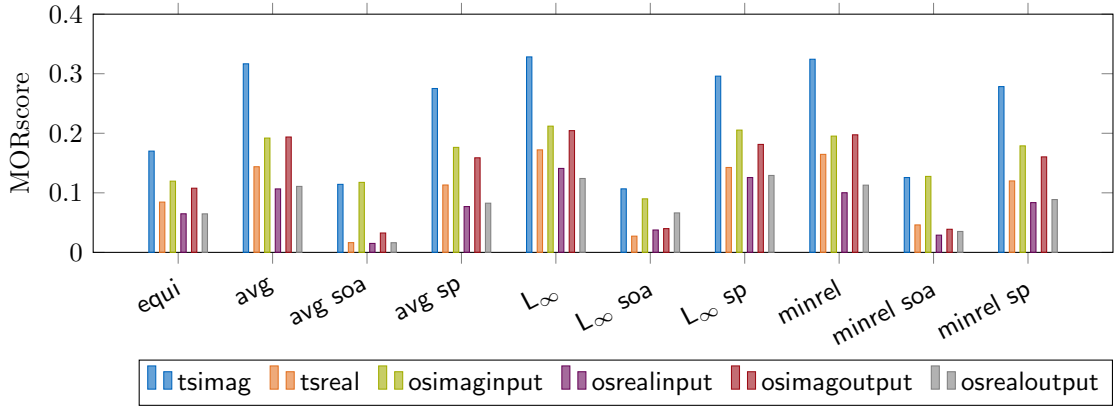
\begin{figure}[t]
	\centering
  \tikzexternalenable%
  \tikzsetnextfilename{radiatterer_morscore}%
  \filemodCmp{graphics/radiatterer_morscore.tex}{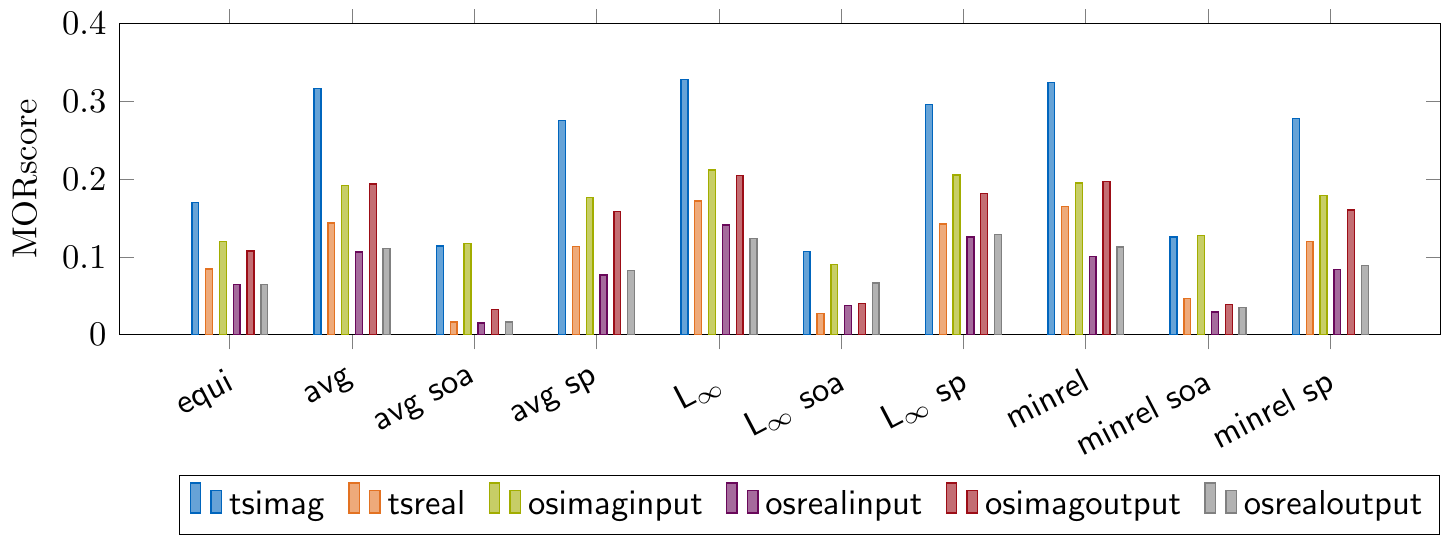}%
  {\tikzset{external/remake next}}{}%
  \begin{tikzpicture}[font = \small]
	\begin{axis}[height=0.2\textheight, width=\textwidth,
		ybar,
		ylabel = {\morscore{}},
		xtick distance = 1,
		ytick distance = .1,
		ymin=0, ymax=.4,
		legend cell align = left,
		legend style={at={(1,-0.45)},anchor=north east, /tikz/every even column/.append style={column sep=1ex}},
		legend columns=-1,
		no markers,
		every axis plot/.append style={fill},
		cycle list name=TUMcolorlist_simple,
		bar width=0.07cm,
		symbolic x coords={\lmorequi{}, \lmoravg{}, \lmoravg{} \lmoraaa{}, \lmoravg{} \lmorsp{}, \lmorlinf{}, \lmorlinf{} \lmoraaa{}, \lmorlinf{} \lmorsp{}, \lmorminrel{}, \lmorminrel{} \lmoraaa{}, \lmorminrel{} \lmorsp{}},
		xticklabel style={rotate=27.5,anchor=north east},
		]
		\addplot+[fill opacity=0.6] table [x=method, y index=1, col sep=comma] {graphics/data/radiatterer_new_fine_res_mor_score_tsimag.csv};
		\addlegendentry{\ltsimag{}}

		\addplot+[fill opacity=0.6] table [x=method, y index=1, col sep=comma] {graphics/data/radiatterer_new_fine_res_mor_score_tsreal.csv};
		\addlegendentry{\ltsreal{}}

		\addplot+[fill opacity=0.6] table [x=method, y index=1, col sep=comma] {graphics/data/radiatterer_new_fine_res_mor_score_osimaginput.csv};
		\addlegendentry{\losimaginput{}}

		\addplot+[fill opacity=0.6] table [x=method, y index=1, col sep=comma] {graphics/data/radiatterer_new_fine_res_mor_score_osrealinput.csv};
		\addlegendentry{\losrealinput{}}

		\addplot+[fill opacity=0.6] table [x=method, y index=1, col sep=comma] {graphics/data/radiatterer_new_fine_res_mor_score_osimagoutput.csv};
		\addlegendentry{\losimagoutput{}}

		\addplot+[fill opacity=0.6] table [x=method, y index=1, col sep=comma] {graphics/data/radiatterer_new_fine_res_mor_score_osrealoutput.csv};
		\addlegendentry{\losrealoutput{}}

	\end{axis}
\end{tikzpicture}%
  \tikzexternaldisable%

	\caption{\morscores{} of all employed reduction for the scattering problem,
    with the maximum accuracy $\epsilon=\num{1e-16}$ and maximum order
    $r_{\rm{max}} = 200$.}
	\label{fig:radiatterer_morscore}
\end{figure}

\begin{figure}[t]
	\centering
  \tikzexternalenable%
  \tikzsetnextfilename{radiatterer_new_err_hinf}%
  \filemodCmp{graphics/radiatterer_new_err_hinf.tex}{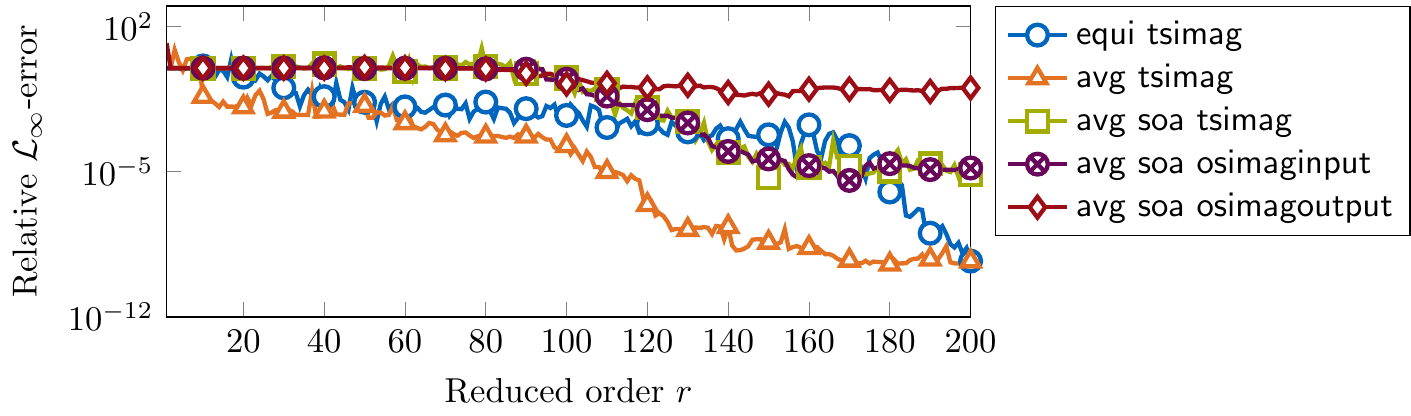}%
  {\tikzset{external/remake next}}{}%
  \begin{tikzpicture}[font = \small]
\begin{axis}[height=0.2\textheight, width=.65\textwidth, ymode=log,
	xlabel = {Reduced order $r$},
	ylabel = {\errlabel},
	xmin=1, xmax=200,
	ymin=1e-12, ymax=1e3,
	axis on top,
	legend pos = outer north east, legend cell align = left,
	cycle list name=TUMcolorlist,
	mark repeat={10},
	mark phase=10,
	]
	
	\addplot+ table[x=r, y=strint_equi, col sep=comma] {graphics/data/radiatterer_new_fine_res_hinfrelerr_tsimag.csv};
	\addlegendentry{\lmorequi{} \ltsimag{}}
	
	\addplot+ table[x=r, y=strint_avg, col sep=comma] {graphics/data/radiatterer_new_fine_res_hinfrelerr_tsimag.csv};
	\addlegendentry{\lmoravg{} \ltsimag{}}

	\addplot+ table[x=r, y=strint_avg_aaaa, col sep=comma] {graphics/data/radiatterer_new_fine_res_hinfrelerr_tsimag.csv};
	\addlegendentry{\lmoravg{} \lmoraaa{} \ltsimag{}}

	\addplot+ table[x=r, y=osimaginput, col sep=comma] {graphics/data/radiatterer_new_fine_res_hinfrelerr_strint_avg_aaaa.csv};
	\addlegendentry{\lmoravg{} \lmoraaa{} \losimaginput{}}

	\addplot+ table[x=r, y=osimagoutput, col sep=comma] {graphics/data/radiatterer_new_fine_res_hinfrelerr_strint_avg_aaaa.csv};
	\addlegendentry{\lmoravg{} \lmoraaa{} \losimagoutput{}}

\end{axis}
\end{tikzpicture}%
  \tikzexternaldisable%

	\caption{Relative \Linf-errors of reduced models of the radiation 
		problem computed by several reduction methods.}
	\label{fig:radiatterer_hinf}
\end{figure}
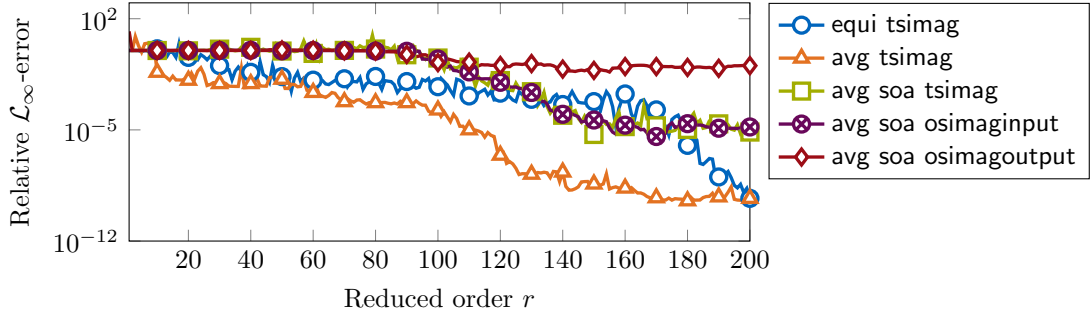

The \morscores{} for each applied reduction method are given in 
\Cref{fig:radiatterer_morscore}.
It can be seen that again two-sided interpolation with complex-valued bases
outperforms the other projection methods.
The lower rate of approximation when using presampling based on \moraaa{} is
also in line with the observations from the previous experiments. 
\morequi{} shows also a considerably worse performance than the methods based 
on presampling.
It can be seen in the error-per-order plot \Cref{fig:radiatterer_hinf} that the
error of reduced-order models computed with \morequi{} stagnates until
approximately $r = 140$ before it drops to the same level as the other methods.
This suggests that important features of the system response have not been
captured by the smaller reduction bases.
The oscillating behavior of the relative error in the region of $150 < r < 190$
is a sign that crucial parts of the transfer function are missed by sampling
with equidistantly distributed expansion points. 

\begin{figure}[t]
	\centering
  \tikzexternalenable%
  \tikzsetnextfilename{radiatterer_tf_red}%
  \filemodCmp{graphics/radiatterer_tf_red.tex}{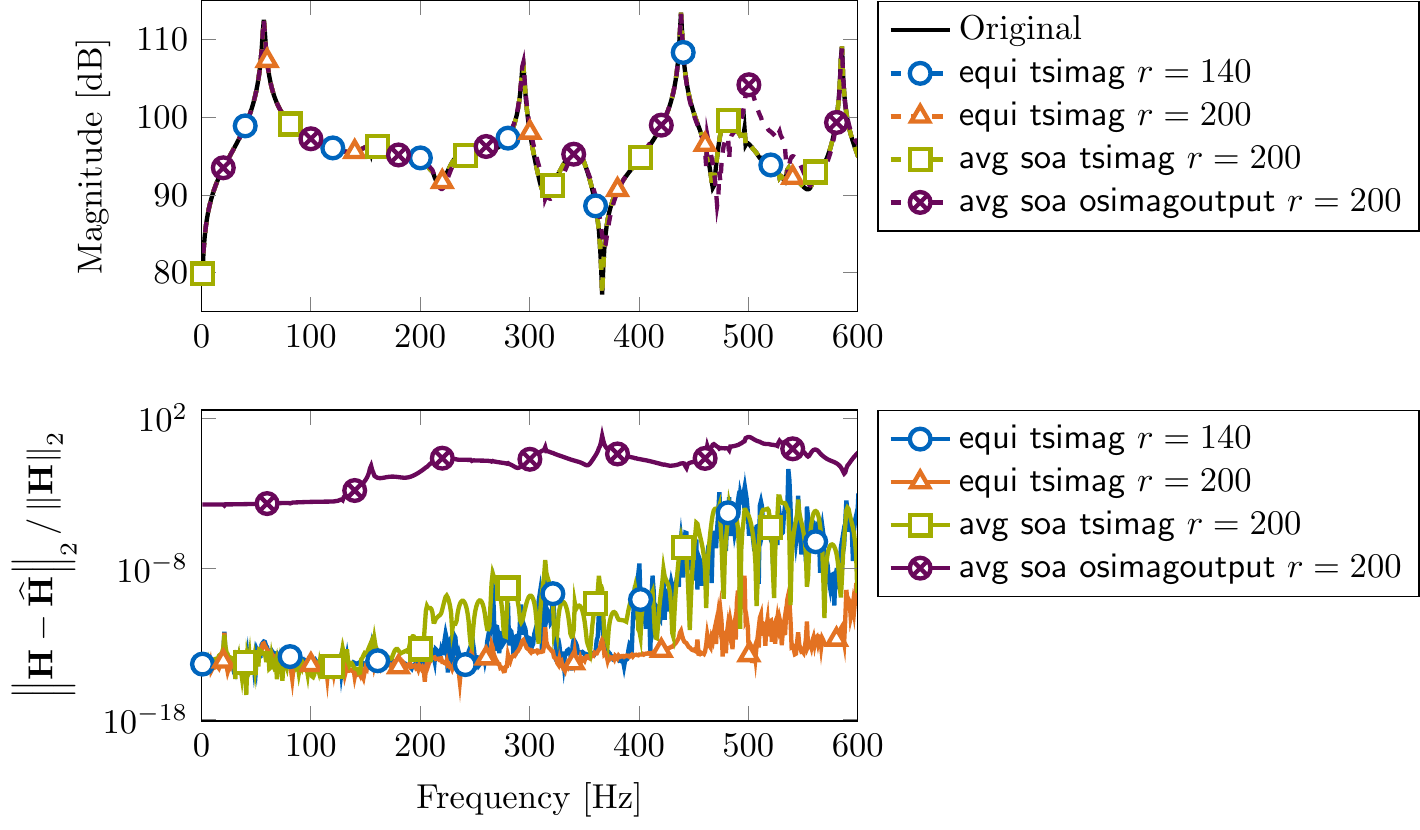}%
  {\tikzset{external/remake next}}{}%
  \begin{tikzpicture}[baseline=(current bounding box.north), font = \small]
	\begin{axis}[name=radiatterer_tf,
		height=0.2\textheight, width=.55\textwidth, 
		ylabel = {Magnitude \ubraces{\si{\decibel}}},
		xmin=0, xmax=600,
		xtick distance = 100,
		ymin=75, ymax=115,
		axis on top,
		legend cell align = left,
		legend pos = outer north east,
		cycle list name=TUMcolorlist,
		mark repeat={80}
		]
		\addplot[color=black, line width=\mylinewidth] table[x=s, y expr={10*log10(\thisrow{res}/1e-9)}, col sep=comma] {graphics/data/radiatterer_new_fine_res_tf_strint_avg_aaaa_tsimag.csv};
		\addlegendentry{Original}
		\pgfplotsset{cycle list shift=-1}
		
		\addplot+[dashed, mark phase={40}] table[x=s, y expr={10*log10(\thisrow{res_r}/1e-9)}, col sep=comma] {graphics/data/radiatterer_new_fine_res_tf_strint_equi_tsimag_r_140.csv};
		\addlegendentry{\lmorequi{} \ltsimag{} $r=140$}
		
		\addplot+[dashed, mark phase={60}] table[x=s, y expr={10*log10(\thisrow{res_r}/1e-9)}, col sep=comma] {graphics/data/radiatterer_new_fine_res_tf_strint_equi_tsimag.csv};
		\addlegendentry{\lmorequi{} \ltsimag{} $r=200$}
	
		\addplot+[dashed] table[x=s, y expr={10*log10(\thisrow{res_r}/1e-9)}, col sep=comma] {graphics/data/radiatterer_new_fine_res_tf_strint_avg_aaaa_tsimag.csv};
		\addlegendentry{\lmoravg{} \lmoraaa{} \ltsimag{} $r=200$}
	
		\addplot+[dashed, mark phase={20}] table[x=s, y expr={10*log10(\thisrow{res_r}/1e-9)}, col sep=comma] {graphics/data/radiatterer_new_fine_res_tf_strint_avg_aaaa_osimagoutput.csv};
		\addlegendentry{\lmoravg{} \lmoraaa{} \losimagoutput{} $r=200$}

	\end{axis}
	\begin{axis}[at={($(radiatterer_tf.south)-(0,1cm,)$)},anchor=north,
		height=0.2\textheight, width=.55\textwidth, ymode=log,
		xlabel = {Frequency \ubraces{\si{\hertz}}},
		ylabel = {$\twonorm{\tf - \tfr} / \twonorm {\tf}$},
		xmin=0, xmax=600,
		xtick distance = 100,
		axis on top,
		legend cell align = left,
		legend pos = outer north east,
		cycle list name=TUMcolorlist,
		mark repeat={80}
		]
		\addplot+ table[x=s, y=hinferr, col sep=comma] {graphics/data/radiatterer_new_fine_res_tf_strint_equi_tsimag_r_140.csv};
		\addlegendentry{\lmorequi{} \ltsimag{} $r=140$}
		
		\addplot+[mark phase={20}] table[x=s, y=hinferr, col sep=comma] {graphics/data/radiatterer_new_fine_res_tf_strint_equi_tsimag.csv};
		\addlegendentry{\lmorequi{} \ltsimag{} $r=200$}
		
		\addplot+[mark phase={40}] table[x=s, y=hinferr, col sep=comma] {graphics/data/radiatterer_new_fine_res_tf_strint_avg_aaaa_tsimag.csv};
		\addlegendentry{\lmoravg{} \lmoraaa{} \ltsimag{} $r=200$}

		\addplot+[mark phase={60}] table[x=s, y=hinferr, col sep=comma] {graphics/data/radiatterer_new_fine_res_tf_strint_avg_aaaa_osimagoutput.csv};
		\addlegendentry{\lmoravg{} \lmoraaa{} \losimagoutput{} $r=200$}

	\end{axis}
\end{tikzpicture}%
  \tikzexternaldisable%

	\caption{Original and reduced transfer functions as well as relative errors 
		for the radiation problem reduced with different methods.}
	\label{fig:radiatterer_tf_red}
\end{figure}
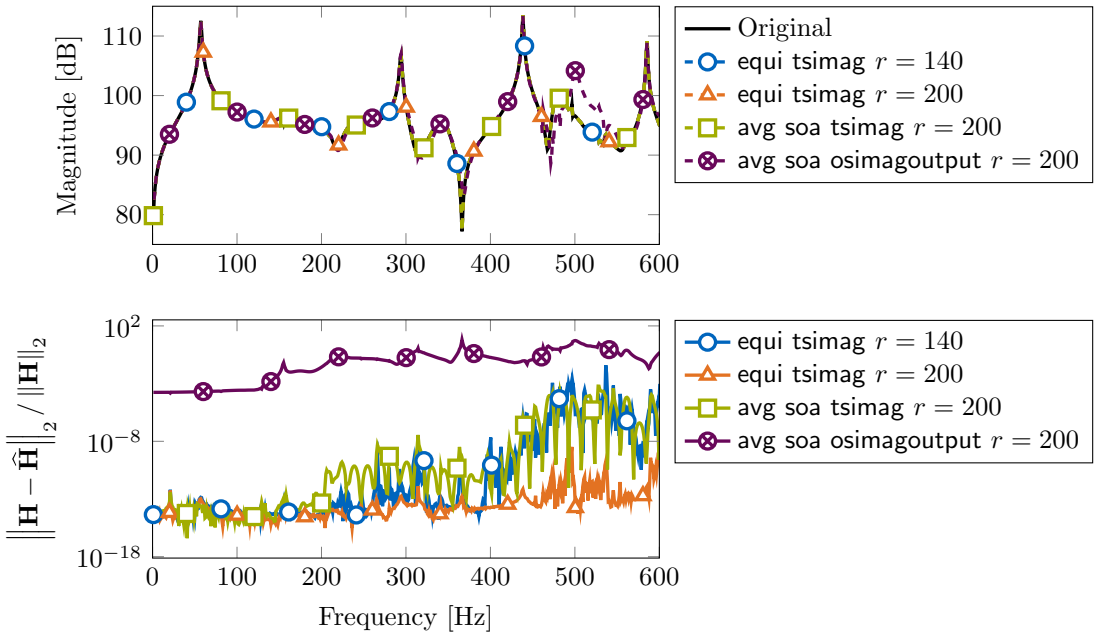

These observations are supported by \Cref{fig:radiatterer_tf_red} plotting the
relative errors of reduced-order models computed by \morequi{} with orders
$r = 140$ and $r = 200$.
While the larger reduced-order model shows a very small error over the 
complete frequency range, the smaller does not for the frequency region above 
$\SI{450}{\Hz}$.
It is, however, also apparent, that the approximation quality for all methods is 
better in the lower frequency range, presumably because of a large number of 
modes in the region above $\SI{450}{\Hz}$.
If this is known a priori, the locations of the expansion points can be altered
appropriately.
If this is not possible, the presampling involved in the methods \moravg{},
\morlinf{} and \morminrel{} shows its benefit.
At the cost of computing a larger intermediate reduction basis, the most
relevant information from this basis is chosen, allowing smaller reduced models
with better accuracy.
Choosing standard presampling or \morsp{} presampling yields accurate
reduced-order models with acceptable high \morscores.

If \moraaa{} presampling is employed, the resulting reduced-order models are 
less accurate than using the other presampling methods. 
For some projections, the reduced-order models do not accurately approximate the
original transfer function, cf. \Cref{fig:radiatterer_hinf}.
The error graph for \moravg{}~\moraaa{}~\tsimag{} in
\Cref{fig:radiatterer_tf_red} shows characteristic spikes at the locations of
the expansion points in the presampling basis.
This suggests that the employed second-order Krylov subspace does not contain
enough information to enable an as accurate approximation as the other
presampling methods.
A remedy would be to increase the size of the Krylov space, which would in turn
increase the size of the presampling basis. 
Note that increasing the size of the Krylov space is up to a certain degree 
less computationally expensive than establishing a completely new shift. 
Projection regarding the system output using a \moraaa{} presampling, however, 
does not yield a good approximation of the original system at all.

\subsection{Acoustic cavity with poroelastic layer}

\begin{figure}[t]
  \centering
  \begin{subfigure}[b]{.3\textwidth}
    \centering
    \vspace{.3cm}
	  \def\svgwidth{\textwidth}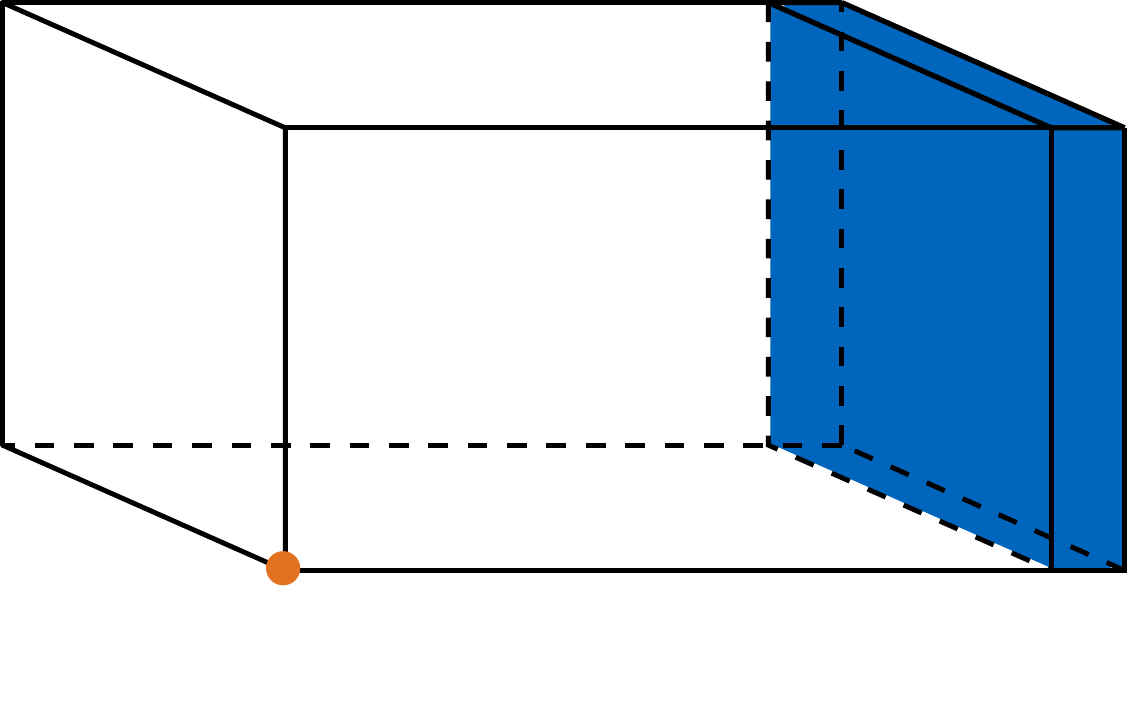
	  \subcaption{Sketch of the poroacoustic system with acoustic point source 
	  	$q$. The porous material covers the right wall.}
	  \label{fig:rumpler_sketch}
  \end{subfigure}%
  \hfill%
  \begin{subfigure}[b]{.65\textwidth}
    \centering
  \tikzexternalenable%
  \tikzsetnextfilename{rumpler2014}%
  \filemodCmp{graphics/rumpler2014.tex}{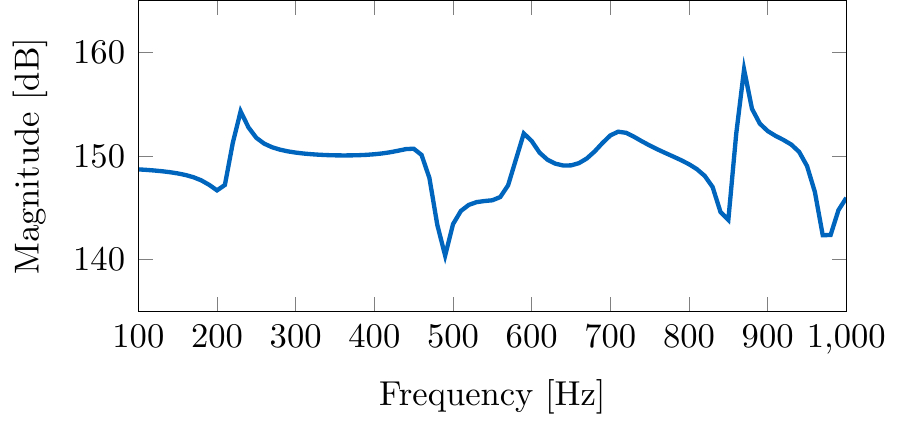}%
  {\tikzset{external/remake next}}{}%
  \begin{tikzpicture}[font = \small]
\begin{axis}[height=0.2\textheight, width=.9\textwidth, 
	xlabel = {Frequency \ubraces{\si{\hertz}}},
	ylabel = {Magnitude \ubraces{\si{\decibel}}},
	xmin=100, xmax=1000,
	xtick distance = 100,
	ytick distance = 10,
	ymin=135, ymax=165,
	axis on top,
	legend pos = outer north east, legend cell align = left,
	no markers
	]
	\addplot[mark=none, color=TUMBlau, line width=\mylinewidth] table[x index={0}, y index={1}, col sep=comma] {graphics/data/rumpler2014_fine_db.csv};
\end{axis}
\end{tikzpicture}%
  \tikzexternaldisable%

    \subcaption{Mean sound pressure level in the acoustic domain.}
    \label{fig:rumpler_tf}
  \end{subfigure}

  \caption{Sketch and transfer function of the poroacoustic model. }
\end{figure}

An acoustic cavity with dimensions $0.75\times 0.6\times \SI{0.4}{\meter}$ is 
examined in the following example.
One wall is covered by a \SI{0.05}{\meter} thick poroelastic layer acting as a
sound absorber.
The poroelastic material is described by the Biot theory~\cite{Bio56}.
The system is excited by an acoustic point source located in the corner opposite
to the porous material.
A sketch of the system is given in \autoref{fig:rumpler_sketch}.
The geometry and material parameters are taken from~\cite{RumGD14}, and the
discretized finite element model has the order $n=\num{386076}$.
We evaluate the model in the frequency range \SIrange{100}{1000}{\hertz}.
The material's frequency dependent dissipation mechanism and the coupling
between solid and fluid phase inside the material are modeled with in total six
complex-valued functions.
Due to the acoustic source, the transfer function also has a frequency-dependent
input vector. 
Thus, the system can be described by a Case~C transfer function with 
non-symmetric and complex-valued system matrices.
The transfer function measures the sound pressure level averaged over the
acoustic domain and is given in \Cref{fig:rumpler_tf}.

The model is reduced using all methods except \morbt{}, which is not applicable 
to this system because of the transfer function structure.
The standard presampling for \morminrel, \moravg and \morlinf{}
considers $n_{\rm{s}} = 200$ frequency shifts distributed linearly in
$2\pi\i [100, 1000]$.
The analytic derivatives of the frequency-dependent functions vanish for orders
larger than $6$ in the considered frequency range such that \morsp{} yields 7
columns for each shift.
Using $n_{\rm{s}} = 29$ shifts linearly distributed in the same range yields the
corresponding intermediate reduction basis with $q = 203$.
For \moraaa{} presampling, a local order of $10$ is chosen for each of
the $n_{\rm{s}} = 20$ shifts, which are also linearly distributed in 
$2\pi\i [100, 1000]$.
This results in an intermediate reduction basis of order $q = 200$.

\begin{figure}[t]
	\centering
  \tikzexternalenable%
  \tikzsetnextfilename{rumpler_morscore}%
  \filemodCmp{graphics/rumpler_morscore.tex}{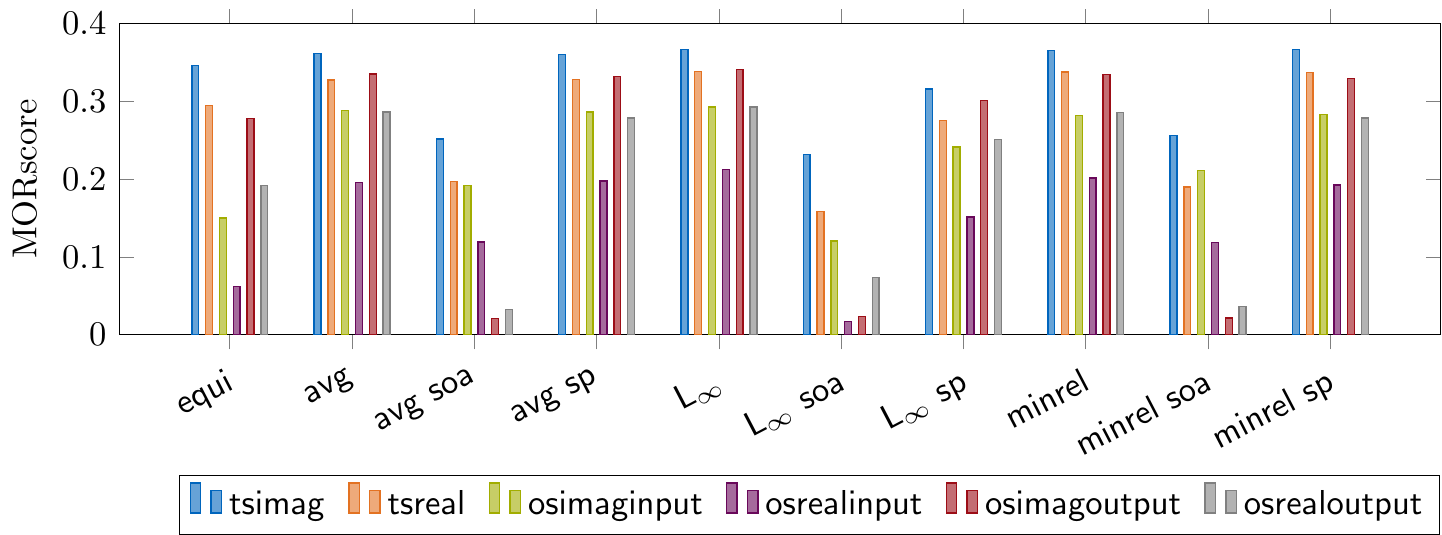}%
  {\tikzset{external/remake next}}{}%
  \begin{tikzpicture}[font = \small]
	\begin{axis}[height=0.2\textheight, width=\textwidth,
		ybar,
		ylabel = {\morscore{}},
		ymin=0, ymax=.4,
		legend cell align = left,
		legend style={at={(1,-0.45)},anchor=north east, /tikz/every even column/.append style={column sep=1ex}},
		legend columns=-1,
		no markers,
		every axis plot/.append style={fill},
		cycle list name=TUMcolorlist_simple,
		bar width=0.07cm,
		symbolic x coords={\lmorequi{}, \lmoravg{}, \lmoravg{} \lmoraaa{}, \lmoravg{} \lmorsp{}, \lmorlinf{}, \lmorlinf{} \lmoraaa{}, \lmorlinf{} \lmorsp{}, \lmorminrel{}, \lmorminrel{} \lmoraaa{}, \lmorminrel{} \lmorsp{}},
		xticklabel style={rotate=27.5,anchor=north east},
		]
		\addplot+[fill opacity=0.6] table [x=method, y index=1, col sep=comma] {graphics/data/rumpler2014_fine_mean_mor_score_tsimag.csv};
		\addlegendentry{\ltsimag{}}
		\addplot+[fill opacity=0.6] table [x=method, y index=1, col sep=comma] {graphics/data/rumpler2014_fine_mean_mor_score_tsreal.csv};
		\addlegendentry{\ltsreal{}}
		\addplot+[fill opacity=0.6] table [x=method, y index=1, col sep=comma] {graphics/data/rumpler2014_fine_mean_mor_score_osimaginput.csv};
		\addlegendentry{\losimaginput{}}
		\addplot+[fill opacity=0.6] table [x=method, y index=1, col sep=comma] {graphics/data/rumpler2014_fine_mean_mor_score_osrealinput.csv};
		\addlegendentry{\losrealinput{}}
		\addplot+[fill opacity=0.6] table [x=method, y index=1, col sep=comma] {graphics/data/rumpler2014_fine_mean_mor_score_osimagoutput.csv};
		\addlegendentry{\losimagoutput{}}
		\addplot+[fill opacity=0.6] table [x=method, y index=1, col sep=comma] {graphics/data/rumpler2014_fine_mean_mor_score_osrealoutput.csv};
		\addlegendentry{\losrealoutput{}}

	\end{axis}
\end{tikzpicture}%
  \tikzexternaldisable%

	\caption{\morscores{} of all employed reduction methods for 
		the poroacoustic system, with the maximum accuracy
		$\epsilon=\num{1e-16}$ and maximum order $r_{\rm{max}} = 100$.}
	\label{fig:rumpler_morscore}
\end{figure}
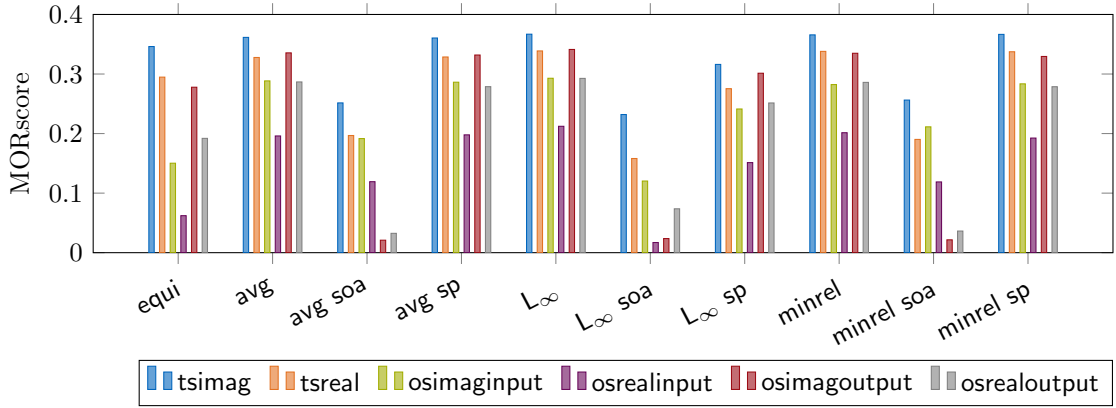

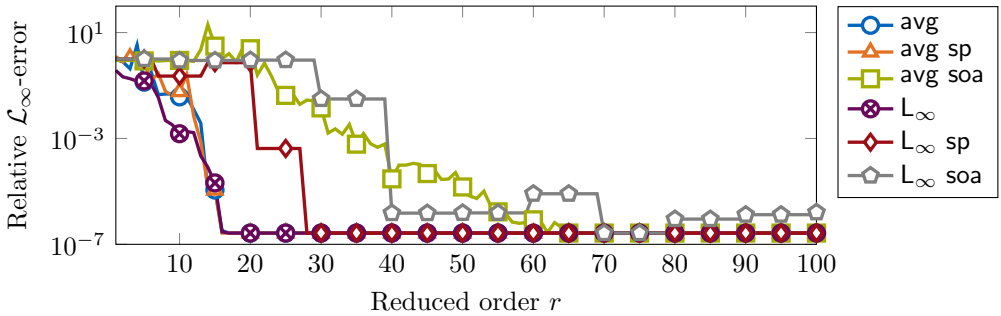
\begin{figure}[t]
	\centering
  \tikzexternalenable%
  \tikzsetnextfilename{rumpler_err_hinf}%
  \filemodCmp{graphics/rumpler_err_hinf.tex}{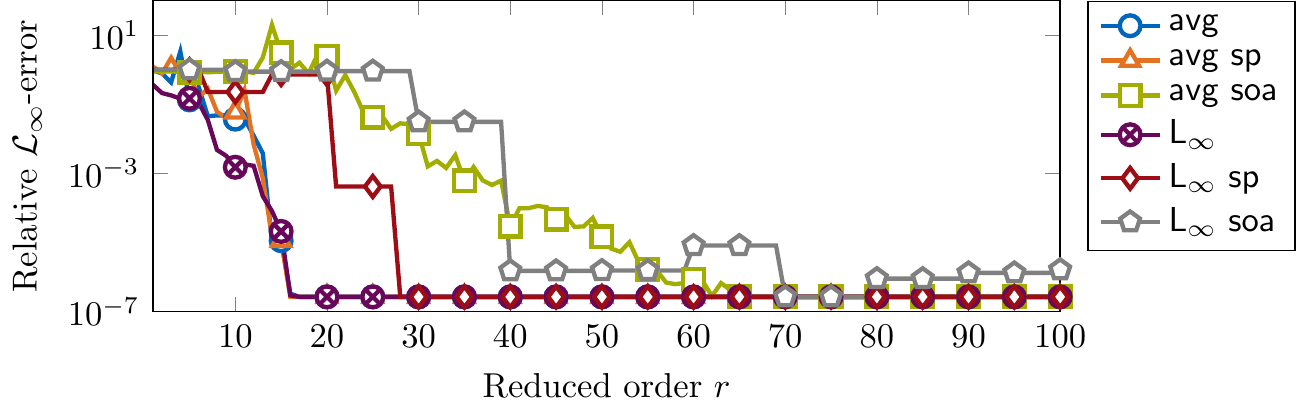}%
  {\tikzset{external/remake next}}{}%
  \begin{tikzpicture}[font = \small]
\begin{axis}[height=0.2\textheight, width=.72\textwidth, ymode=log,
	xlabel = {Reduced order $r$},
	ylabel = {\errlabel},
	xmin=1, xmax=100,
	ymin=1e-7, ymax=100,
	axis on top,
	legend pos = outer north east, legend cell align = left,
	cycle list name=TUMcolorlist,
	mark repeat={5},
	mark phase=5,
	]
	\addplot+ table[x=r, y=strint_avg, col sep=comma] {graphics/data/rumpler2014_fine_mean_hinfrelerr_tsimag.csv};
	\addlegendentry{\lmoravg}

	\addplot+ table[x=r, y=strint_avg_strprs, col sep=comma] {graphics/data/rumpler2014_fine_mean_hinfrelerr_tsimag.csv};
	\addlegendentry{\lmoravg{} \lmorsp{}}

	\addplot+ table[x=r, y=strint_avg_aaaa, col sep=comma] {graphics/data/rumpler2014_fine_mean_hinfrelerr_tsimag.csv};
	\addlegendentry{\lmoravg{} \lmoraaa{}}

	\addplot+ table[x=r, y=strint_linf, col sep=comma] {graphics/data/rumpler2014_fine_mean_hinfrelerr_tsimag.csv};
	\addlegendentry{\lmorlinf{}}
	\addplot+ table[x=r, y=strint_linf_strprs, col sep=comma] {graphics/data/rumpler2014_fine_mean_hinfrelerr_tsimag.csv};
	\addlegendentry{\lmorlinf{} \lmorsp{}}
	\addplot+ table[x=r, y=strint_linf_aaaa, col sep=comma] {graphics/data/rumpler2014_fine_mean_hinfrelerr_tsimag.csv};
	\addlegendentry{\lmorlinf{} \lmoraaa{}}
\end{axis}
\end{tikzpicture}%
  \tikzexternaldisable%

	\caption{Relative \Linf-errors of reduced models of the poroacoustic 
		system computed by several reduction methods. A two-sided projection 
		with complex-valued bases is considered for all cases.}
	\label{fig:rumpler_hinf}
\end{figure}

The \morscores{} of employed methods are reported in \Cref{fig:rumpler_morscore}
and show a good performance of nearly all reduction methods.
Even \morlinf{}~\morsp{}, whose reduced-order models are incremented in steps of
$r = 7$ has a comparable \morscore{}.
Additionally, it reaches an error as low as the other good performing methods
around $r = 28$, cf. \Cref{fig:rumpler_hinf}.
As already observed, the projections considering only the system output yield 
worse results if used in combination with \moraaa{} presampling.
To use \moraaa{} for this experiment, the approximations of the nonlinear
frequency-dependent functions are truncated after the quadratic term so that the
second-order Krylov subspace can be used. 
This has an impact on the approximation quality of the reduced-order models and 
results in a slower convergence of the approximation errors compared to the
other presampling strategies.
However, a comparably small error can also be achieved with \moraaa{}, if $r$ is
chosen high enough.
\Cref{fig:rumpler_hinf} compares the approximated relative $\Linf$-errors over
the reduced order for the different presampling methods.
The good performance of \morsp{} is evident here.
It is also interesting to note, that an \osimagoutput{} projection yields 
models, which accuracy is comparable to \tsreal{} for all methods except 
\moraaa{}.

\section{Conclusions}
\label{sec:conclusions}

In this work, we described structure-preserving model order reduction 
methods based on rational interpolation and balanced truncation applied 
to models of vibro-acoustic systems.
The benchmark examples were chosen such that their transfer functions exhibit
different properties, for example, complex-valued and/or frequency-dependent
system matrices or a frequency-dependent excitation.
Each benchmark case represents a relevant class of vibro-acoustic problems.
We also presented a strategy to incorporate higher-order frequency-dependent
terms in a standard second-order reduction framework.

The interpolation-based methods have been applicable to all considered models
and have been able to compute reduced-order models of reasonable accuracy and
small size. 
Second-order balanced truncation also has succeeded in computing compact
reduced models.
However, it is not applicable to systems with non-standard second-order 
transfer functions, which strongly restricts its application to vibro-acoustic
systems.
The methods based on oversampling the frequency response and extracting the most
relevant information have been shown to be the most successful. 
Strategies to leverage the initial cost of computing the presampling data have
been proposed and showed in many cases comparable results.


\addcontentsline{toc}{section}{Acknowledgments}
\section*{Acknowledgments}

The authors gratefully acknowledge the computational and data resources 
provided by the Leibniz Supercomputing Centre (\url{www.lrz.de}).

Parts of this work were carried out while Aumann was at the Technical 
University of Munich, Germany, and Werner was at the Max Planck Institute
for Dynamics of Complex Technical Systems in Magdeburg, Germany.
	
This research did not receive any specific grant from funding agencies in the 
public, commercial, or not-for-profit sectors.


\addcontentsline{toc}{section}{References}
\bibliographystyle{plainurl}
\bibliography{bibtex/myref}

\end{document}